\documentclass[12pt]{report}
\usepackage[colorlinks,linkcolor=linkcol,pdfpagelabels,hyperfootnotes=false]{hyperref}
\usepackage[dvipsnames]{color}
\definecolor{linkblue}{named}{MidnightBlue}
\hypersetup{linktocpage=false}

\usepackage{wsuthesis}

\usepackage{amssymb,amsthm,graphicx,multicol,amsmath,amsfonts}
\usepackage{cite}
\usepackage{graphicx}
\usepackage{wrapfig}
\usepackage{mathrsfs}
\usepackage{upgreek}
\usepackage[inline]{enumitem}
\usepackage{booktabs}
\usepackage[nameinlink,noabbrev,capitalize]{cleveref}
\usepackage[all]{hypcap}
\usepackage{subfigure}
\usepackage{commath}
\usepackage{tikz}
\usepackage[section]{placeins}
\usepackage[doublespacing,nodisplayskipstretch]{setspace}
\usepackage{import}
\usepackage[stable]{footmisc}

\usetikzlibrary{calc}
\usetikzlibrary{patterns}
\usetikzlibrary{positioning}
\usetikzlibrary{matrix}
\crefalias{subequation}{equation}

\definecolor{linkcol}{rgb}{0,0,0.5}

\newtheorem{theorem}{Theorem}[section]
\newtheorem{thm}{Theorem}[section]
\newtheorem{lemma}[theorem]{Lemma}
\newtheorem{lem}[theorem]{Lemma}
\newtheorem{corollary}[theorem]{Corollary}
\newtheorem{cor}[theorem]{Corollary}
\newtheorem{conjecture}[theorem]{Conjecture}
\newtheorem{subclaim}{Claim}[theorem]
\newtheorem{clm}{Claim}
\newtheorem{definition}[theorem]{Definition}
\newtheorem{dfn}[theorem]{Definition}
\newtheorem{remark}[theorem]{Remark}
\newtheorem{rem}[theorem]{Remark}

\newcommand{\D}{\operatorname{D}}
\newcommand{\spn}{\operatorname{span}}
\newcommand{\spt}{\operatorname{spt}}
\newcommand{\perimeter}{\operatorname{perimeter}}
\newcommand{\per}{\operatorname{P}}
\newcommand{\diam}{\operatorname{diameter}}
\newcommand{\inrad}{\operatorname{inradius}}
\newcommand{\intrr}{\operatorname{Int}}
\newcommand{\vol}{\operatorname{V}}
\newcommand{\mass}{\operatorname{M}}
\newcommand{\skel}[1]{\operatorname{skel}_{#1}}
\newcommand{\current}[1]{\mathcal{D}_{#1}}
\newcommand{\ccurrent}[1]{\mathcal{E}_{#1}}
\newcommand{\form}[1]{\mathcal{D}^{#1}}

\DeclareMathOperator{\transverse}{\cap\kern-7.75pt\top}
\DeclareMathOperator{\area}{Area}
\DeclareMathOperator{\argmax}{arg max}
\DeclareMathOperator{\lip}{Lip}

\newcommand{\boundary}{\partial}
\newcommand{\R}{ {\mathbb R} }
\newcommand{\Z}{ {\mathbb Z} }

\newcommand{\F}{ {\mathbb F} }
\newcommand{\va}{ \mathbf{a} }
\newcommand{\vb}{ \mathbf{b} }

\newcommand{\ve}{ \mathbf{e} }
\newcommand{\vs}{ \mathbf{s} }
\newcommand{\vt}{ \mathbf{t} }
\newcommand{\vu}{ \mathbf{u} }
\newcommand{\vv}{ \mathbf{v} }
\newcommand{\vw}{ \mathbf{w} }
\newcommand{\vx}{ \mathbf{x} }
\newcommand{\vy}{ \mathbf{y} }
\newcommand{\vz}{ \mathbf{z} }
\newcommand{\vzero}{ \mathbf{0} }
\newcommand{\vone}{ \mathbf{1} }
\newcommand{\st}{\mbox{s.t. }}
\newcommand{\dmsfn}{{\sc DMSFN}}

\input xy
\xyoption{all}

\begin{document}
\title{\uppercase{Data-inspired advances in geometric measure theory:\\generalized surface and shape metrics}}
\author{Sharif N. Ibrahim}
\chair{Kevin R. Vixie}
\committee{Bala Krishnamoorthy}{Thomas J. Asaki}{}
\degree{Doctor of Philosophy}
\dept{Department of Mathematics}
\thesistype{dissertation}
\submitdate{August}
\copyrightyear{2014}
\beforepreface
\indent Modern geometric measure theory, developed largely to solve the Plateau problem, has generated a great deal of technical machinery which is unfortunately regarded as inaccessible by outsiders.
Consequently, its ideas have not been incorporated into other fields as effectively as possible.
Some of these tools (e.g., distance and decompositions in generalized surface space using the flat norm) hold interest from a theoretical perspective but computational infeasibility prevented practical use.
Others, like nonasymptotic densities as shape signatures, have been developed independently as useful data analysis tools (e.g., the integral area invariant).
Here, geometric measure theory has promise to help close the gaps in our understanding of these ideas.

The flat norm measures distance between currents (or generalized surfaces) by decomposing them in a way that is robust to noise.
One new result here is that the flat norm can be suitably discretized and approximated on a simplicial complex by means of a simplicial deformation theorem.
While not surprising given the classical (cubical) deformation theorem or, indeed, Sullivan's convex cellular deformation theorem (which includes simplicial deformation as a special case), the bounds on the deformation can be made smaller and more practical by focusing on the simplicial case.

Computationally, the discretized flat norm can be expressed as a linear programming problem and thus solved in polynomial time.
Furthermore, the solution is guaranteed to be integral if the complex satisfies a simple topological condition (absence of relative torsion).
This discretized integrality result (with some work) yields a similar statement for the continuous case: the flat norm decomposition of an integral 1-current in the plane can be taken to be integral, something previously unknown for 1-currents which are not boundaries of 2-currents.

Nonasymptotic densities (integral area invariants) taken along the boundary of a shape are often enough to reconstruct the shape.
This result is easy when the densities are known for arbitrarily small radii but that is not generally possible in practice.
When only a single radius is used, variations on reconstruction results (modulo translation and rotation) of polygons and (a dense set of) smooth curves are presented.

\newpage
%You can put acknowledgments or whatever here
%\newpage
%\listoffigures
\afterpreface

\chapter{Introduction}
\label{ch:intro}
\newcommand{\measurerestr}{%
  \,\raisebox{-.127ex}{\reflectbox{\rotatebox[origin=br]{-90}{$\lnot$}}}\,%
}

\section{Overview}
This dissertation applies and extends geometric measure theory tools used for currents and densities.
In particular, the flat norm is used to measure currents and provides a useful metric in surface space.
This notion is discretized to obtain the multiscale simplicial flat norm and a simplicial deformation theorem (\cref{ch:msfn}, based on \cite{IbKrVi2013}) which approximates currents with chains on a simplicial complex via small deformations (as measured by the flat norm).

The multiscale simplicial flat norm can be computed efficiently and, for integral inputs, has guaranteed integral minimizers in several important cases (in particular, for codimension 1 chains).
This statement is stronger than what was known for the continuous case (where the statement was for codimension 1 boundaries).
Bridging the gap between these statements and extending the discrete results to the continuous case is the goal of \cref{ch:ic} (based on \cite{IbKrVi2014}) where it is shown for 1-currents in $\R^2$ with a framework for establishing the result in general assuming suitable triangulation results.

Lastly, the notion of nonasymptotic densities (also known as the integral area invariant) is developed in the plane in \cref{ch:nad} (based on \cite{ibrahim-2014-1}) where uniqueness questions are addressed in light of a certain useful regularity condition (tangent cone graph-like).

This research was supported in part by the National Science Foundation through grants DMS-0914809 and CCF-1064600.

%\section{Preface}
%TODO
\section{Measure theory}
A few concepts from measure theory prove useful in our development.
The Hausdorff measure allows us to sensibly measure $m$-dimensional sets in $\R^n$.

\begin{definition}[Hausdorff measure]
Given a set $A \subset \R^n$, the $m$-dimensional Hausdorff measure of $A$ is an outer measure defined by
\[
\mathcal{H}^m(A) = \lim_{\delta \downarrow 0} \del{\inf_{\mathcal{S}} \sum_{S_j \in \mathcal{S}} \alpha_m \del{\frac{\diam{S_j}}{2}}^m}
\]
where $\alpha_m$ is the volume of the unit ball in $\R^m$ and the infimum is taken over all countable coverings $\mathcal{S} = \{S_1, S_2, \dots\}$ of $A$ with every $S_j \in \mathcal{S}$ having diameter at most $\delta$.
\end{definition}

The Hausdorff measure approximates $A$ locally by covering it with small sets which in turn have their $m$-dimensional volumes approximated by balls of the same radius in $\R^m$.
This is the natural way to measure $m$-dimensional volume in $\R^n$ and agrees with intuitive notions of what this should mean, for example, for an $m$-dimensional manifold embedded in $\R^n$.
It also provides sensible results for any nonnegative real dimension by extending the unit ball volume via the $\Gamma$ function: $\alpha_m = \pi^{m/2}/\Gamma(m/2+1)$.
For any particular nonempty set $A$, there is a ``correct'' dimension to use when measuring it with the Hausdorff measure in the sense that using any other value yields a trivial result.

\begin{definition}[Hausdorff dimension]
The Hausdorff dimension of a nonempty set $A \subseteq \R^n$ is the unique nonnegative real number $m$ such that $\mathcal{H}^p(A) = 0$ for all $p > m$ and $\mathcal{H}^q(A) = \infty$ whenever $m > q$ and $q \geq 0$.
\end{definition}

Knowing that the set $A$ has Hausdorff dimension $m$ places no restrictions on $\mathcal{H}^m(A)$.
That is, one can construct examples with any desired measure in the interval $[0,\infty]$.

\begin{definition}[Rectifiable sets]
\label{def:rectifiableset}
A set $A \subseteq \R^n$ is called an $m$-dimensional rectifiable set if $\mathcal{H}^m(A) < \infty$ and there exists a set $E$ such that $\mathcal{H}^m(A-E) = 0$ and $E$ is the union of the images of countably many Lipschitz functions from $\R^m$ to $\R^n$.
\end{definition}

\begin{definition}[Density]
\label{def:densityset}
Given a set $A \subseteq \R^n$ and $1 \leq m \leq n$, the $m$-dimensional density of $A$ at a point $x \in \R^n$ is given by
\[
\upvartheta^m(A, x) = \lim_{r \downarrow 0} \frac{\mathcal{H}^m(A \cap B(x,r))}{\alpha_m r^m}
\]
where $B(x,r)$ is the closed ball in $\R^n$ with center $x$ and radius $r$ and $\alpha_m$ is the volume of the unit ball in $\R^m$.
\end{definition}

\begin{definition}[Density of measures]
Given a measure $\mu$ on $\R^n$, $1 \leq m \leq n$, and $x \in \R^n$, we define the $m$-dimensional measure of $\mu$ at $x$ by
\[
\uptheta^m(\mu, x) = \lim_{r \downarrow 0} \frac{\mu(B(x,r)}{\alpha_m r^m}.
\]
\end{definition}

Density of a set in \cref{def:densityset} is a special case of density of measures using the Hausdorff measure restricted to $A$ (denoted $\mathcal{H}^m \measurerestr A$ and defined by $(\mathcal{H}^m \measurerestr A)(B) = \mathcal{H}^m(A \cap B)$).

\section{Currents}
The following is a brief introduction to currents, largely following Federer\cite{Federer1969}, Krantz and Parks\cite{KrantzParks2008}, and Morgan\cite{Morgan2008} which are recommended as references for some of the details in descending order of difficulty.
Currents are the primary objects of study in \cref{ch:msfn,ch:ic} where the definition of various types of currents (general, normal, and integral) and the flat norm on currents play a central role.
There is significant machinery to develop which can obscure the intuition which is that (suitably nice) $m$-currents can be thought of like oriented submanifolds of dimension $m$.

\begin{definition}[$m$-covectors]
\label{def:mcovec}
Given $n$ and $m$, the set of $m$-covectors of $\R^n$ is denoted by $\wedge^m(\R^n)$ and contains all $\phi$ such that
\begin{itemize}
\item $\phi$ maps a collection of $m$ vectors in $\R^n$ to a real number: $\phi : (\R^n)^m \rightarrow \R$.
\item $\phi$ is $m$-linear; that is, linear in each of its $m$ arguments.
In particular,
\begin{align*}
&\phi(\vu_1, \vu_2, \dots, \vu_{\ell-1}, \alpha \vv + \beta \vw, \vu_{\ell+1}, \dots, \vu_m)\\
&=\alpha\phi(\vu_1, \dots, \vu_{\ell-1}, \vv, \vu_{\ell+1}, \dots, \vu_m)\\
&\quad +\beta\phi(\vu_1, \dots, \vu_{\ell-1}, \vw, \vu_{\ell+1}, \dots, \vu_m)
\end{align*}
whenever $1 \leq \ell \leq m$, $\alpha, \beta \in \R$, and $\vv, \vw, \vu_i \in \R^n$.
\item $\phi$ is alternating: transposing any two arguments changes the sign.
If $1 \leq i < j \leq m$ and $\vu_1, \dots, \vu_m \in \R^n$, then we have
\begin{align*}
&\phi(\vu_1, \dots, \vu_{i-1}, \vu_i, \vu_{i+1}, \dots, \vu_{j-1}, \vu_j, \vu_{j+1}, \dots \vu_m)\\
&=-\phi(\vu_1, \dots, \vu_{i-1}, \vu_j, \vu_{i+1}, \dots, \vu_{j-1}, \vu_i, \vu_{j+1}, \dots \vu_m).
\end{align*}
\end{itemize}
\end{definition}

The most well-known function with these properties is the determinant applied to $n$-by-$n$ matrices.
It is easy to show that the determinant is (up to multiplication) the only member of 
$\wedge^n(\R^n)$.

Given the standard basis vectors $\ve_i$ for $\R^n$, we define dual basis vectors $\dif x_j$ linearly by
\[
\dif x_j ( \ve_i) = \begin{cases}
1 \text{ if } i = j,\\
0 \text{ if } i \neq j
\end{cases}
\]
and note that any 1-covector can be written in this basis.

\begin{definition}[Exterior product, simple covectors]
Given $a_1, \dots, a_m \in \wedge^1(\R^n)$, we denote the exterior or wedge product of these 1-covectors to be the $m$-covector
\[
a_1 \wedge a_2 \wedge \dots \wedge a_m
\]
which is defined by
\[
(a_1 \wedge a_2 \wedge \dots \wedge a_m)(\vu_1, \dots, \vu_m) = {\det\del{\begin{pmatrix}\vline&&\vline\\
\va_1&\dots&\va_m\\
\vline&&\vline\end{pmatrix}^T\begin{pmatrix}\vline&&\vline\\
\vu_1&\dots&\vu_m\\
\vline&&\vline\end{pmatrix}}}
\]
where the vector $\va_i$ is the representation of the 1-covector $a_i$ in the dual basis $[\dif x_1, \dots, \dif x_n ]$.
Any element of $\wedge^m(\R^n)$ that can be written as a wedge product of 1-covectors $a_i$ is called simple and every $m$-covector can be expressed as the sum of simple $m$-covectors.
The wedge product extends to higher degree covectors by means of this decomposition and a distributive law.
\end{definition}

The wedge product is $m$-linear and is negated when any two covectors are transposed because it relies on the determinant.
For the same reason, if a particular 1-covector appears more than once in the wedge product, the result is 0.
Working an  example, we have
\begin{align*}
&\del{4\dif x_1 \wedge \dif x_3 + 3 \dif x_4 \wedge \dif x_3} \wedge \del{2\dif x_1 \wedge \dif x_2 - \dif x_1 \wedge \dif x_3 }\\
&= 8 \dif x_1 \wedge \dif x_3 \wedge \dif x_1 \wedge \dif x_2  - 4\dif x_1 \wedge \dif x_3 \wedge \dif x_1 \wedge \dif x_3\\
&\quad + 6\dif x_4 \wedge \dif x_3 \wedge \dif x_1 \wedge \dif x_2 - 3\dif x_4 \wedge \dif x_3 \wedge \dif x_1 \wedge \dif x_3\\
&= 6\dif x_4 \wedge \dif x_3 \wedge \dif x_1 \wedge \dif x_2.
\end{align*}

\begin{definition}[Forms]
\label{def:form}
Given open $U \subseteq \R^n$, a differential $m$-form on $U$ is a function $\phi : U \rightarrow \wedge^m(\R^n)$.
The set of all $m$-forms on $U$ is denoted by $\form{m}(U)$.
We say that $\phi \in \form{m}(U)$ is $C^k$ if $\phi(\vx)$ applied to $\vv_1\wedge \vv_2 \wedge \dots \wedge \vv_m$ is a $C^k$ function in $\vx$ for any fixed vectors $\vv_i \in \R^n$.
%TODO: Support?
%The support of $\phi$ is the closure of the set $\{x \in U \mid \phi(x) \neq 0\}$.
%Given $m < n$, the set of $m$-forms in $\R^n$ is denoted by $\form{m}$.
%Each element $\phi \in \form{m}$ associates with each point $x \in \R^n$ a multilinear, alternating map $\phi(x) : \wedge^m \R^n \rightarrow \R$.
%The support of $\phi$ 
%supported on the set $U \subseteq \R^n$ is denoted by $\form{m}(U)$.
\end{definition}

Observe that any function $f : U \rightarrow \R$ can be considered as a 0-form.
Differential $m$-forms can be used as integrands over $m$-dimensional surfaces as they can vary both based on location of a point and its tangent plane; this serves as a useful generalization of integration of 1-forms over a curve.

\begin{definition}[Exterior differentiation]
Suppose $U \subset \R^n$ is open and $f : U \rightarrow \R$ is $C^1$.
The exterior derivative of the 0-form $f$ is the 1-form $\dif f$ defined by
\[
\dif f = \pd{f}{x_1}\dif x_1 + \pd{f}{x_2}\dif x_2 + \dots + \pd{f}{x_n}\dif x_n.
\]
The exterior derivative of the simple $m$-form $\phi = f \dif x_{i_1} \wedge \dif x_{i_2} \wedge \dots \wedge \dif x_{i_m}$ (where the $i_k$ are integers from 1 to $n$) is given by the $(m+1)$-form
\[
\dif \phi = \dif f \wedge \dif x_{i_1} \wedge \dif x_{i_2} \wedge \dots \wedge \dif x_{i_m}.
\]
For all other $C^1$ $m$-forms, the definition is extended by linearity.
\end{definition}

\begin{theorem}[Properties of exterior differentiation, \cite{KrantzParks2008} p. 163]
\label{thm:formprop}
If $\phi$ and $\psi$ are $C^1$ $m$-forms and $\theta$ is a $C_1$ $\ell$-form, then we have:
\begin{itemize}
\item $\dif(\phi + \psi) = \dif \phi + \dif \psi$
\item $\dif(\phi \wedge \theta) = (\dif \phi) \wedge \theta + (-1)^m\phi \wedge (\dif \theta)$
\item If $\phi$ is $C^2$, then $\dif \dif \phi = 0$.
\end{itemize}
\end{theorem} 

\begin{definition}[Currents]
\label{def:current}
The space of $m$-currents in $\R^n$ is denoted by $\current{m}(\R^n)$ and is dual to the set of $C^\infty$ differential forms of compact support.
The mass of a current $T \in \current{m}(\R^n)$ is given by 
\[
\mass(T) = \sup\{T(\phi) \mid \phi \in \form{m}(\R), \phi\text{ has compact support, } \norm{\phi} \leq 1\}.
\]
\end{definition}

Whereas differential forms correspond to integrands, (suitably nice) currents can be intuitively thought of as the linear integration operator itself, representing and generalizing the oriented submanifolds over which we can integrate differential forms.
When $T \in \current{m}(\R^n)$ represents an oriented submanifold, the mass is simply its $m$-dimensional volume, counting multiplicities (this is the intuition take take away from this definition).

\begin{definition}[Boundary of a current]
\label{def:boundcurrent}
The boundary of an $m$-current $T \in \current{m}(U)$ is defined in terms of the exterior derivative on differential forms.
Namely, for $m > 0$, we let $\boundary T \in \mathcal{D}_{m-1}(U)$ be the linear operator on $(m-1)$-forms defined by $\boundary T(\phi) = T(\dif \phi)$ for all $\phi \in \form{m-1}(U)$.
For $m = 0$, we let $\boundary T = 0$ as a 0-current.
\end{definition}

%For all currents $T$, $\boundary \boundary T = 0$ because $d d \phi = 0$.
This definition along with facts about exterior differentiation immediately provides us with some useful properties:
\begin{subequations}
\begin{align}
\label{eq:bbcurrentzero}
\boundary \boundary T &= 0,\\
\boundary(\alpha T_1 + \beta T_2) &= \alpha \boundary T_1 + \beta \boundary T_2.
\end{align}
\end{subequations}
Defining $\boundary T$ for 0-currents is not universal, but doing so allows us to simplify various statements slightly by omitting special cases (\cref{eq:bbcurrentzero}, for example).

\begin{definition}[Support of a current]
\label{def:support}
The support of a current $T$ in $\current{m}$ is the complement of the largest open set $U$ such that $T(\phi) = 0$ whenever $\phi \in \form{m}(U)$.
\end{definition}

Currents can be created from any oriented $m$-dimensional rectifiable set $R$.
Define $S$ as the map from points $x \in R$ to unit $m$-vectors corresponding to the oriented tangent plane to $R$ at $x$.
By this, we mean that $S(x)$ is the wedge product of $m$ orthonormal tangent vectors to $R$ at $x$.
% and the function and is made continuous by suitable choice of orientations at each point.
Then for any differential form $\phi \in \form{m}(\R^n)$, define an $m$-current $T$ by
\[
T(\phi) = \int_R \langle \phi(x), S(x)\rangle \dif \mathcal{H}^m.
\]
We allow $T$ to carry integer multiplicities by introducing a function $\eta : R \to \Z$ with $\int_R \eta(x) \dif \mathcal{H}^m < \infty$ to obtain
\[
T(\phi) = \int_R \langle \phi(x), S(x)\rangle \eta(x) \dif \mathcal{H}^m.
\]
The currents which can be constructed via this procedure are called rectifiable currents and their existence justifies the statement that currents generalize oriented submanifolds.

\begin{definition}[Rectifiable currents]
\label{def:rectifiablecurrent}
A rectifiable current is a current with compact support associated with a rectifiable set with integer multiplicities and finite total measure (counting multiplicities).
\end{definition}

\begin{definition}[Normal currents]
\label{def:intronormalcurrent}
An $m$-current $T$ is normal if and only if $\mass(T) + \mass(\boundary T) < \infty$ and the support of $T$ is compact.
\end{definition}

Note that nothing prevents normal currents from being ``smeared'' out in space.
Morgan\cite[p. 48]{Morgan2008} gives an example of a normal 1-current $S_2$ which covers the unit square in $\R^2$ but with the concentration of mass of the 2-dimensional Hausdorff measure so it has finite mass and boundary.

\begin{definition}[Integral currents]
\label{def:integralcurrent}
A current $T$ is an integral current if $T$ and $\boundary T$ are rectifiable.
\end{definition}

As an aid to understanding the various classes of currents, note that all integral currents are both normal and rectifiable (in fact, this can be taken as the definition of integral current by the closure theorem\cite[4.2.16]{Federer1969}).
Furthermore, integral and rectifiable currents have integer multiplicities while normal and general currents need not.

\begin{figure}[ht!]
\centering
\includegraphics[scale=0.7, trim=1.25in 8.3in 2.8in 1in, clip]{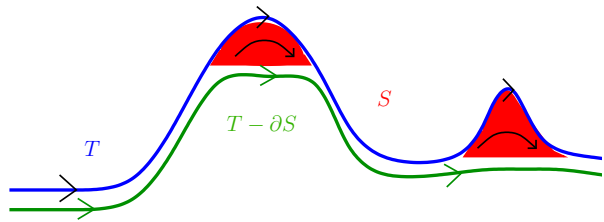}
\caption[Flat norm decomposition]{\label{fig:ic1dcurrent} The flat norm decomposes the 1D current $T$ into (the boundary of) a 2D piece $S$ and the 1D piece $X = T - \boundary S$. The resulting current is shown slightly separated from the input current for clearer visualization.}
\label{fig:introflatnorm}
\end{figure}

\begin{figure}
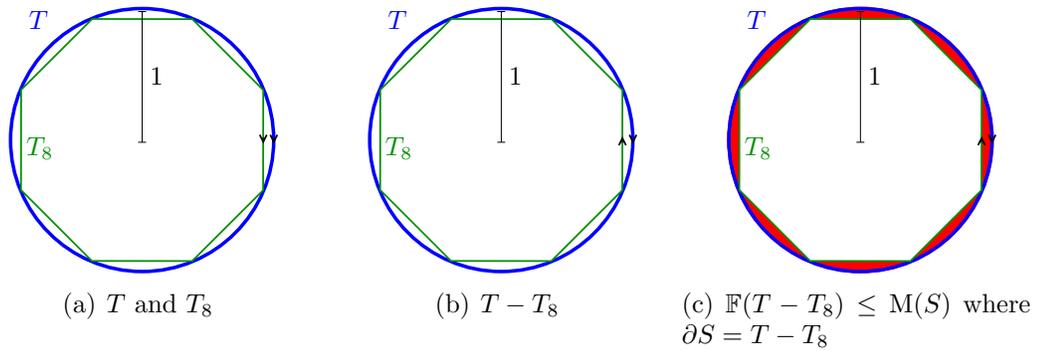

\begin{center}
\subfigure[$T$ and $T_8$]{\label{subfig:flatdistance1} \makebox[\textwidth/7*2][c]{\input{polygonexampleb1.pspdftex}}}
\subfigure[$T - T_8$]{\label{subfig:flatdistance2} \makebox[\textwidth/7*2][c]{\input{polygonexampleb2.pspdftex}}}
\subfigure[$\F(T-T_8) \leq \mass(S)$ where $\boundary S = T - T_8$]{\label{subfig:flatdistance3} \makebox[\textwidth/7*2][c]{\input{polygonexampleb3.pspdftex}}}
%\subfigure[]{\subimport{figs}{polygonexampleb1.pspdftex}}
%\subfigure[]{\input{figs/polygonexampleb2.pspdftex}}
%\subfigure[]{\input{figs/polygonexampleb3.pspdftex}}
\end{center}
\caption[Flat norm convergence]{The flat norm indicates the unit circle $T$ and inscribed $n$-gon $T_n$ are close because the region they bound has small area.}
\end{figure}

Suppose $T$ is the current representing the unit circle in $\mathbb{R}^2$ and $T_n$ is an inscribed regular $n$-gon, both oriented clockwise (see \cref{subfig:flatdistance1}).
As $n$ gets large, it is clear that $T_n$ intuitively approximates $T$ arbitrarily well.
Thus it would be desirable to have a notion of convergence for which $T_n \rightarrow T$.
In particular, the mass norm is not useful on its own here: the current $T_n - T$ has mass $\mass(T_n) + \mass(T) \rightarrow 4\pi$ since $T_n$ and $T$ do not exactly coincide (so there is no cancellation) except on a measure 0 subset.

\begin{definition}[Flat norm]
\label{def:flatnorm}
Given an $m$-current $T \in \current{m}(\R^n)$, we define its flat norm $\F(T)$ to be the least cost decomposition of $T$ into two pieces: the boundary of an $(m+1)$-current $S$ and the $m$-current $X = T - \boundary S$ (see \cref{fig:introflatnorm}).
%and $m$-current $X$ and the boundary of an $(m+1)$-current $S$.
The cost of a particular decomposition $T = X + \boundary S$ is given by $\mass(X) + \mass(S)$.
Formally,
\begin{align*}
\F(T) &= \min\{\mass(X) + \mass(S) \mid T = X + \boundary S, \, X \in \ccurrent{m},\, S \in \ccurrent{m+1}\}
\end{align*}
where $\ccurrent{m} \subset \current{m}$ is the set of $m$-currents with compact support.
\end{definition}

The flat norm is usually defined as a supremum over forms but this definition is equivalent and more immediately useful for our purposes.
Of note is that the minimum exists and is attained whenever $\F(T) < \infty$.
This is proved using the Hahn-Banach theorem\cite[p. 367]{Federer1969}.
If $\mass(T) < \infty$, then $\mass(S) + \mass(\boundary S) < \infty$ so $S$ is normal by \cref{def:intronormalcurrent}.

%TODO: lower semicontinuous with respect to mass

\chapter{Multiscale simplicial flat norm\footnote{Previously published as \cite{IbKrVi2013}}}
\label{ch:msfn}
%%%%%%%%%%%%%%%%%%%%%%%%%%%%%%%%%%%%%%%%%%%%%%
%% Simplicial Flat Norm with Scale - for arXiv
%%%%%%%%%%%%%%%%%%%%%%%%%%%%%%%%%%%%%%%%%%%%%%

% Format cells in comparison table
\newcommand{\twolinecelleq}[2]{
% Redefining the command to be just this will display the equation on a single line
%$#1#2$
%
% Avoid adding extra vertical space with the parbox 
\setlength{\abovedisplayskip}{0pt}
\setlength{\belowdisplayskip}{0pt}
\parbox{2.0in}{
\centering
\begin{multline*}
#1\\
#2
\end{multline*}
}
}

\newcommand{\MSFN}{multiscale simplicial flat norm}
\newcommand{\OHCP}{optimal homologous chain problem}
\newcommand{\OBCP}{optimal bounding chain problem}

%% Start of paper
%% --------------

\section{Introduction} \label{sec:introduction} 

{\em Currents} are standard objects studied in geometric measure
theory, and are named so by analogy with electrical currents that have
a kind of magnitude and direction at every point. Intuitively, one
could think of currents as generalized surfaces with orientations and
multiplicities. The mathematical machinery of currents has been used
to tackle many fundamental questions in geometric analysis, such as
the ones related to area minimizing surfaces, isoperimetric problems,
and soap-bubble conjectures \cite{Morgan2008}. 

To formally define $d$-currents in $\R^n$, we first let
$\mathcal{D}^d$ be the set of $C^\infty$ differentiable $d$-forms with
compact support.  Then the set of $d$-currents is given by the dual
space of $\mathcal{D}^d$ (denoted $\mathcal{D}_d$) with the weak
topology.  We denote by $\mathcal{R}_m$ the set of rectifiable
currents, which contains all currents that represent oriented
rectifiable sets (i.e., sets which are almost everywhere the countable
union of images of Lipschitz maps from $\R^m$ to $\R^n$) with integer
multiplicities and finite total mass (with multiplicities).

The {\em mass} $\mass(T)$ of a $d$-dimensional current $T$ can be
thought of intuitively as the weighted $d$-dimensional volume of the
generalized object represented by $T$. For instance, the mass of a
$2$-dimensional current can be taken as the area of the surface it
represents.  Formally, the mass of $T$ is given by $M(T) = \sup_{\phi
  \in \mathcal{D}^d} \{T(\phi) \mid \sup \| \phi(x) \| \leq 1\}.$

\begin{figure}[ht!]  
\centering
\includegraphics[scale=0.7, trim=1.25in 8.3in 2.8in 1in, clip]{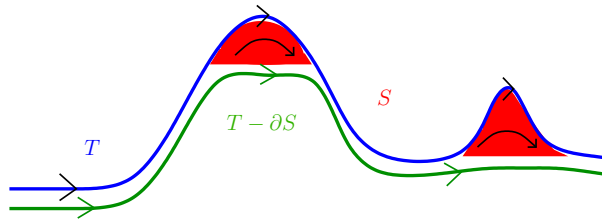}
\caption[Flat norm decomposition]{\label{fig:1dcurrent} 1D current $T$, and flat norm
decomposition $T-\boundary S$ at appropriate scale $\lambda$. The
resulting current is shown slightly separated from the input current
for clearer visualization.}
\end{figure}

The boundary $\boundary T$ of a current $T$ is defined by duality with
forms.  That is, we have $\boundary T(\phi) = T(d \phi)$ for every
differential form $\phi \in \mathcal{D}^d$.  Note that when $T$
represents a smooth oriented manifold with boundary, this corresponds
to the usual definition of boundary.  We restrict our attention to
integral currents $T$ that are rectifiable currents with a rectifiable
boundary (i.e., $T \in \mathcal{R}_m$ and $\boundary T \in
\mathcal{R}_{m-1}$).
%have the property that $T \in \mathcal{R}_m$ and $\boundary T \in \mathcal{R}_{m-1}$.
The {\em flat norm} of a $d$-dimensional current $T$ is given by
\begin{equation}
\label{eq:flatnorm}
\F(T) = \min_S \{ \mass(T-\boundary S) + \mass(S)~\big|~ T-\boundary S 
	\in {\mathscr E}_d, \, S \in {\mathscr E}_{d+1}\},
\end{equation}
where ${\mathscr E}_d$ is the set of $d$-dimensional currents with
compact support. One also uses flat norm to measure the ``distance''
between two $d$-currents. More precisely, the flat norm distance
between two $d$-currents $T$ and $P$ is given by
\begin{equation} \label{eq:flatdist}
\F(T,P) = \inf \{\mass(Q) + \mass(R) \, \big| \, T-P = Q+\boundary R,
\, Q \in {\mathscr E}_d, \, R \in {\mathscr E}_{d+1}\}.
\end{equation}
Morgan and Vixie \cite{MoVi2007} showed that the $L^1$ total variation
functional ($L^1$TV) introduced by Chan and Esedo\=glu \cite{ChEs2005}
computes the flat norm for boundaries $T$ with integer
multiplicity. Given this correspondence, and the use of {\em scale} in
$L^1$TV, Morgan and Vixie defined \cite{MoVi2007} the flat norm with
scale $\lambda \in [0,\infty)$ of an oriented $d$-dimensional set $T$
as
\begin{equation}
\label{eq:flatnormwscale}
\F_{\lambda}(T) \equiv \min_S \{ \vol_d(T-\boundary S) + \lambda \vol_{d+1}(S) \},
\end{equation}
where $S$ varies over oriented $(d+1)$-dimensional sets, and $\vol_d$
is the $d$-dimensional volume, used in place of mass. Figure
\ref{fig:1dcurrent} illustrates this definition. Flat norm of the 1D
current $T$ is given by the sum of the length of the resulting
oriented curve $T - \boundary S$ (shown separated from the input curve
for clarity) and the area of the 2D patch $S$ shown in red. Large
values of $\lambda$, above the curvature of both humps in the curve
$T$, preserve both humps. Values of $\lambda$ between the two
curvatures eliminate the hump on the right. Even smaller values
``smooth out'' both humps as illustrated here, giving a more ``flat''
curve, as $S$ can now be comprised of much bigger 2D patches.

%\begin{wrapfigure}{r}{3in} 
 %\end{wrapfigure}

Figure \ref{fig:1dcurrent} illustrates the utility of flat norm for
deblurring or smoothing applications, e.g., in 3D terrain maps or 3D
image denoising. But efficient methods for computing flat norm are
known only for certain types of currents in two dimensions. For $d=1$,
Under the setting where $T$ is a boundary, i.e., a loop, embedded in
$\R^2$ and the minimizing surface $S \in \R^2$ as well, the flat norm
could be calculated efficiently, for instance, using graph cut methods
\cite{KoZa2002} -- see the work of Goldfarb and Yin \cite{GoYi2009}
and Vixie et al.~\cite{Vietal2010}, and references therein. Motivated
by applications in image analysis, these approaches usually worked
with a grid representation of the underlying space ($\R^2$). Pixels in
the image readily provide such a representation.

While it is computationally convenient that $L^1$TV minimizers give us
the scaled flat norm for the input images, this approach restricts us
to currents that are boundaries of codimension 1.  Correspondingly,
the calculation of flat norm for $1$-boundaries embedded in higher
dimensional spaces, e.g., $\R^3$, or for input curves that are not
necessarily boundaries has not received much attention so
far. Similarly, flat norm calculations for higher dimensional input
sets have also not been well-studied. Such situations often appear in
practice -- for instance, consider the case of an input set $T$ that
is a curve sitting on a manifold embedded in $\R^3$, with choices for
$S$ restricted to this manifold as well. Further, computational
complexity of calculating flat norm in arbitrary dimensions has not
been studied. But this is not a surprising observation, given the
continuous, rather than combinatorial, setting in which flat norm
computation has been posed so far.

Simplicial complexes that triangulate the input space are often used
as representations of manifolds. Such representations use triangular
or tetrahedral meshes \cite{Edbook2006} as opposed to the uniform
square or cubical grid meshes in $\R^2$ and $\R^3$. Various simplicial
complexes are often used to represent data (in any dimension) that
captures interactions in a broad sense, e.g., the Vietoris--Rips
complex to capture coverage of coordinate-free sensor networks
\cite{SiGh2006,SiGh2007a}. It is natural to consider flat norm
calculations in such settings of simplicial complexes for denoising or
regularizing sets, or for other similar tasks. At the same time,
requiring that the simplicial complex be embedded in high dimensional
space modeled by regular square grids may be cumbersome, and
computationally prohibitive in many cases.

\subsection{Contributions} \label{ssec:ourconts}

We define a {\em simplicial} flat norm (SFN) for an input set $T$
given as a subcomplex of the finite oriented simplicial complex $K$
triangulating the set, or underlying space $\Omega$. More generally,
$T$ is the simplicial representation of a rectifiable current with
integer multiplicities. The choices of the higher dimensional sets $S$
are restricted to $K$ as well. We extend this definition to the
multiscale simplicial flat norm (MSFN) by including a scale parameter
$\lambda$. The simplicial flat norm is thus a special case of the
\MSFN{} with the default value of $\lambda=1$.

This discrete setting lets us address the worst case complexity of
computing flat norm. Given its combinatorial nature, one would expect
the problem to be difficult in arbitrary dimensions.  Indeed, we show
the problem of computing the \MSFN{} is NP-complete by reducing the
optimal bounding chain problem (OBCP), which was recently shown to be
NP-complete \cite{DuHi2011}, to a special case of the \MSFN{}
problem. We cast the problem of finding the optimal $S$, and thus
calculating the \MSFN{}, as an integer linear programming (ILP)
problem.  Given that the original problem is NP-complete, instances of
this ILP could be hard to solve. Utilizing recent work \cite{DHK2010}
on the related optimal homologous chain problem (OHCP), we provide
conditions on $K$ under which this ILP problem can in fact be solved
in polynomial time. In particular, the \MSFN{} can be computed in
polynomial time when $T$ is $d$-dimensional, and $K$ is
$(d+1)$-dimensional and orientable, for all $d \geq 0$. A similar
result holds for the case when $T$ is $d$-dimensional, and $K$ is
$(d+1)$-dimensional and embedded in $\R^{d+1}$, for all $d \geq 0$.

Our most significant contribution is the {\em simplicial deformation
  theorem} (Theorem \ref{thm:simpldeform}), which states that given an
arbitrary $d$-current in $\abs{K}$ (underlying space), we are assured
of an approximating current in the $d$-skeleton of $K$. This result is
a substantial modification and generalization of the classical
deformation theorem for currents on to square grids. Our deformation
theorem explicitly specifies the dependence of the bounds of
approximation on the regularity and size of the simplices in the
simplicial complex. Hence it is immediate from the theorem that as we
refine the simplicial complex $K$ while preserving the bounds on
simplicial regularity, the flat norm distance between an arbitrary
$d$-current in $\abs{K}$ and its deformation onto the $d$-skeleton of
$K$ vanishes. More importantly, such refinement of $K$ does not affect
the efficient computability of the \MSFN{} by solving the associated
ILP in many cases, e.g., when $K$ is orientable or when it is
full-dimensional.

\subsection{Work on Related Problems} \label{ssec:relprobs}

The problem of computing \MSFN{} is closely related to two other
problems on chains defined on simplicial complexes -- the optimal
homologous chain problem (OHCP) and the optimal bounding chain problem
(OBCP). Given a $d$-chain $\vt$ of the simplicial complex $K$, the
\OHCP{} is to find a $d$-chain $\vx$ that is homologous to $\vt$ such
that $\norm{\vx}_1$ is minimal. In the \OBCP{}, we are given a
$d$-chain $\vt$ of $K$, and the goal is to find a $(d+1)$-chain $\vs$
of $K$ whose boundary is $\vt$ and $\norm{\vs}_1$ is minimal. The
\OBCP{} is closely related to the problem of finding an
area-minimizing surface with a given boundary
\cite{Morgan2008}. Computing the \MSFN{} could be viewed, in a simple
sense, as combining the objectives of the corresponding optimal
homologous chain and optimal bounding chain problem instances, with
the scale factor determining the relative importance of one objective
over the other.

When $\vt$ is a cycle and the homology is defined over $\Z_2$, Chen
and Freedman showed that the \OHCP{} is NP-hard \cite{ChFr2010a}. Dey,
Hirani, and Krishnamoorthy \cite{DHK2010} studied the original version
of the \OHCP{} with homology defined over $\Z$, and showed that the problem
is in fact solvable in polynomial time when $K$ satisfies certain
conditions (when it has no relative torsion).  Recently, Dunfield and
Hirani \cite{DuHi2011} have shown that the \OHCP{} with homology defined
over $\Z$ is NP-complete. We will use their results to show that the
problem of computing the \MSFN{} is NP-complete (see Section
\ref{ssec:npcompMSFN}). These authors also showed that the \OBCP{} with
homology defined over $\Z$ is NP-complete as well. Their result builds
on the previous work of Agol, Hass, and Thurston \cite{AgHaTh2006},
who showed that the knot genus problem is NP-complete, and a slightly
different version of the least area surface problem is NP-hard.

The standard simplicial approximation theorem from algebraic topology
describes how continuous maps are approximated by simplicial maps that
satisfy the star condition \cite[\S14]{Munkres1984}. Our simplicial
deformation theorem applies to currents, which are more general
objects than continuous maps. More importantly, we present explicit
bounds on the expansion of mass of the current resulting from
simplicial approximation. In his PhD thesis, Sullivan
\cite{Sullivan1990} considered deforming currents on to the boundary
of convex sets in a cell complex, which are more general than the
simplices we work with. But simplicial complexes admit efficient
algorithms more naturally than cell complexes. We adopt a different
approach for deformation from Sullivan and obtain new bounds on the
approximations (see Section \ref{ssec:compbnds}). Along with the
\MSFN{}, our deformation theorem also establishes how the \OHCP{} and
\OBCP{} could be used on general continuous inputs by taking
simplicial approximations, thus expanding widely the applicability of
this family of techniques.

\section{Definition of Simplicial Flat Norm} \label{sec:defsfn}

Consider a finite $p$-dimensional simplicial complex $K$ triangulating
the set $\Omega$, where the simplices are oriented, with $p \geq
d+1$. The set $T$ is defined as the integer multiple of an oriented
$d$-dimensional subcomplex of $K$, representing a rectifiable
$d$-current with integer multiplicity. Let $m$ and $n$ be the number
of $d$- and $(d+1)$-dimensional simplices in $K$, respectively. The
set $T$ is then represented by the $d$-chain $\sum_{i=1}^m t_i
\sigma_i$, where $\sigma_i$ are all $d$-simplices in $K$ and $t_i$ are
the corresponding {\em weights}. We will represent this chain by the
vector of weights $\vt \in \Z^m$. We use bold lower case letters to
denote vectors, and the corresponding letter with subscript to denote
components of the vector, e.g., $\vx = [x_j]$. For $\vt$ representing
the set $T$ with integer multiplicity of one, $t_i \in \{-1,0,1\}$
with $-1$ indicating that the orientations of $\sigma_i$ and $T$ are
opposite. But $t_i$ can take any integer value in general. Thus, $\vt$
is the representation of $T$ in the elementary $d$-chain basis of $K$.
We consider $(d+1)$-chains in $K$ modeling sets $S$ representing
rectifiable $(d+1)$-currents with integer multiplicities, and denote
them similarly by $\sum_{j=1}^n s_j \tau_j$ in the elementary
$(d+1)$-chain basis of $K$ consisting of the individual simplices
$\tau_j$. We denote the chain modeling such a set $S$ using the
corresponding vector of weights $\vs \in \Z^n$.

Relationships between the $d$- and $(d+1)$-chains of $K$ are captured
by its $(d+1)$-boundary matrix $[\boundary_{d+1}]$, which is an $m
\times n$ matrix with entries in $\{-1,0,1\}$. If the $d$-simplex
$\sigma_i$ is a face of the $(d+1)$-simplex $\tau_j$, then the $(i,j)$
entry of $[\boundary_{d+1}]$ is nonzero, otherwise it is zero. This
nonzero value is $+1$ if the orientations of $\sigma_i$ and $\tau_j$
agree, and is $-1$ when their orientations are opposite. The $d$-chain
representing the set $T - \boundary_{d+1} S$ is then given as 
\begin{equation*}
\label{eq:xeqTminusbdyS}
\vx = \vt - [\boundary_{d+1}] \vs.
\end{equation*}
Notice that $\vx \in \Z^m$. We define the simplicial flat norm (SFN)
of $T$ represented by the $d$-chain $\vt$ in the $(d+1)$-dimensional
simplicial complex $K$ as
\begin{equation}
\label{eq:SFNdef}
F_S(T) = \min_{\vs \in \Z^n} \left\{ \sum_{i=1}^m \vol_d(\sigma_i)
\abs{x_i} + \sum_{j=1}^n \vol_{d+1}(\tau_j) \abs{s_j} ~\big|~ \vx = \vt -
[\boundary_{d+1}] \vs,\, \vx \in \Z^m \right\}.
\end{equation}
Since $\vx$ and $\vs$ are chains in a simplicial complex, the masses
of the currents they represent (as given in
Equation~\ref{eq:flatnorm}) are indeed given by the weighted sums of
the volumes of the corresponding simplices. The integer restrictions
$\vx \in \Z^m$ and $\vs \in \Z^n$ are important in this definition as
we are studying currents with integer multiplicities.  The simplicial
flat norm is intuitively the problem of deforming an input chain to
another chain of least cost, where cost is determined both by the mass
of the resulting chain and the size of the deformation (constrained to
the complex) used to get it.  For instance, in a triangulation of a
manifold, we constrain ourselves to only use deformations on the
manifold.  We generalize the definition of SFN to define a {\em
  multiscale} simplicial flat norm (MSFN) of $T$ in the simplicial
complex $K$ by including a scale parameter $\lambda \in [0,\infty)$.
\begin{equation}
\label{eq:MSFNdef}
F^{\lambda}_S(T) = \min_{\vs \in \Z^n} \left\{ \sum_{i=1}^m
\vol_d(\sigma_i) \abs{x_i} + \lambda \left( \sum_{j=1}^n \vol_{d+1}(\tau_j)
\abs{s_j} \right) ~\big|~ \vx = \vt - [\boundary_{d+1}] \vs,\, \vx
\in \Z^m \right\}.
\end{equation}
This definition is the simplicial version of the multiscale flat norm
defined in Equation (\ref{eq:flatnormwscale}). The default, or
nonscale, simplicial flat norm in Equation~(\ref{eq:SFNdef}) is a
special case of the multiscale simplicial flat norm with the default
value of $\lambda=1$.

The (non-simplicial) flat norm with scale $\lambda > 0$ of a
$d$-dimensional current $T$ can be rewritten as $\F_{\lambda}(T) =
\lambda^d \cdot \F_{1}(T/\lambda)$.  Thus the flat norm with scale can
be thought of as the traditional flat norm applied to a scaled copy of
the input current.  An equivalent statement can be made for the
simplicial flat norm, but crucially requires that the simplicial
complex be similarly scaled.  To avoid this complex scaling issue
especially when considering all possible scales, and to simplify our
notation, we henceforth study the more general multiscale simplicial
flat norm (which also allows us to consider the $\lambda = 0$ case).

We assume the $d$- and $(d+1)$-dimensional volumes of simplices to be
any nonnegative values. For example, when $\sigma_i$ is a $1$-simplex,
i.e., edge, $\vol_1(\sigma_i)$ could be taken as its Euclidean
length. Similarly, $\vol_2(\tau_j)$ for a triangle $\tau_j$ could be
its area. For ease of notation, we denote $\vol_d(\sigma_i)$ by $w_i$
and $\vol_{d+1}(\tau_j)$ by $v_j$, with the dimensions $d$ and $d+1$
evident from the context.

\begin{remark}
\label{rem:MSFNmin}
The minimum in the definition of the \MSFN{} (Equation~\ref{eq:MSFNdef})
indeed exists. The function 
\begin{equation}
\label{eq:funcmsfn}
f^{\lambda}(T,S) = \sum_{i=1}^m w_i \abs{x_i} + \lambda \, (
\sum_{j=1}^n v_j \abs{s_j} )~~\,\mbox{ with }~~\,\vx = \vt -
[\boundary_{d+1}] \vs\,
\end{equation}
is lower bounded by zero, as it is the sum of nonnegative entries (we
have $\lambda \geq 0$). Notice that $F^{\lambda}_S(T) = \min_{S}
f^{\lambda}(T,S)$. Further, we only consider integral $\vs$ defined on
the finite simplicial complex $K$, and hence there are only a finite
number of values for this function. Hence its minimum indeed exists,
which defines the \MSFN{} of $\vt$. On the other hand, the proof of
existence of minimum in the original definition of flat norm for
rectifiable currents employs the Hahn--Banach theorem
\cite[pg.~367]{Federer1969}.
\end{remark}

\begin{figure}[ht!]
  \centering \includegraphics[scale=0.6, trim=1in 5.5in 0.3in 0.8in,
  clip] {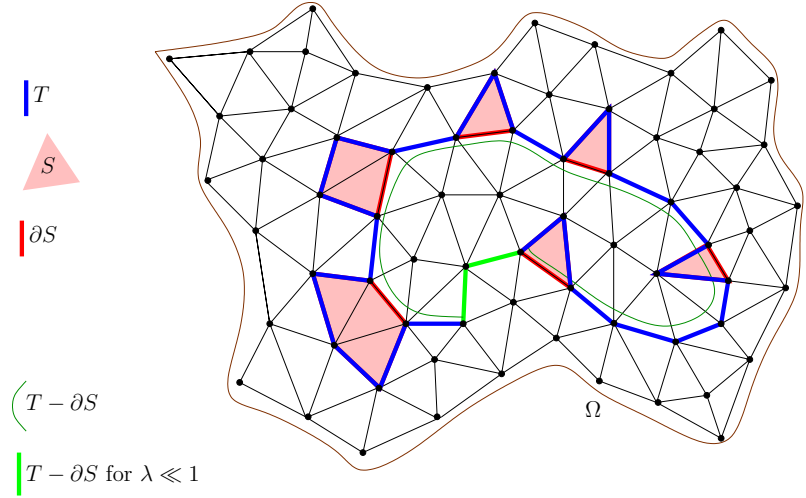}
% \centering \includegraphics[scale=0.6]{Fig_CurveOnMesh.eps}
%
  \caption[Multiscale simplicial flat norm]{\label{fig:curveonmesh} The \MSFN{} illustrated for two different
  scales $\lambda = 1$ and $\lambda \ll 1$. See text for explanation.}
\end{figure}

We illustrate the optimal decompositions to compute the \MSFN{} for
two different scales ($\lambda=1$ and $\lambda \ll 1$) in Figure
\ref{fig:curveonmesh}. Notice that the input set $T$, shown in blue,
is not a closed loop here. It is a subcomplex of the simplicial
complex triangulating $\Omega$. The underlying set $\Omega$ need not
be embedded in $\R^2$ -- it could be sitting in $\R^3$ or any higher
dimension. We do not show the orientations of individual simplices and
chains so as not to clutter the figure. We could take each triangle to
be oriented counterclockwise (CCW), with $T$ oriented CCW as well, and
each edge oriented arbitrarily. When scale $\lambda=1$, we get the
default SFN of $T$, where the $S$ chosen (shown in light pink) is such
that the resulting optimal $T - \boundary S$ (indicated by the thin
curve in dark green) is devoid of all the ``kinks'', but is similar to
$T$ in overall form.  This removal of the tightest ``kinks'' is a
discrete analogue of how the $\lambda$ in the flat norm relates to the
curvature in the continuous case.  For $\lambda \ll 1$, the second
term in the definition (Equation~\ref{eq:MSFNdef}) contributes much
less to the \MSFN{}. As such, the optimal $T- \boundary S$ consists of
a short chain of two edges (shown in light green), which closes the
original $T$ curve to form a loop. $S$ in this case includes the
triangles in the former choice of $S$, and all other triangles
enclosed by the original curve $T$ and the resulting $T -\boundary S$.

\subsection{Complexity of \MSFN{}} \label{ssec:npcompMSFN}

To study the complexity of computing the \MSFN{}, we consider a decision
version of the problem, termed {\em decision}-MSFN or \dmsfn. The
function $f^{\lambda}(T,S)$ used here is defined in Equation
\ref{eq:funcmsfn}, with the modification that $w_i$ and $v_j$ are
assumed to be rational for purposes of analyses of complexity.
\begin{definition}[{\rm \dmsfn}] 
  Given a $p$-dimensional finite simplicial complex $K$ with $p \geq
  d+1$, a set $T$ defined as a $d$-subcomplex of $K$, a scale $\lambda
  \in [0,\infty)$, and a rational number $f_0 \geq 0$, does there
    exist a $(d+1)$-dimensional subcomplex $S$ of $K$ such that
    $f^{\lambda}(T,S) \leq f_0$?
\end{definition}
The related optimal homologous chain problem (OHCP) was recently shown
to be NP-complete \cite[Theorem 1.4]{DuHi2011}. We reduce OHCP to a
special case of \dmsfn, thus showing that \dmsfn \ is NP-complete as
well. The default optimization version of MSFN consequently turns out
to be NP-hard.
\begin{theorem} \label{thm:dmsfnnpc}
  \dmsfn \ is NP-complete, and {\rm MSFN} is NP-hard.
\end{theorem}
\begin{proof}
  \dmsfn \ lies in NP as we can calculate $f^{\lambda}(T,S)$ in
  polynomial time when given a pair of $d$- and $(d+1)$-chains $\vt$
  and $\vs$, respectively, of the simplicial complex $K$.  On the
  other hand, given an instance of the optimal homologous chain
  decision problem, we can reduce it to the \dmsfn{} by taking
  $\lambda = 0$ and $w_i = 1$ for $1 \leq i \leq m$.  Since the
  \OHCP{} was recently shown to be NP-complete~\cite[Theorem
    1.4]{DuHi2011}, the result follows.
\end{proof}

\begin{remark}
  Although we showed MSFN is NP-hard in general, the case for any
  particular $\lambda > 0$ is not known.  For $\lambda$ large enough,
  the problem in fact becomes easy-- when the $(d+1)$-simplices have
  positive volumes and $\lambda > (\sum w_i) / \min v_j$, then
  optimality occurs when $\vs$ is the empty $(d+1)$-chain.
\end{remark}

We now consider attacking the \MSFN{} problem using techniques from
the area of discrete optimization. Even though the problem is NP-hard,
this approach helps us to identify special cases in which we can
compute the \MSFN{} in polynomial time.

\section{Multiscale Simplicial Flat Norm and Integer Linear Programming} \label{ssec:MSFNILP}

The problem of finding the \MSFN{} of the $d$-chain $\vt$
(Equation~\ref{eq:MSFNdef}) can be cast formally as the following
optimization problem.
\begin{equation}
\label{eq:optprobMSFN}
\begin{array}{ll}
\mbox{minimize} & \sum_{i=1}^m w_i \abs{x_i} + \lambda
     ( \sum_{j=1}^n v_j \abs{s_j} ) \\
\vspace*{-0.13in} \\
\mbox{subject to} & ~~~~~~\vx = \vt - [\boundary_{d+1}] \vs, \\
 & ~~~~~~\vx \in \Z^m,~ \vs \in \Z^n.
\end{array}
\end{equation}
The objective function is piecewise linear in the integer variables
$\vx$ and $\vs$. Using standard modeling techniques from linear
optimization \cite[pg.~18]{BeTs1997}, we can reformulate the problem
as the following integer {\em linear} program (ILP).
\begin{equation}
\label{eq:ILPMSFN}
\begin{array}{ll}
\min & \sum_{i=1}^m w_i (x^+_i + x^-_i) \,+ \,\lambda
     \left( \sum_{j=1}^n v_j (s^+_j + s^-_j) \right) \\
\vspace*{-0.13in} \\
\st & \vx^+ - \vx^-  = \vt - [\boundary_{d+1}] (\vs^+ - \vs^-) \\
 & \vx^+, \vx^- \geq \vzero,~~\vs^+, \vs^- \geq \vzero \\
 & \vx^+, \vx^- \in \Z^m,~~ \vs^+, \vs^- \in \Z^n.
\end{array}
\end{equation}
The objective function coefficients need to be nonnegative for this
formulation to work -- indeed, we have $w_i, v_j$, and $\lambda$
nonnegative. Integer linear programming is NP-complete
\cite{Schrijver1986}. The linear programming relaxation of the ILP 
above is obtained by ignoring the integer restrictions on the
variables.
\begin{equation}
  \label{eq:LPMSFN}
  \begin{array}{ll}
    \min & \sum_{i=1}^m w_i (x^+_i + x^-_i) \,+ \,\lambda
    \left( \sum_{j=1}^n v_j (s^+_j + s^-_j) \right) \\
    \vspace*{-0.13in} \\
    \st & \vx^+ - \vx^-  = \vt - [\boundary_{d+1}] (\vs^+ - \vs^-) \\
        & \vx^+, \vx^- \geq \vzero,~~\vs^+, \vs^- \geq \vzero 
  \end{array}
\end{equation}

We are interested in instances of this linear program (LP) that have
integer optimal solutions, which hence are optimal solutions for the
original ILP (Equation~\ref{eq:ILPMSFN}) as well. Totally unimodular
matrices yield a prime class of linear programming problems with
integral solutions.  Recall that a matrix is totally unimodular if all
its subdeterminants equal $-1,0,$ or $1$; in particular, each entry is
$-1,0,$ or $1$. The connection between total unimodularity and linear
programming is specified by the following theorem.
\begin{theorem}
\label{thm:IPandTU}
{\rm \cite{VeDa1968}} Let $A$ be an $m \times n$ totally unimodular
matrix, and $\vb \in \Z^m$. Then the polyhedron ${\cal P} = \{ \vx \in
\R^n \, | \, A \vx = \vb, \, \vx \geq \vzero\}$ has integral vertices.
\end{theorem}

Notice that the feasible set of the \MSFN{} LP (Equation~\ref{eq:LPMSFN})
has the form specified in the theorem above, with the variable vector
$(\vx^+, \vx^-, \vs^+, \vs^-)$ in place of $\vx$. The corresponding
equality constraint matrix $A$ has the form $\begin{bmatrix} I & -I &
B & -B \end{bmatrix}$, where $I$ is the identity matrix and $B =
[\boundary_{d+1}]$. The input $d$-chain $\vt$ is in place of the
right-hand side vector $\vb$. In order to use Theorem
\ref{thm:IPandTU} for computing the \MSFN{}, we connect the total
unimodularity of constraint matrix $A$ and that of boundary matrix
$B$.

\begin{lemma}
\label{lem:TUconmat}
If $B = [\boundary_{d+1}]$ is totally unimodular, then so is the
matrix $A = \begin{bmatrix} I & -I & B & -B \end{bmatrix}$.
\end{lemma}

\begin{proof}
Starting with $B$, we get the matrix $A$ by appending columns of $B$
scaled by $-1$ to its right, and appending columns with a single
nonzero entry of $\pm 1$ to its left. Both these classes of operations
preserve total unimodularity \cite[pg.~280]{Schrijver1986}.
\end{proof}

\noindent Consequently, we get the following result on polynomial time
computability of the \MSFN{}.

\begin{theorem}
\label{thm:MSFNpolytime}
If the boundary matrix $[\boundary_{d+1}]$ of the finite oriented
simplicial complex $K$ is totally unimodular, then the multiscale
simplicial flat norm of the set $T$ specified as a $d$-chain $\vt \in
\Z^m$ of $K$ can be computed in polynomial time.
\end{theorem}

\begin{proof}
The problem of computing the \MSFN{} of $T$ (Equation~\ref{eq:MSFNdef})
is cast as the optimization problem given in
Equation~(\ref{eq:optprobMSFN}). This problem is reformulated as an
instance of ILP (Equation~\ref{eq:ILPMSFN}). We get the \MSFN{} LP
(Equation~\ref{eq:LPMSFN}) by relaxing the integrality constraints of
this ILP.  As noted in Remark \ref{rem:MSFNmin}, the optimal cost of
this LP is finite. The polyhedron of this LP has at least one vertex,
given that all variables are nonnegative \cite[Cor.~2.2]{BeTs1997}. 
By Lemma \ref{lem:TUconmat}, the constraint matrix of this LP is
totally unimodular, as $[\boundary_{d+1}]$ is so. Hence by
Theorem~\ref{thm:IPandTU}, all vertices of the feasible region of the \MSFN{}
LP are integral, since $\vt \in \Z^m$.

An optimal solution $( \vx^+_*, \vx^-_*, \vs^+_*, \vs^-_* )$ of the
\MSFN{} LP can be found in polynomial time using an interior point
method \cite[Chap.~9]{BeTs1997}. If it happens to be a unique optimal
solution, then it will be a vertex, and hence will be integral by
Theorem \ref{thm:IPandTU}. Hence it is an optimal solution to the ILP
(Equation~\ref{eq:ILPMSFN}).

If the optimal solution is not unique, then $( \vx^+_*, \vx^-_*,
\vs^+_*, \vs^-_* )$ may be nonintegral. But since the optimal cost is
finite, there must exist a vertex in its polyhedron that has this
minimum cost. Given a nonintegral optimal solution obtained by an
interior point method, one can find such an integral optimal solution
at a vertex in polynomial time \cite{GuHeRoTeTs1993}. Hence the \MSFN{}
ILP can be solved in polynomial time in this case as well.
\end{proof}

\begin{remark}
\label{rem:strgpoly}
We point out that since the boundary matrix $B = [\boundary_{d+1}]$
has entries only in $\{-1,0,1\}$, the constraint matrix of the \MSFN{} LP
(Equation~\ref{eq:LPMSFN}) also has entries only in
$\{-1,0,1\}$. Hence the \MSFN{} LP can be solved in strongly polynomial
time \cite{Tardos1986}, i.e., the time complexity is independent of
the objective function and right-hand side coefficients, and depends
only on the dimensions of the problem.
\end{remark}

\begin{remark}
\label{rem:zeroonesoln}

Components of variables $\vx^+,\vx^-,\vs^+,\vs^-$ in the \MSFN{} ILP
(Equation~\ref{eq:ILPMSFN}) could assume values other than
$\{-1,0,1\}$, indicating integer multiplicities higher than $1$ for
the corresponding simplices in the optimal decomposition. The
definition of \MSFN{} (Equation~\ref{eq:MSFNdef}) does allow such
larger multiplicities. At the same time, if one insists on using each
$(d+1)$-simplex at most {\em once} when calculating the \MSFN{}, and
insists on similar restrictions on $d$-simplices in the optimal
decomposition, we can modify the ILP such that Theorem
\ref{thm:MSFNpolytime} still holds.

Denoting the entire variable vector by $\vx = (\vx^+, \vx^-, \vs^+,
\vs^-) \in \Z^{2m+2n}$, we add the upper bound constraints $\vx \leq
\vone$, where $\vone$ is the $(2m+2n)$-vector of ones. These
inequalities could be converted to the set of equations $\vx + \vy =
\vone$, where $\vy$ is the $(2m+2n)$-vector of slack variables that
are nonnegative. These modifications give an ILP whose polyhedron is
in the same form as described in Theorem \ref{thm:IPandTU}, with the
equations denoted as $A'\vx' = \vb'$ for the variable vector $\vx' =
(\vx,\vy)$. The new constraint matrix $A'$ is related to the
constraint matrix $A$ of the original \MSFN{} ILP given in Lemma
\ref{lem:TUconmat} as 
\[
A' = \begin{bmatrix} A & O \\ I & I \end{bmatrix},
\]
where $I$ is the $2m+2n$ identity matrix, and $O$ is the $m \times
(2m+2n)$ zero matrix. Hence $A'$ is obtained from $A$ by first adding
$2m+2n$ rows with a single nonzero entry of $+1$, and then adding to
the resulting matrix $2n+2m$ more columns with a single nonzero entry
of $+1$. These operations preserve total unimodularity
\cite[pg.~280]{Schrijver1986}, and hence the new constraint matrix
$A'$ is totally unimodular when $[\boundary_{d+1}]$ is so. The new
right-hand side vector $\vb' \in \Z^{3m+2n}$ consists of the input
chain $\vt$ and the vector of ones from the new upper bound
constraints.

\end{remark}

Since the efficient computability of the \MSFN{} depends on the total
unimodularity of the boundary matrix, we study the conditions under
which total unimodularity of boundary matrices can be guaranteed.

\section{Simplicial Complexes and Relative Torsion}

Dey, Hirani, and Krishnamoorthy \cite{DHK2010} have given a simple
characterization of the simplicial complex whose boundary matrix is
totally unimodular. In short, if the simplicial complex does not have
relative torsion then its boundary matrix is totally unimodular. We
state this and other related results here for the sake of
completeness, and refer the reader to the original paper
\cite{DHK2010} for details and proofs. The simplicial complex $K$ in
these results has dimension $d+1$ or higher. Recall that a
$d$-dimensional simplicial complex is {\em pure} \ if it consists of
$d$-simplices and their faces, i.e., there are no lower dimensional
simplices that are not part of some $d$-simplex in the complex.

\begin{theorem} \label{thm:TUreltorfree}
  {\rm \cite[Theorem 5.2]{DHK2010}} The boundary matrix
  $[\boundary_{d+1}]$ of a finite simplicial complex $K$ is totally
  unimodular if and only if $H_d(L,L_0)$ is torsion-free for all pure
  subcomplexes $L_0,L$ of $K$, with $L_0 \subset L$.
\end{theorem}

These authors further describe situations in which the absence of
relative torsion is guaranteed. The following two special cases
describe simplicial complexes for which the boundary matrix is always
totally unimodular.

\begin{theorem} \label{thm:TUorimfld}
  {\rm \cite[Theorem 4.1]{DHK2010}} The boundary matrix
  $[\boundary_{d+1}]$ of a finite simplicial complex triangulating a
  compact orientable $(d+1)$-dimensional manifold is totally
  unimodular.
\end{theorem}

\begin{theorem} \label{thm:TUembedinRdp1}
  {\rm \cite[Theorem 5.7]{DHK2010}} The boundary matrix
  $[\boundary_{d+1}]$ of a finite simplicial complex embedded in
  $\R^{d+1}$ is totally unimodular.
\end{theorem}

\noindent For simplicial complexes of dimension $2$ or lower, the
boundary matrix is totally unimodular when the complex does not have a
M\"obius subcomplex.

\begin{theorem} \label{thm:TUnomobius}
  {\rm \cite[Theorem 5.13]{DHK2010}} For $d \leq 1$, the boundary
  matrix $[\boundary_{d+1}]$ is totally unimodular if and only if the
  finite simplicial complex has no $(d+1)$-dimensional M\"obius
  subcomplex.
\end{theorem}

It is appropriate to mention here that the connection between total
unimodularity of boundary matrices and torsion in the complex has been
observed as early as in 1895 by
Poincar\'{e}\cite{Poincare2010}. However, the result in~\cite{DHK2010}
connecting the total unimodularity with relative torsion is different
and has led to a polynomial time algorithm for the OHCP
problem. Notice that a complex can be torsion-free, but have
non-trivial relative torsion. The M\"{o}bius strip is such an example.

We illustrate the implications of the results above for the efficient
computation of the \MSFN{} by considering certain sets. When the input
set $T$ is of dimension 1, and is described on an orientable
$2$-manifold to which the choices of $2$-dimensional set $S$ are also
restricted, we can always compute its \MSFN{} by solving the \MSFN{}
LP (Equation~\ref{eq:LPMSFN}) in polynomial time.  A similar result
holds when $T$ is a set of dimension $2$ described as a subcomplex of
a $3$-complex sitting in $\R^3$. For a $1$-dimensional set $T$ with
choices of $S$ restricted to a $2$-complex $K$, we can always compute
the \MSFN{} of $T$ efficiently as long as $K$ does not have a
$2$-dimensional M\"obius subcomplex. Notice that $K$ itself need not
be embedded in $\R^3$ for this result to work -- it could be sitting
in some higher dimensional space.

\section{Simplicial Deformation Theorem} \label{sec:simpdefthm}

When can we use the multiscale simplicial flat norm as a discrete
surrogate for the traditional flat norm?  That is, if we wish to solve
a flat norm problem (for which there are no practical algorithms in
general), can we discretize the problem and find a problem close
enough to the original one which we can solve?

The deformation theorem \cite[Sections 4.2.7--9]{Federer1969} is one
of the fundamental results of geometric measure theory, and more
particularly of the theory of currents. It approximates an integral
current by deforming it onto a cubical grid of appropriate mesh
size. On the other hand, we have been studying currents or sets in the
setting of simplicial complexes, rather than on square grids. Our
proof is a substantial modification of the classical proof of the
deformation theorem. We found the presentation of the latter proof by
Krantz and Parks \cite[Section 7.7]{KrantzParks2008} especially
helpful. Our proof mimics their proof when possible. The gist of this
theorem is the assertion that we may approximate a current with a
simplicial current.

Recall that $\vol_d(\sigma)$ denotes the $d$-dimensional volume of a
$d$-simplex $\sigma$. The {\em perimeter} of $\sigma$ is the set of
all its $(d-1)$-dimensional faces, denoted as $\perimeter(\sigma) = \{
\cup_j \tau_j \,|\, \tau_j \in \sigma, \dim(\tau_j) = d-1\}$. We will
also refer to the $(d-1)$-dimensional volume of $\perimeter(\sigma)$
as the perimeter of $\sigma$, but denote it as $\per(\sigma) =
\sum_{\tau_j \in \perimeter(\sigma)} \vol_{d-1}(\tau_j)$. We let
$\diam(\sigma)$ be the diameter of $\sigma$, which is the largest
Euclidean distance between any two points in $\sigma$.

\begin{theorem}[{\bf Simplicial Deformation Theorem}] \label{thm:simpldeform}
  Let $K$ be a $p$-dimensional simplicial complex embedded in $\R^q$,
  with $p=d+k$ for $k \geq 1$ and $q \geq p$. Suppose that for every
  simplex $\sigma \in K$
  \[ \frac{ \diam(\sigma) \per(\sigma)} {\vol_d({\cal B_\sigma})} \leq
  \upkappa_1 < \infty,\]
  \[ \frac{\diam(\sigma)}{r_{\sigma}} \leq \upkappa_2 < \infty,\] and
  \[\diam(\sigma) \leq \Updelta\] 
  hold, where $\hat{\cal B}_\sigma$ is the largest ball inscribed in
  $\sigma$, ${\cal B_\sigma}$ is the ball with half the radius and
  same center as $\hat{\cal B}_\sigma$, and $r_{\sigma}$ is the radius
  of ${\cal B_\sigma}$. Let $T$ be a $d$-dimensional current in $\R^q$
  such that the support of $T$ is a subset of the underlying
  space of $K$. Suppose that $T$ satisfies
  \[ \mass(T) + \mass(\boundary T) < \infty.\] 
  Then there exists a simplicial $d$-current $P$ supported in the
  $d$-skeleton of $K$ whose boundary $\boundary P$ is supported in the
  $(d-1)$-skeleton of $K$ such that
  \[ T-P = Q + \boundary R,\] 
  and the following controls on mass $M$ hold: 
  \begin{align}
    \mass(P) & \leq (4\upvartheta_K)^{k} \mass(T)+ \Updelta(4\upvartheta_{K})^{k+1}\mass(\boundary T),
    \label{eq:simpdefthmMP}\\
    \mass(\boundary P) & \leq (4\upvartheta_K)^{k+1} \mass(\boundary T), 
    \label{eq:simpdefthmMbdyP}\\
    \mass(R) & \leq \Updelta(4\upvartheta_K)^{k} \mass(T), \mbox{ and}
    \label{eq:simpdefthmMR} \\
    \mass(Q) & \leq  \Updelta(4\upvartheta_{K})^{k}(1+4\upvartheta_{K})\mass(\boundary T), 
    \label{eq:simpdefthmMQ}
  \end{align}
  where $\upvartheta_K = \upkappa_1 + \upkappa_2$.
\end{theorem}

\begin{remark}
\label{rem:defthm1}
  It is immediate that the flat norm distance between $T$ and $P$ can
  be made arbitrarily small by subdividing the simplicial complex to
  reduce $\Updelta$ while preserving the regularity of the refinement
  as measured by $\upkappa_1$ and $\upkappa_2$.
\end{remark}
\begin{remark}
  Note that this theorem combines the unscaled and scaled versions of
  the original deformation theorem \cite[Theorems 7.7.1 and
  7.7.2]{KrantzParks2008} into one theorem through the explicit form
  of the constraints. In our proof of Theorem \ref{thm:simpldeform},
  we replace certain pieces of the original proof as presented by
  Krantz and Parks~\cite[Pages 211--222]{KrantzParks2008} without
  reproducing all the other details of their proof. We found their
  exposition quite well-structured, making it easier to identify the
  modifications needed to get our theorem.
\end{remark}
\begin{remark}
  The bound for $\mass(P)$ in Theorem \ref{thm:simpldeform} is larger
  than the classical bound. We get this large bound because we
  generate $P$ through retractions alone, and not using the usual
  Sobolev-type estimates \cite[Pages 220--222]{KrantzParks2008}. And
  of course, the $\Updelta$ in the coefficient of the extra term means
  that it becomes unimportant as the simplicial complex is
  appropriately subdivided.
\end{remark}

\subsection{Proof of the Simplicial Deformation Theorem}

At the heart of the modification of the deformation theorem (from
cubical grid to simplicial complex settings) is the recalculation of
an integral over the current and its boundary. This integral appears
in a bound on the {\em Jacobian of the retraction}, which measures the
expansion in mass of the current resulting from the process of
retracting it on to the simplices of the simplicial complex. To do
this recalculation, we consider the retraction $\phi$ one step at a
time, building it through independent choices of centers to project
from in every simplex and its every face.

We first describe the general set up of retraction within a
simplex. We then present certain bounds on the mass expansion
resulting from the retraction in Lemmas \ref{lem:bnddjacob},
\ref{lem:bndJxfixa}, and \ref{lem:bndJafixx}. In particular, we obtain
bounds on the expansion that are independent of the choice of points
from which we project. These bounds are independent of the particular
current that we retract on to the simplicial complex. But we employ
these bounds to subsequently bound the overall expansion of mass of
the current resulting from the retraction.

\subsubsection{Retracting from a center inside a simplex} \label{sssec:retrsetup}

\begin{figure}[htp!]
  \centering
  \includegraphics[scale=1, trim=1in 7.2in 3in 1in, clip]{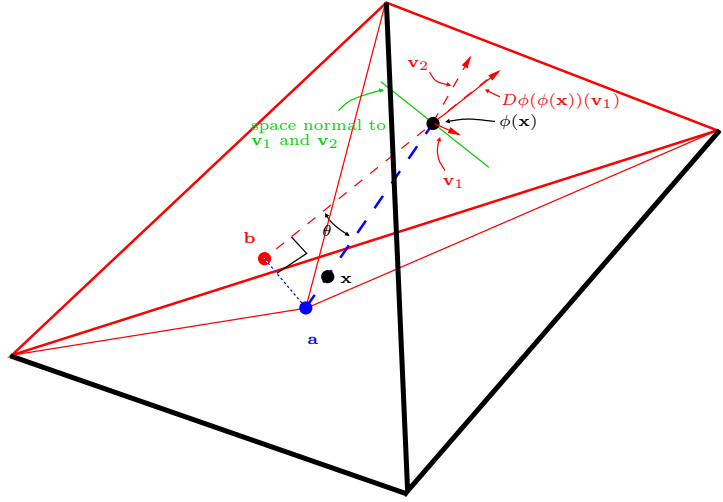}
  \caption[Retraction for a $3$-simplex]{Illustration of the dilation and nonorthogonal projection
  involved in retraction for a $3$-simplex.}  
  \label{fig:jacob}
  \end{figure}

\begin{figure}[ht!]
  \centering
  \includegraphics[scale=0.9, trim=1in 6.9in 2.2in 1in, clip]{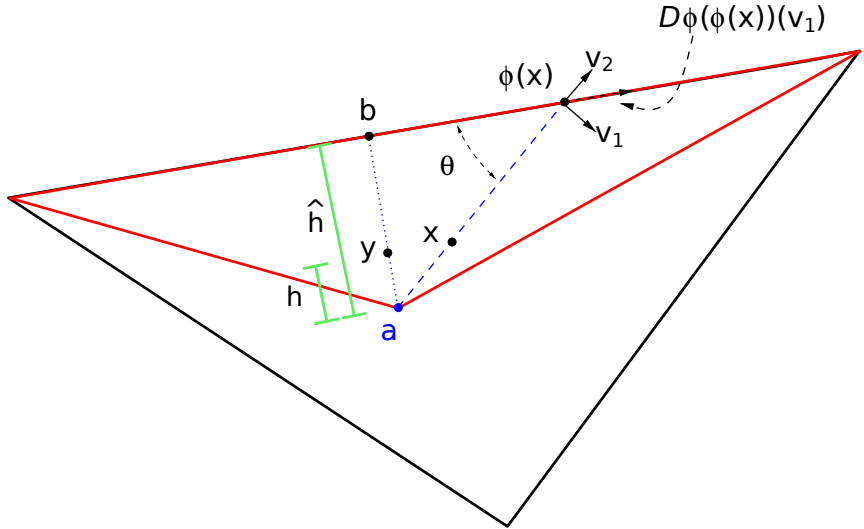}
  \caption[$2$-dimensional illustration of dilation calculation]{A $2$-dimensional illustration of the dilation
  calculation.}
  \label{fig:2d-diagram} 
  \end{figure}

  We describe the details of retraction for an $\ell$-simplex $\sigma$
  in the $p$-dimensional simplicial complex $K$. This set up is valid
  for any $\ell$, but in particular, we will use the bounds thus
  obtained for $d \leq \ell \leq p$ when retracting a $d$-current onto
  $K$. We pick a center $\va \in \intrr(\sigma)$, the interior of
  $\sigma$, and project every $\vx \in \intrr(\sigma) \setminus
  \{\va\} $ along the ray $(\vx-\va)/\norm{\vx-\va}$ to
  $\perimeter(\sigma)$. Denoting this map as $\phi(\vx,\va)$, we get
  \begin{equation} \label{eq:phixa}
  \phi(\vx,\va) = (\phi_{\pi} \circ \phi_{\delta})(\vx,\va),
  \end{equation}

  where $\phi_{\delta}(\vx,\va)$ is a dilation of $\R^\ell$ by the
  factor $\norm{\phi(\vx,\va)-\va}/\norm{\vx-\va}$ and
  $\phi_{\pi}(\vx, \va)$ is a nonorthogonal projection along $(\vx -
  \va)/ \norm{\vx-\va}$ onto $\tau_{\vx}$, the $(\ell-1)$-dimensional
  face of $\sigma$ containing $\phi(\vx,\va)$. We denote $\hat{r} =
  \norm{\phi(\vx,\va) - \va}$ and $r = \norm{\vx-\va}$. Let $E_{\ell}$
  be the $\ell$-hyperplane that contains $\sigma$ and $E_{\ell-1}$ the
  $(\ell-1)$-hyperplane that contains $\tau_{\vx}$. Denote the
  orthogonal projection of $\va$ onto $E_{\ell-1}$ by $\vb$, and let
  $\hat{h} = \norm{\vb - \va}$. For any point $\vy = \va + (\vb - \va)
  \gamma$ with $0 < \gamma < 1$, we get $\phi(\vy,\va) = \vb$. In
  particular, we consider the point of intersection of line connecting
  $\va$ and $\vb$ with the $(\ell-1)$-hyperplane parallel to
  $\tau_{\vx}$ that contains $\vx$. Naming this point $\vy$, we define
  $h = \norm{\vy - \va}$. Let $\vz \in E_{\ell}$ denote either normal
  to $\tau_{\vx}$ at $\phi(\vx,\va)$ (either of the two possibilities
  work). Let $\vv_2 = (\vx - \va)/ \norm{\vx-\va}$, and let $\vv_1$ be
  the vector in $\spn ( \vz, \vv_2 )$ that is normal to $\vv_2$ and
  points into $\sigma$. We illustrate this construction on a
  $3$-simplex in Figure~\ref{fig:jacob}, where the cone of $\va$ with
  face $\tau$ is shown in red and the other points and vectors are
  labeled. We also illustrate the corresponding slice spanned by
  $\vv_1$ and $\vv_2$ in Figure~\ref{fig:2d-diagram}.

  Choose an orthogonal basis $\{\vw_1,...,\vw_{\ell-2}\}$ for
  $\spn(\vv_1,\vv_2)^{\perp}$.  Note that $\spn(\vw_1, ...,$
  $\vw_{\ell-2}) \subset E_{\ell-1}$. Let $\vw'$ be a unit vector in
  $\spn(\vw_1,...,\vw_{\ell-2})^{\perp}\cap E_{\ell-1}$ parallel
  to $\,\phi(\vx) - \vb$. Then $\{\vv_1, \vw_1, ..., \vw_{\ell-2},$
  $\vv_2\}$ is an orthogonal basis for $\R^{\ell}$, and $\phi_{\pi}$
  is given by
  \begin{equation} \label{eq:orthbasisRl}
     \begin{array}{rcl}		
     \phi_{\pi}(\vv_1) & = &  \alpha \vw',\\
     \phi_{\pi}(\vw_i) & = & \vw_i, \; i \in \{1,...,\ell-2\},~~\mbox{ and }\\
     \phi_{\pi}(\vv_2) & = & 0, 
     \end{array}
  \end{equation}
  where $\alpha = \hat{r}/\hat{h}$. Notice that the above set up works
  everywhere except when $\phi(\vx) = \vb$, in which case we obtain an
  orthogonal projection for $\phi_{\pi}(\vx)$ along
  $\vb-\va$. Choosing coordinates for the tangent spaces of $\sigma$
  and $\tau_{\vx}$ to be $\{\vv_1, \vw_1, ..., \vw_{\ell-2}, \vv_2\}$
  and $\{\vw', \vw_1,..., $ $\vw_{\ell-2}\}$, respectively, we get
  from Equation~(\ref{eq:orthbasisRl}) that $D \phi_{\pi}(\vx,\va)$ is
  the $(\ell-1) \times \ell$ matrix given as
  \begin{equation} \label{eq:Dphimatrix}
  D\phi_{\pi}(\vx,\va) = 
   \begin{bmatrix} %\left( \begin{array}[htp!]{llllll} 
    \alpha & 0 & 0 & ... & 0 & 0 \\
	 0 & 1 & 0 & ... & 0 & 0 \\ 
    \vdots & &\ddots & & \vdots & \vdots \\ 
      0 & \hdots  & & & 1 & 0
    \end{bmatrix}. 
  \end{equation}
%    \end{array} \right). \]
%
\subsubsection{Bounding the Integral of the Jacobian}

We now present a series of bounds on integrals of the dilation of
$d$-volumes induced by the retraction. Since $\ell=d$ implies we are
already in the $d$-skeleton and no retraction is needed, we can assume
that $\ell > d$. We start with a bound on the maximum dilation of
$d$-volumes under the retraction $\phi$. $D\phi$ will denote the
\emph{tangent map} or \emph{Jacobian map} of $\phi$.
\begin{definition} 
  Let $J_d\phi(\vx,\va)$ be the maximum dilation of $d$-volumes
  induced by $D\phi(\vx,\va)$ at $\vx$.
\end{definition}

We will use the definitions and results on $D\phi_{\pi}(\vx,\va)$ in
$\ell$-dimension given above. In particular, recall that
$\diam(\sigma)$ is the diameter of $\sigma$, $\hat{h} = \|b-a\|$ and $h = \|y-a\|$.

  \begin{lemma} \label{lem:bnddjacob} For any center $\va$ and any
  point $\vx \neq \va$ in the $\ell$-simplex $\sigma$ with $d < \ell
  \leq p=d+k$,
  \[ J_d \phi(\vx,\va) \leq \left(\frac{\hat{h}}{h}\right)^d
  \frac{\diam(\sigma)}{\hat{h}}. \]
  \end{lemma}
  \begin{proof}
  Following Equation~(\ref{eq:phixa}), we seek bounds on $D
  \phi_{\delta}(\vx,\va)$ and $D \phi_{\pi}(\vx,\va)$. Since
  $D\phi_{\delta} (\vx,\va)$ simply scales by $\hat{r}/r = \hat{h}/h$,
  the expansion of $d$-volume of any $d$-hyperplane by
  $D\phi_{\delta}(\vx, \va)$ is by a factor of $(\hat{h}/h)^d$.  On
  the other hand, bounding the dilation that $D\phi_{\pi}(\vx,\va)$
  can cause in $d$-hyperplanes is a little more involved. We seek a
  bound on
  \begin{equation} \label{eq:dilratio}
  \frac{\sqrt{\det(\, (D\phi_{\pi}(\vx,\va) U)^T (D\phi_{\pi}(\vx,\va)
  U) \,)}} {\sqrt{\det(U^T U)}}
  \end{equation}
  for all $\ell \times d$ matrices $U$. Using the generalized
  Pythagorean theorem \cite[Section 1.5]{KrantzParks2008}, we get
  \[ \det(U^T U) =  \sum_{\lambda\in\Lambda} (\det(U_{\lambda}))^2 \]
  where submatrix $U_{\lambda}$ consists of the $d$ rows of $U$
  specified by the set of index maps $\Lambda$ given as
  \[ 
  \lambda \in \Lambda\equiv \{f | f:[1,...,d] \rightarrow
  [1,...,\ell], ~f \text{ is one to one and increasing}\}.
  \]
  A similar result holds for $\det( (D\phi_{\pi}(\vx,\va) U)^T
  (D\phi_{\pi}(\vx,\va) U) )$, with the functions $f$ considered
  mapping $[1,...,d]$ to $[1,...,\ell-1]$.

  Observe that multiplying by $D\phi_{\pi}(\vx,\va)$ (Equation
  \ref{eq:Dphimatrix}) just scales the first row of $U$ by $\alpha$
  and removes the last row. Thus $\alpha\det(U_{\lambda}) \geq
  \det((D\phi_{\pi}(\vx,\va)U)_{\lambda})$, which implies that
  $\alpha$ is a bound on the ratio in Equation~(\ref{eq:dilratio}).
  Thus we have that
  \[ J_d \phi(\vx,\va) = \left( \frac{\hat{r}}{r}\right)^d
  \frac{\norm{\phi(\vx,\va) - \va}} {\norm{\vb-\va}} = \left(
  \frac{\hat{h}}{h} \right)^d \frac{\norm{\phi(\vx,\va) -
  \va}}{\norm{\vb-\va}} \leq \left(\frac{\hat{h}}{h}\right)^d
  \frac{\diam(\sigma)}{\hat{h}} \]
  holds for all $\vx$ and $\va$ in $\sigma$, where $\diam(\sigma)$ is
  the diameter of the $\ell$-simplex $\sigma$. \end{proof}

  Next we describe a bound on the integral of $J_d \phi(\vx,\va)$ over
  the entire $\ell$-simplex, for a fixed center $\va$. We will find
  that this bound is independent of the position of $\va$. Recall that
  $\perimeter(\sigma)$ and $\per(\sigma)$ denote the perimeter of
  $\ell$-simplex $\sigma$ and the $(\ell-1)$-dimensional volume of the
  perimeter, respectively, and $\intrr(\sigma)$ its interior.

  \begin{lemma} \label{lem:bndJxfixa} 
    For any fixed center $\va$ in the $\ell$-simplex $\sigma$ with $d
    < \ell \leq p=d+k$,
    \[ \int_{\intrr(\sigma)} J_d \phi(\vx,\va) \, {\rm d}{\mathcal
    L}^{\ell}(\vx) \, \leq \, \diam(\sigma) \per(\sigma). \]
  \end{lemma}
  \begin{proof} 
  Consider the $(\ell-1)$-dimensional faces $\tau_j$ of $\sigma$, with
  $\perimeter(\sigma) = \{ \cup_j \tau_j \,|\, \tau_j \in \sigma,\,
  \dim(\tau_j)=\ell-1\}$. Let $\sigma_j$ denote the $\ell$-simplex
  generated by $\va$ and $\tau_j$. Then
  \[\int_{\intrr(\sigma)} J_d \phi(\vx,\va) \, {\rm d}{\mathcal
  L}^{\ell}(\vx) = \sum_j \int_{\intrr(\sigma_j)} J_d \phi(\vx,\va) \,
  {\rm d}{\mathcal L}^{\ell}(\vx).\]
  Let $\tau_j(h)$ denote the $(\ell-1)$-simplex formed by the
  intersection of $\sigma_j$ and the $(\ell-1)$-hyperplane parallel to
  $\tau_j$ at a distance $h$ from $\va$. Thus, $\tau_j(\hat{h})$ is
  $\tau$ itself. We observe that our bound on $J_d \phi(\vx,\va)$ is
  constant in $\tau_j(h)$ for any $h$. The $(\ell-1)$-dimensional
  volume of $\tau_j(h)$ is given by
  \[ \vol_{\ell-1}(\tau_j(h)) =
  \left(\frac{h}{\hat{h}}\right)^{\ell-1} \vol_{\ell-1} (\tau_j).\]
  Using the bound on $J_d \phi(\vx,\va)$ from Lemma
  \ref{lem:bnddjacob}, and noting that $\diam(\sigma_j) \leq
  \diam(\sigma)~ \forall \,j$, we get
  \begin{align*} \int_{\intrr(\sigma_j)} J_d \phi(\vx,\va) \, {\rm d}{\mathcal
  L}^{\ell}(\vx) &\leq \int_0^{\hat{h}}
  \left(\frac{h}{\hat{h}}\right)^{\ell-1} \vol_{\ell-1} (\tau_j) \,
  \left(\frac{\hat{h}}{h}\right)^d \frac{\diam(\sigma)}{\hat{h}} \,
  {\rm d}h\\ &= \frac{\vol_{\ell-1}(\tau_j) \diam(\sigma)}{\ell-d}.
  \end{align*}
  Summing this quantity over all $\tau_j \in \perimeter(\sigma)$ and
  replacing $\ell-d \geq 1$ with $1$ gives the overall
  bound.  \end{proof}

  We now bound the integral of $J_d \phi(\vx,\va)$ over centers $\va$
  with a fixed $\vx$ that we are retracting onto $\perimeter(\sigma)$.
  Examination of the corresponding proof for the original deformation
  theorem \cite[Section 7.7]{KrantzParks2008} shows that symmetry of
  the cubical mesh plays a very special role, which cannot be
  duplicated in the case of simplicial complex. In particular, we must
  avoid integrating over $\va$ close to the perimeter of
  $\sigma$. Hence we integrate over as big a region as we can while
  still avoiding a neighborhood of the perimeter. As in the statement
  of the main Theorem \ref{thm:simpldeform}, let $\hat{\cal B}_\sigma$
  be the largest ball inscribed in $\sigma$, ${\cal B}_\sigma$ be the
  ball with half the radius and same center as $\hat{\cal B}_\sigma$,
  and $r_{\sigma}$ be the radius of ${\cal B}_\sigma$.
  \begin{lemma} \label{lem:bndJafixx}
    For any point $\vx$ in the $\ell$-simplex $\sigma$ with $d < \ell
    \leq p=d+k$,
     \[ \int_{{\cal B}_\sigma} J_d \phi(\vx,\va) \, {\rm d}{\cal
     L}^{\ell} (\va) \, \leq \, \diam(\sigma) \per(\sigma) +
     \vol_{\ell}({\cal B}_\sigma) \frac{\diam(\sigma)}{r_{\sigma}}. \]
  \end{lemma}
  \begin{proof} 
  Similar to the subsimplices of $\sigma$ considered in the Proof of
  Lemma \ref{lem:bndJxfixa}, let $\sigma_j$ now denote the
  $\ell$-simplex formed by $\vx$ and $\tau_j \in \perimeter(\sigma)$.
  In order to derive an upper bound, we integrate instead over regions
  that are by construction bigger than these subsimplices of $\sigma$.
  Denoting the simplex $\sigma_j$ as Region 1, we define Regions 2 and
  3 as follows. We refer the reader to Figure~\ref{fig:integrate-a}
  for an illustration of this construction. Let $\sigma'_j$ be the
  reflection of $\sigma_j$ through $\vx$, and similarly, let $\tau'_j$
  be the reflection through $\vx$ of $\tau_j$. We define $\sigma'_j$
  as Region 2. Notice that unlike Region 1, Region 2 need not be
  contained fully in $\sigma$. As defined in Section
  \ref{sssec:retrsetup}, let $\vz$ be the unit vector normal to the
  $(\ell-1)$-hyperplane containing $\tau_j$ pointing into $\sigma$. We
  define Region 3 as the $\ell$-dimensional set $\,\tau'_j +
  [0,\diam(\sigma)]\vz$, as illustrated in
  Figure~\ref{fig:integrate-a}.

  Note that the union of all Region 2's and Region 3's cover $\sigma$.
  By an argument almost identical to that above, we have the following
  upper bound on the integrand in question.
  \begin{equation*} \label{eq:int-a-bound} 
  \left(\frac{h'}{h}\right)^d \frac{\norm{\phi(\vx,\va) -
  \va}}{\hat{h}} \, \leq \, \left( \frac{h'}{h} \right)^d
  \frac{\diam(\sigma)}{\hat{h}}.
  \end{equation*}
  \begin{figure}[ht!]
  \centering 
  \includegraphics[scale=0.6, trim=1in 3.2in 1.5in 1in,
  clip]{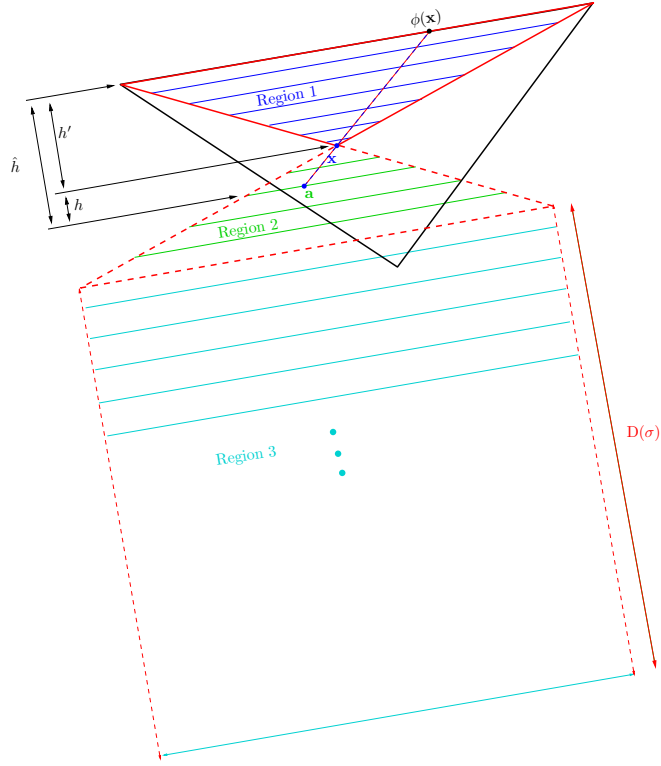} 
  \caption[Integrating Jacobian bound over centers]{Illustration for integration of Jacobian bound over centers
  $\va$ instead of $\vx$'s. For the case of the triangle shown, there
  will be 3 sets of 3 regions. In general there will be 3 regions for
  every face of the simplex.}
  \label{fig:integrate-a} 
  \end{figure}
  Integrating the second of these two terms over Region 2 and summing
  the integral over all such Regions 2 for all faces $\tau_j$, we get
  the upper bound of $\diam(\sigma)\per(\sigma)$. Here we use the same
  arguments as the ones employed in Lemma \ref{lem:bndJxfixa}. Region
  2 alone is not guaranteed to cover ${\cal B}_{\sigma}$ as some of
  ${\cal B_\sigma}$ may occupy parts of Region 3. Since $\va \in {\cal
    B}_\sigma$, we have $\hat{h} > r_{\sigma}$, and $h' \leq h$ when
  $\va \in $Region 3, so that
  \[\left(\frac{h'}{h}\right)^d \frac{\diam(\sigma)}{\hat{h}} \leq
  \frac{\diam(\sigma)}{r_{\sigma}}.\]
  Combining the above estimates while integrating over all such
  Regions 2 and 3 gives us the bound specified in the
  Lemma. \end{proof}

\subsubsection{Bounding the pushforwards of the current}

  We consider the $d$-current $T$, and employ the bounds on the
  Jacobian of retraction described above to the pushforwards of $T$
  and its boundary $\boundary T$ on to the simplicial complex $K$.
  Our treatment of the pushforwards essentially follows the
  corresponding results of Krantz and Parks for the case of square
  grid \cite[Pages 218--219]{KrantzParks2008}. We denote by $\norm{T}$
  the total variation measure of the current $T$, which is determined
  by the identity
  \[
  \norm{T}(W) = \sup_{\text{\parbox{9em}{\centering$\omega \in
        \mathcal{D}^d$, $\norm{\omega} = 1$,\\$\spt\omega \subset
        W$}}} T(\omega).
  \]
  \begin{lemma} \label{lem:retTbnd} 
  Suppose $K$ is a $p$-dimensional simplicial complex with $p = d+k$
  for $k \geq 1$. Consider the stepwise retraction of the $d$-current
  $T\subset K$ (the $(d-1)$-current $\boundary T \subset K$) onto the
  $d$-skeleton of $K$ (respectively, the $(d-1)$-skeleton of
  $K$). Each step of the retraction on to the perimeter of an
  $\ell$-simplex $\sigma$ for $d < \ell \leq p$ (respectively, $d \leq
  \ell \leq p $) increases the mass of $T$ or $\boundary T$ by at most
  a factor of
  \[ 4 \upvartheta_K = 4(\upkappa_1 + \upkappa_2) = 4 \max_{\sigma \in
  K} \left( \frac{\diam(\sigma)\per(\sigma)}{\vol_{\ell}({\cal
  B}_\sigma)} + \frac{\diam(\sigma)}{r_{\sigma}} \right).\]
  \end{lemma}

  \begin{proof} 
  Using Fubini's theorem \cite[Page 26]{KrantzParks2008} and applying
  the bound in Lemma \ref{lem:bndJafixx}, we get
  \[ \int_{{\cal B}_\sigma} \int_{\sigma} J_d \phi(\vx,\va) \, {\rm
  d}\norm{T}(\vx) \, {\rm d} {\cal L}^{\ell}(\va) = \int_{\sigma}
  \int_{{\cal B}_\sigma} J_d \phi(\vx,\va) \, {\rm d}{\cal
  L}^{\ell}(\va) \, {\rm d} \norm{T}(\vx) \, \leq \,
  \upvartheta_{\sigma} \mass(T|_{\sigma}), \]
  where $\upvartheta_{\sigma} = \diam(\sigma)\per(\sigma) +
  \vol_d({\cal B}_\sigma) (\diam(\sigma)/r_{\sigma})\,$ and
  $T|_{\sigma}\,$ is the portion of the current $T$ restricted to the
  simplex $\sigma$.  Consider the subset of ${\cal B}_{\sigma}$
  defined as
  \[ H_T = \left\{ \va \in {\cal B}_{\sigma} \, \big| \, \int_{\sigma}
  J_d \phi(\vx,\va) \, {\rm d}\norm{T}(\vx) > \frac{
  4\upvartheta_{\sigma} \mass(T|_{\sigma})}{\vol_{\ell}({\cal
  B}_{\sigma})} \right\}. \]
  Then $\vol_{\ell}(H_T) \leq (1/3) \vol_{\ell}({\cal
  B}_{\sigma})$. Similarly we define $H_{\boundary T}$ for the
  pushforward of $\boundary T$ and get $\vol_{\ell}(H_{\boundary T})
  \leq (1/3) \vol_{\ell}({\cal B}_{\sigma})$. Then the set $\,{\cal
  B}_{\sigma} \setminus \{H_T \cup H_{\boundary T} \}\,$ defines a
  subset of ${\cal B}_{\sigma}$ with positive measure, with the
  centers $\va$ in this subset satisfying $\, \int_{\sigma} J_d
  \phi(\vx,\va) \, {\rm d} \norm{T}(\vx) \, \leq \, 4
  \upvartheta_{\sigma} \mass(T|_{\sigma})/\vol_{\ell}({\cal
  B}_{\sigma})\,$ and $\, \int_{\sigma} J_d \phi(\vx,\va) \, {\rm d}
  \norm{\boundary T}(\vx) \, \leq \, 4 \upvartheta_{\sigma}
  \mass(T|_{\sigma})/\vol_{\ell}({\cal B}_{\sigma})$. Hence we can
  choose centers to retract from in each simplex $\sigma$ such that
  the expansion of mass of the current restricted to that simplex is
  bounded by $4 \upvartheta_{\sigma} / \vol_{\ell}({\cal
  B}_{\sigma})$. The bound specified in the Lemma follows when we
  consider retracting the entire current over multiple simplices in
  $K$, and set $\upvartheta_K = \max_{\sigma \in K}
  \upvartheta_{\sigma}$ as the generic upper bound that holds for all
  simplices in $K$.
  \end{proof}

  \paragraph{Bound on complete sequence of retractions.} We can apply
  the bound specified in Lemma \ref{lem:retTbnd} over multiple levels
  $\ell$. Pushing $T$ onto the $d$-skeleton of $p$-complex $K$
  multiplies the mass of $T$ by a factor of at most
  $(4\upvartheta_K)^k$. Likewise, pushing $\boundary T$ on to the
  $(d-1)$-skeleton multiplies the mass of $\boundary T$ by a factor of
  at most $(4\upvartheta_K)^{k+1}$.

\subsubsection{Bounding the distance between the current and its simplicial approximation}

  In the final step, we construct the simplicial current $P$
  approximating the original current $T$, and bound the flat norm
  distance between the two. Since we are now considering retraction
  maps over many simplices simultaneously, we let $\phi_i$ denote the
  global projection from the $(p-i+1)-$skeleton to the
  $(p-i)-$skeleton, suppressing the particular $\vx$ and $\va$.  We
  denote the composition of all these steps as $\psi^1 \equiv \phi_k
  \circ \dots \circ \phi_1$ and hence we map $T$ forward by $\psi^1$,
  picking centers (see Lemma \ref{lem:bndJafixx}) to project from in
  each step and in each simplex.  We pick each of these centers such
  that the retractions map $\boundary T$ with bounded amplification of
  mass as well (see Lemma \ref{lem:retTbnd}).

  The homotopy formula \cite[Section 7.4.3]{KrantzParks2008} states
  that given a smooth homotopy $g$ from $f_0$ to $f_1$ where $f_0, f_1
  : U \subseteq \R^{n_0} \to \R^{n_1}$ are smooth functions with
  $g(0,x) = f_0(x)$ and $g(1,x) = f_1(x)$, if $T$ is a $d$-current and
  $f^{-1}(F) \cap \spt f$ is compact for every compact set $F
  \subseteq \R^{n_1}$, we have that the difference in pushforwards of
  $T$ under $f_1$ and $f_0$ is given by
  \[
    {f_1}_{\#}(T) - {f_0}_{\#}(T) = \boundary g_{\#}([0,1] \times T) +
    g_{\#}([0,1] \times \boundary T).
    \]
  Define the homotopy $g(\gamma,\vx) = \gamma\vx + (1-\gamma)
  \psi^1(\vx)$ for $\gamma \in [0,1]$. Then the homotopy formula gives
  \[ T - \psi^1_{\#}(T) = \boundary g_{\#}([0,1] \times T) +
  g_{\#}([0,1] \times \boundary T).\]
  We define $R = g_{\#}([0,1] \times T)$ and $Q_1 = g_{\#}( [0,1]
  \times \boundary T)$.  Then we get
  \begin{equation} \label{eq:Tminuspsi}
    T - \psi^1_{\#}(T) = \boundary R + Q_1. 
  \end{equation}

  \noindent Finally, we map $\psi^1_{\#}(\boundary T)$ forward to the
  $(d-1)$-skeleton of simplicial complex $K$ with $\phi = \phi_{k+1}$
  to get $\psi^2_{\#}(\boundary T) = \phi_{\#}(\psi^1_{\#}( \boundary
  T))$. For this purpose, consider the homotopy $h(\gamma,\vx)$ from
  $\psi^1_{\#}(\boundary T)$ to $\psi^2_{\#}(\boundary T)$, i.e.,
  \[ h(\gamma,x) = \gamma\psi^1_{\#}(\vx) + (1-\gamma)
  \psi^2_{\#}(\vx) ~ \mbox{ for }~\gamma \in [0,1]. \]
  We define
  \begin{equation} \label{eq:defP}
    P = \psi^1_{\#}(T)- h_{\#}([0,1]\times\psi^1_{\#}(\boundary T)).
  \end{equation}
  $P$ is a $d$-current whose boundary $\boundary P$ is contained in the
  $(d-1)$-skeleton of $K$. Define $Q_2 = h_{\#}([0,1] \times
  \psi^1_{\#}(\boundary T))$.  Using the homotopy formula, we get
   \begin{align*} % \label{eq:homotopy}
    \boundary P & = \boundary\left( \, \psi^1_{\#}(T)-
    h_{\#}([0,1]\times\psi^1_{\#}(\boundary T)) \, \right) \\
              & =  \, \psi^1_{\#}(\boundary T)-
              \boundary h_{\#}([0,1]\times\psi^1_{\#}(\boundary T)) \, \\
              & = \psi^2_{\#}(\boundary T) \subset (d-1)\text{-skeleton of } K.
   \end{align*}
  Equation (\ref{eq:defP}) gives $\psi^1_{\#}(T) = P + Q2$.  Defining
  $Q = Q_1 + Q_2$, Equation (\ref{eq:Tminuspsi}) gives
  \begin{align*}
    T - (P+Q2) & = \boundary R + Q_1, ~~\mbox{hence } \\
    T - P      & = \boundary R + Q.  
  \end{align*}

  Finally, we apply the bounds on the retraction described in Lemma
  \ref{lem:retTbnd} and the paragraph following this Lemma to the
  masses of the pushforwards. Noticing that $\diam(\sigma) \leq
  \Updelta$ for all $\sigma \in K$, we get the following bounds, which
  finish the proof of our simplicial deformation theorem (Theorem
  \ref{thm:simpldeform}). 
  \begin{align*}
  \mass(P) 	    & \leq  (4\upvartheta_{K})^k \mass(T) + 
			\Updelta(4\upvartheta_{K})^{k+1}\mass(\boundary T)\\
                    &  =    (4\upvartheta_{K})^k\left(\mass(T) + 
			\Updelta(4\upvartheta_{K})\mass(\boundary T)\right),\\
  \mass(\boundary P) & \leq  (4\upvartheta_{K})^{k+1}\mass(\boundary T),\\
  \mass(R)          & \leq   \Updelta \mass(\psi^1_{\#}(T)) \\
               	    & \leq   \Updelta (4\upvartheta_{K})^{k} \mass(T),~\mbox{ and } \\
  \mass(Q)          & \leq  \Updelta(4\upvartheta_{K})^{k} \mass(\boundary T)
			+ \Updelta(4\upvartheta_{K})^{k+1} \mass(\boundary T)\\
                    &  =    \Updelta(4\upvartheta_{K})^{k}(1+4\upvartheta_{K})
			\mass(\boundary T).
  \end{align*}
\qed

%  \bigskip

  \begin{remark} 
  \label{rem:subdivisionregularity}
    The influence of {\em Simplicial regularity} as measured by
    $\upkappa_1$ and $\upkappa_2$ is clearly revealed by the statement
    of our deformation theorem (Theorem \ref{thm:simpldeform}).
    Explicit constants are a simple yet useful part of the result; as
    observed above in Remark~\ref{rem:defthm1}, the statement of this
    theorem leads to an easy observation that the flat norm distance
    between $T$ and $P$ can be made a small as desired by subdividing
    the simplicial complex in a manner that keeps the regularity
    constants bounded.  This can be done, for example, by using the
    subdivision algorithm of Edelsbrunner and Grayson \cite{EdGr2000}.
  \end{remark}

  \begin{remark}
    We did not explicitly discuss the case of $0$-dimensional
    currents.  In this case, the bounds on mass expansion are all
    equal to one.
  \end{remark}

\subsection{Comparison of Bounds of Approximation} \label{ssec:compbnds}

Sullivan studied the deformation of integral currents on to the
skeleton of a cell complex, which is composed of compact convex
sets. He presented a deformation theorem for deforming integral
currents on to the boundary of a cell complex \cite[Theorem
  4.5]{Sullivan1990}. For ease of comparison, we use {\em our}
notation to restate the bounds given by Sullivan for deforming a
$d$-current $T$ to a polyhedral current $P$ in the boundary of a cell
complex in $\R^q$. Recall that in our simplicial deformation theorem
(Theorem \ref{thm:simpldeform}), the simplicial complex considered has
dimension $p$ and is embedded in $\R^q$ for $q \geq p$.  Furthermore,
$\upkappa_1$, $\upkappa_2$, $\Updelta$, and $\upvartheta_K$ are
simplicial regularity constants.  We also note that even though
Sullivan stated his results for full-dimensional complexes and the
standard flat norm, it is straightforward to extend them to lower
dimensional complexes and the flat norm with scale:
\begin{align}
  \mass(P) & ~\leq~ {p \choose d} \left( 2d \left(\frac{d+1}{2d}
  \upkappa_2 \right)^{d+1} \right)^{p-d+1}
  \mass(T) \label{eq:SullbndMP}, \\
  \vspace*{-0.1in} \nonumber \\
  \mass(\boundary P) & ~\leq~ {p \choose d-1} \left(2d \left(
  \frac{d+1}{2d} \upkappa_2 \right)^{d} \right)^{p-d+1}
  \mass(\boundary T), ~~~ \mbox{and} \label{eq:SullbndMbdyP} \\
  \vspace*{-0.5in} \nonumber \\
 \F_\lambda(T,P) & ~=~\lambda^d \cdot \F_1(T/\lambda,
 P/\lambda)\nonumber\\ &~\leq~ \lambda^d \cdot (p-d+1) \Updelta \, (
 \, \mass(P/\lambda) + \mass(\boundary P/\lambda) \,)\nonumber\\ &~=
 \lambda^d \cdot (p-d+1) \Updelta \, ( \,\lambda^{-d} \cdot \mass(P) +
 \lambda^{1-d} \cdot \mass(\boundary P) \,)\nonumber\\ &~= (p-d+1)
 \Updelta \, ( \, \mass(P) + \lambda \mass(\boundary P) \,).
\label{eq:SullbndFPbdyP}
\end{align}
Our results corresponding to the first two bounds in Equations
(\ref{eq:SullbndMP}) and (\ref{eq:SullbndMbdyP}) are presented in
Equations (\ref{eq:simpdefthmMP}) and (\ref{eq:simpdefthmMbdyP}) in
Theorem \ref{thm:simpldeform}, which we repeat here with the
substitution $k = p-d$.
  \begin{align}
    \mass(P) & \leq (4\upvartheta_K)^{p-d} \mass(T)+ \Updelta(4\upvartheta_{K})^{p-d+1}\mass(\boundary T),
    \tag{\ref{eq:simpdefthmMP} revisited}\\
    \mass(\boundary P) & \leq (4\upvartheta_K)^{p-d+1} \mass(\boundary T), 
    \tag{\ref{eq:simpdefthmMbdyP} revisited}\\
    \mass(R) & \leq \Updelta(4\upvartheta_K)^{p-d} \mass(T), \mbox{ and}
    \tag{\ref{eq:simpdefthmMR} revisited} \\
    \mass(Q) & \leq  \Updelta(4\upvartheta_{K})^{p-d}(1+4\upvartheta_{K})\mass(\boundary T). 
    \tag{\ref{eq:simpdefthmMQ} revisited}
  \end{align}
To obtain the flat
norm distance corresponding to the third bound given by Sullivan in
Equation (\ref{eq:SullbndFPbdyP}), we use the definition of flat norm
distance between two currents specified in Equation
\eqref{eq:flatdist}. Using $T-P = \boundary Q + R$, we combine two of
our bounds specified in Equations \eqref{eq:simpdefthmMR} and
\eqref{eq:simpdefthmMQ} to get
\begin{align*}
\F_\lambda(T,P) & ~\leq~ \Updelta (4\upvartheta_{K})^{p-d} \, \left( \, \mass(T) + \lambda(1 
	+ 4 \upvartheta_{K}) \mass(\boundary T) \, \right). \hspace*{0.5in}
\end{align*}

To gain a better understanding of how the two sets of bounds compare,
we compute these bounds explicitly for the case of a $2$-current in a
regular tetrahedral complex (thus, $p = 3$ and $d = 2$). Notice that
this instance is close to a best case for Sullivan's bounds, as less
regular complexes affect them more severely. With this point in mind,
we present in Table \ref{table:boundcomparison} our bounds and
Sullivan's bounds on both a regular tetrahedral complex and one on
which we stretch the regular tetrahedra by a factor of 10 in a
direction normal to one of their faces (i.e., turn them into skinny,
spike-like simplices). 
\begin{table}[h]
\centering
\renewcommand{\arraystretch}{2}
\begin{tabular}{@{}rccr@{}}\toprule
\\[-7ex]
%& \multicolumn{2}{|c|}{Regular tetrahedra} & \multicolumn{2}{|c|}{Stretched tetrahedra}\\
Quantity & Sullivan's bound & Our bound\\
\hline

\multicolumn{3}{c}{Regular tetrahedra}\\[-1ex]

$M(P)$ & $(1.2 \times 10^5)\, \mass(T)$ & \twolinecelleq{(1.6 \times 10^3) \, \mass(T)}{ + (2.5 \times 10^6)\, \Updelta\mass(\boundary T)}
\\

$M(\boundary P)$ & $(8.7 \times 10^3)\,\mass(\boundary T)$ & $(2.5 \times 10^6)\,\mass(\boundary T)$ \\

$\F_\lambda(T, P)$ & \twolinecelleq{(2.4 \times 10^5)\,\Updelta\mass(T)}{+(1.7 \times
10^3)\,\Updelta\lambda\mass(\boundary T)} & \twolinecelleq{(1.6 \times
10^3)\,\Updelta\mass(T)}{ + (2.5 \times 10^6) \,
\Updelta\lambda\mass(\boundary T)} \vspace*{0.1in} \\[-2ex]

\multicolumn{3}{c}{Stretched tetrahedra}\\[-1ex]

$M(P)$ & $(5.5 \times 10^9)\,\mass(T)$ & \twolinecelleq{(3.7 \times 10^4)\,\mass(T) }{+ (1.4 \times 10^9)\, \Updelta\mass(\boundary T)} \\

$M(\boundary P)$ & $(1.1 \times 10^7)\,\mass(\boundary T)$ & $(1.4 \times 10^9)\,\mass(\boundary T)$ \\

$\F_\lambda(T, P)$ & 
\twolinecelleq{(1.1 \times 10^{10})\Updelta\mass(T)}{ + (2.3 \times 10^7)\Updelta\lambda\mass(\boundary T)}
%\parbox{1.8in}
%{
%\begin{align*}
%&(1.1 \times 10^{10})\Updelta\mass(T)\\
%+\,&(2.3 \times 10^7)\Updelta\lambda\mass(\boundary T)
%\end{align*}
%}
%$(1.09778 \times 10^{10})\Updelta\mass(T) + (2.26170 \times 10^7)\Updelta\lambda\mass(\boundary T)$
& \twolinecelleq{(3.7 \times 10^4)\,\Updelta\mass(T)}{ + (1.4 \times 10^9)\Updelta\lambda\mass(\boundary T)}\\
\bottomrule
\end{tabular}
\caption[Comparison of simplicial deformation theorem bounds with Sullivan]{Comparison of our bounds with those obtained by Sullivan for
  a 2-current in a (1) 3-complex of congruent regular tetrahedra and
  (2) a 3-complex of congruent stretched tetrahedra which are created
  by taking regular tetrahedra and multiplying their height by a
  factor of 10.}
\label{table:boundcomparison}
\end{table}

For the regular tetrahedral complex and the $\mass(P)$ bound, our
coefficient of $\mass(T)$ is more than $74$ times better, but we do
have a second term that can be quite large, but diminishes in
importance if the complex is subdivided appropriately (see Remark
\ref{rem:subdivisionregularity}).  In the stretched complex, our
coefficient on $\mass(T)$ is $1.5 \times 10^5$ times better,
indicating that our bound is better behaved for irregular complexes.
Our bound on $\mass(\boundary P)$ is about $290$ times worse than
Sullivan's for the regular tetrahedra, and about $120$ times worse for
the stretched complex.  For the flat norm bound in the regular
complex, we are about $148$ times better on the $\mass(T)$ term and
about $145$ times worse on the $\mass(\boundary T)$ term.  On the
stretched complex, our $\mass(T)$ coefficient is about $3 \times 10^5$
times better, and our $\mass(\boundary T)$ coefficient is about $60$
times worse.  We also note that in the case of the flat norm with
scale, our larger $\mass(\boundary T)$ coefficient becomes less
important for small $\lambda$.

\begin{remark}
  For the important case where $\boundary T$ is empty, i.e., when $T$
  is a cycle, we have $\mass(\boundary T) = 0$, and hence our bounds
  are uniformly better than Sullivan's.
\end{remark}
As compared to Sullivan, we are able to take advantage of our
simplicial setting to get better bounds on the mass expansion of $T$.
While our mass expansion bounds involving $\boundary T$ are currently
inferior to Sullivan's, we suspect our arguments can be tightened and
modified to obtain bounds that are better in all cases.  More
importantly, our bounds are less sensitive to simplicial irregularity.
Given the challenges inherent in creating meshes without slivers even
in three dimensions \cite{ChDeSh2012}, bounds that behave well in
their presence are highly desirable.

\section{Computational Results} \label{sec:compresults}

\begin{figure}[htp!]
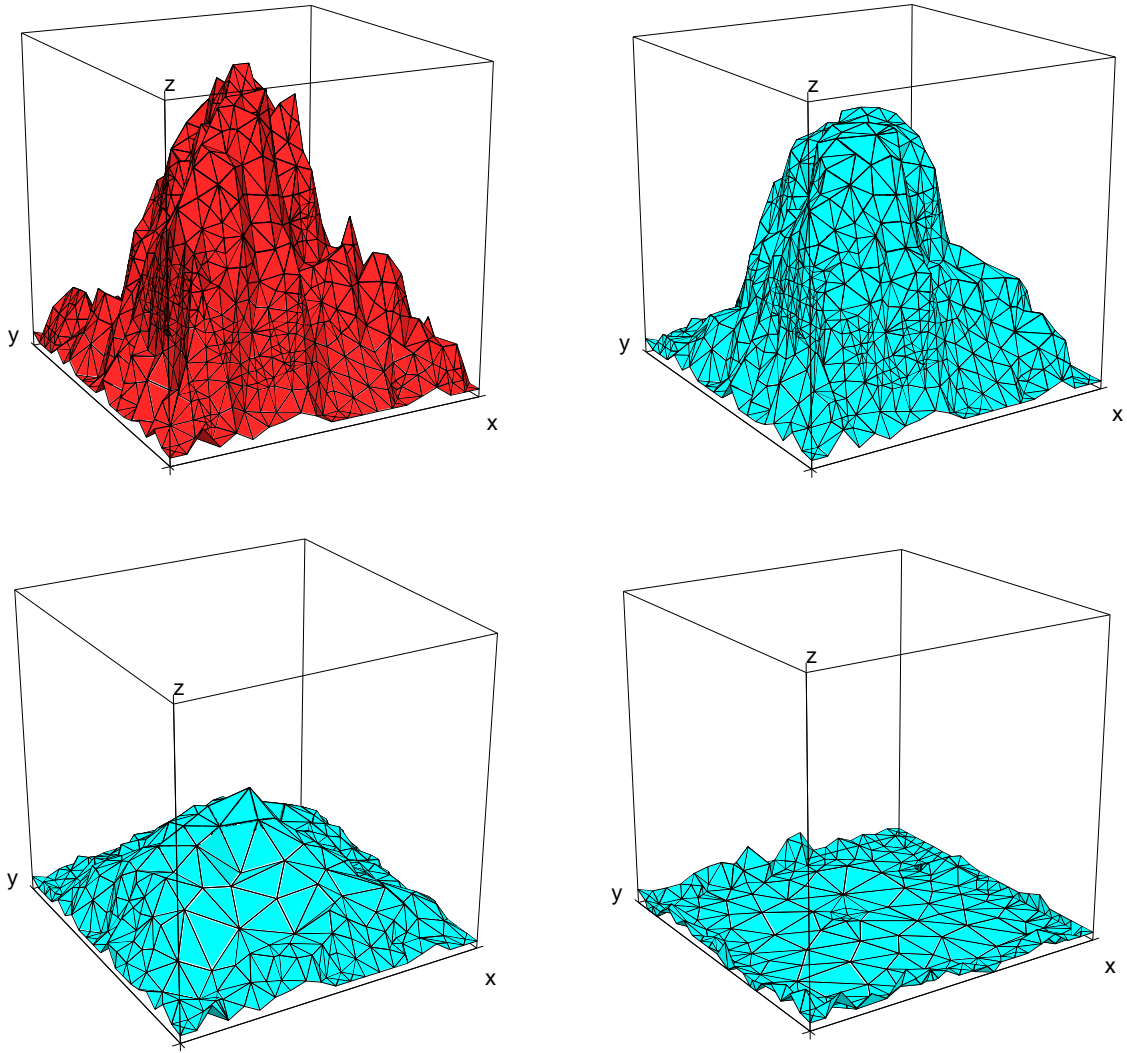

\centering 
\vspace*{0.5in}
 \includegraphics[scale=0.51, trim=1.2in 4.1in 1.5in 1.5in, clip]{../Figs/View1_Pyramid.pdf} 
%\includegraphics[scale=0.6, trim=1.1in 1.1in 1in 1.5in, clip]{View1_Pyramid.eps} 
%\hspace*{0.05in}
\hfill
 \includegraphics[scale=0.51, trim=1.2in 4.1in 1.5in 1.5in, clip]{../Figs/MSFNSurf_l6.pdf} \\
 \includegraphics[scale=0.51, trim=1.2in 4.1in 1.5in 1.2in, clip]{../Figs/MSFNSurf_l2.pdf} 
%\includegraphics[scale=0.55, trim=1.1in 1in 1in 1.4in, clip]{MSFNSurf_l2.eps} 
%\hspace*{0.05in}
\hfill
 \includegraphics[scale=0.51, trim=1.2in 4.1in 1.5in 1.2in, clip]{../Figs/MSFNSurf_l0p01.pdf} 
\vspace*{0.1in}
\caption[3-dimensional flat norm decomposition example]{\label{fig:pyrsurf} Top left: A view of original 
pyramidal surface in three dimensions. The remaining three figures
show the flat norm decomposition for scales $\lambda=6$ (top right),
$\lambda=2$ (bottom left), and $\lambda=0.01$ (bottom right). See text
for further explanation. The images were generated using the package
TetView \cite{tetview}.}
\end{figure}

We illustrate computations of the \MSFN{} by describing the flat norm
decompositions of a $2$-manifold with boundary embedded in $\R^3$ (see
Figure \ref{fig:pyrsurf}). The input set has the underlying shape of a
pyramid, to which several peaks and troughs of varying scale, as well
as random noise, have been added. We model this set as a piecewise
linear $2$-manifold with boundary, and find a triangulation of the
same as a subcomplex of a tetrahedralization of the $2 \times 2 \times
2$ cube centered at the origin, within which the set is located. We
use the method of constrained Delaunay tetrahedralization
\cite{Si2010} implemented in the package TetGen \cite{tetgen} for this
purpose. We then compute the \MSFN{} decomposition of the input set at
various scale ($\lambda$) values. At high values, e.g., when
$\lambda=6$, the optimal decomposition resembles the input set with
the small kinks due to random noise smoothed out. At the other end,
for $\lambda=0.01$, the optimal decomposition resembles a flat
``sheet''. For intermediate values of $\lambda$, the optimal
decomposition captures features of the input set at varying scales.

The entire $3$-complex mesh modeling the cube in question consisted of
14,002 tetrahedra and 28,844 triangles. For each $\lambda$,
computation of the \MSFN{} described above took only a few minutes on a
regular PC using standard functions from MATLAB. This example
demonstrates the feasibility of efficiently computing flat norm
decompositions of large datasets in high dimensions, for the purposes
of denoising or to recover scale information of the data.

\section{Discussion}

Our result on simplicial deformation (Theorem \ref{thm:simpldeform})
places the definition of the \MSFN{} into clear context. If a current
lives in the underlying space of a simplicial complex, we can deform
it to be a simplicial current on the simplicial complex, and do so
with {\em controlled} error. In fact, by subdividing the simplicial
complex carefully, we can move this error as close to zero as we
like. Since the \MSFN{} could be computed efficiently when the
simplicial complex does not have relative torsion, one could naturally
use our approach to compute the flat norm of a large majority of
currents in arbitrarily large dimensions. An important open question
in this context is whether the \MSFN{} of a current on a simplicial
complex {\em with} relative torsion could be approximated efficiently
by {\em coarsening} the complex so that the relative torsion is
removed. For instance, it has been observed recently that edge
contractions could remove existing relative torsion while preserving
the homology groups of the simplicial complex in certain cases
\cite{DeHiKrSm2013}.

The \MSFN{} problem, similar to the recent results on the \OBCP{}
\cite{DuHi2011}, apply notions from algebraic topology and discrete
optimization to problems from geometric measure theory such as flat
norm of currents and area-minimizing hypersurfaces. What other classes
of problems from the broader area of geometric analysis could we
tackle using similar approaches? One such question appears to be the
following: under what conditions is the flat norm decomposition of an
integral current guaranteed to be another integral current?  Working
in the setting of simplicial complexes, results on the existence of
integral optimal solutions for instances of ILPs with integer
right-hand side vectors may prove useful in answering this question.

While $L^1$TV and flat norm computations have been used widely on data
in two dimensions, such as images, the \MSFN{} opens up the
possibility of utilizing flat norm computations for higher dimensional
data. Similar to the flat norm-based signatures for distinguishing
shapes in two dimensions \cite{Vietal2010}, could we define shape
signatures using \MSFN{} computations to characterize the geometry of
sets in arbitrary dimensions? The sequence of optimal \MSFN{}
decompositions of a given set for varying values of the scale
parameter $\lambda$ captures all the scale information of its
geometry. Could we represent all this information in a compact manner,
for instance, in the form of a barcode?

\section*{Acknowledgments}

We acknowledge the financial support from the National Science
Foundation (NSF) through grants DMS-0914809 and CCF-1064600.

\chapter{Flat norm decomposition of integral currents\footnote{Based on \cite{IbKrVi2014}}}
\label{ch:ic}

% margins for LaTeX-dvipdf/ps
%\setlength{\oddsidemargin}{0in} 
%\setlength{\evensidemargin}{0in}
% \setlength{\topmargin}{-0.5in} %submitted to arXiv
%\setlength{\topmargin}{0.2in}
%\setlength{\textwidth}{6.5in} 
%\setlength{\textheight}{9in}

% margins for PdfLaTeX
%\usepackage[margin=1in]{geometry}

% for PdfLaTeX
%\usepackage[small,bf]{caption}
%\usepackage[colorlinks,linkcolor=black,bookmarksopen,
%bookmarksnumbered,citecolor=black,urlcolor=black]{hyperref}

% Definitions from Bala
\newcommand{\currentrestr}{%
  \,\raisebox{-.127ex}{\reflectbox{\rotatebox[origin=br]{-90}{$\lnot$}}}\,%
}

\newenvironment{subproof}{\noindent{\it Proof of claim.}}{}
% Theorem-like declarations

%% Start of paper
%% --------------

\title{{Flat norm decomposition of integral currents}}

%\author{Sharif~Ibrahim,\thanks{\affil{Department of Mathematics, Washington State University, Pullman, WA 99164-3113},\newline
%\email{\{math.msfn,bkrishna,vixie\}@\{sharifibrahim.com,math.wsu.edu,speakeasy.net\}}}\,\,
%  Bala~Krishnamoorthy,\footnotemark[1]\,\,\footnote{Corresponding author}\,\,
%  Kevin R.~Vixie\footnotemark[1]
%}

\section{Introduction}
In geometric measure theory, currents represent a generalization of oriented surfaces with multiplicities.
Currents were developed in the context of Plateau's problem and have also found application in isoperimetric problems and soap bubble conjectures\cite{Morgan2008}.

Given a $d$-dimensional current $T$, we can consider decompositions $T = X + \boundary S$ where $X$ is a $d$-dimensional current and $S$ is a $(d+1)$-dimensional current.
Over all such decompositions, the minimum total mass (volume) of the two pieces (i.e., $\mass(X) + \mass(S)$) is the flat norm $\F(T)$.
More recently, the $L^1$TV functional (introduced in the form most relevant to us by Chan and Esedo\=glu\cite{ChEs2005}) was shown to be related to the flat norm\cite{MoVi2007}.
This connection suggested the flat norm with scale (yielding the objective $\mass(X) + \lambda \mass(S)$ for any fixed scale $\lambda$) and a geometric interpretation for the optimal decompositions: varying $\lambda$ controls the scale of features isolated in the decomposition.

One natural question: must currents in a particular regularity class (in this paper, integral currents) have an optimal flat norm decomposition in the same class?
The $L^1$TV connection shows this is true for boundaries of codimension 1 (i.e., boundaries of $(d+1)$-currents in $\R^{d+1}$) since the $L^1$TV functional applied to binary (or step function) input is known to have binary (step function) minimizers\cite{ChEs2005}.
This may be taken one step further in the discretized problem where the boundary requirement can be dropped\cite{IbKrVi2013}.

In the present work, we present a framework to bridge the gap between the continuous and discrete cases, assuming a suitable triangulation result.
This allows us to drop the requirement that integral $d$-currents in $\R^{d+1}$ be boundaries to have a guaranteed integral optimal decomposition.
The necessary triangulation result is proved in $\R^2$ by means of Shewchuk's Terminator algorithm\cite{Sh2002} for subdividing planar straight line graphs.
This algorithm simultaneously bounds the smallest angles in the complex and tells us where they can occur, allowing us to tailor a simplicial complex to a given set of input currents.
We then obtain a simplicial deformation theorem with constant bounds for these currents and simplicial complex, ensuring the sequence of aprroximating discretized problems are well-behaved and solve the continuous problem in the limit.
Assuming a suitable triangulation result for higher dimensions (see \cref{con:boundreg}), we show that codimension 1 integral currents have an integral optimal flat norm decomposition.

For the related problem of least area with a given boundary (which can be considered as the flat norm problem with $X$ constrained to be empty), counterexamples of Young\cite{Young1963}, White\cite{White1984}, and Morgan\cite{Morgan1984} provide instances in which the minimizer is not integral.
These negative results are of codimension 3 (i.e., 1-dimensional curves in $\R^4$) which may translate into a limit on the flat norm question.

\subsection{Definitions}

To formally define $d$-currents in $\R^n$, let $\form{d}$ be the set of $C^\infty$ differentiable $d$-forms with compact support.
The set of $d$-currents (denoted $\current{d}$) is the dual space of $\form{d}$ with the weak topology.

Currents have mass and boundary that correspond (for rectifiable currents, at least) to one's intuition for what these should mean for $d$-dimensional surfaces in $\R^n$ with care taken to respect orientation and multiplicities.
For more general classes of current, these concepts are still defined but may not have the same geometric significance.
The mass of a $d$-current $T$ is formally given by $\sup_{\phi \in \form{d}} \{T(\phi) \mid \norm{\phi} \leq 1\}$ and the boundary is defined when $d \geq 1$ by $\boundary T(\psi) = T(\dif \psi)$ for all $\psi \in \form{d-1}$.
When $T$ is a 0-current, we let $\boundary T = 0$ as a 0-current.
%This last convention is not universal, but allows us to simplify our presentation slightly.
The boundary operator on currents is linear and nilpotent (i.e., $\boundary \boundary T = 0$ for any current $T$), inheriting these properties from exterior differentiation of forms (which are linear and satisfy $\dif \dif \phi = 0$).

{\it Normal $d$-currents} have compact support and finite mass and boundary mass (i.e., $\mass(T) + \mass(\boundary T) < \infty$).
The set $\mathcal{R}_d$ denotes the {\it rectifiable $d$-currents} and contains all currents with compact support that represent oriented rectifiable sets with integer multiplicities and finite mass.
That is, sets which are almost everywhere the countable union of images of Lipschitz maps from $\R^d$ to $\R^n$.
Lastly, the set $\mathcal{I}_d$ represents {\it integral $d$-currents} and contains all currents that are both rectifiable and normal (formally, it is the set of rectifiable currents with rectifiable boundary, but this definition is equivalent by the closure theorem\cite[4.2.16]{Federer1969}).

\begin{figure}[ht!]  
\centering
\includegraphics[scale=0.7, trim=1.25in 8.3in 2.8in 1in, clip]{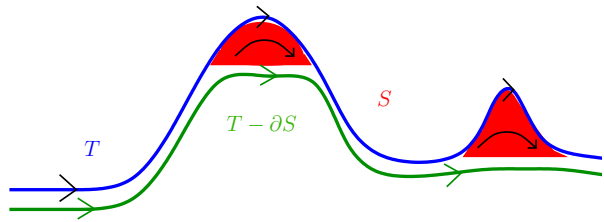}
\caption[Flat norm decomposition]{\label{fig:ic1dcurrent} The flat norm decomposes the 1D current $T$ into (the boundary of) a 2D piece $S$ and the 1D piece $X = T - \boundary S$. The resulting current is shown slightly separated from the input current for clearer visualization.}
\label{fig:icflatnorm}
\end{figure}

\begin{figure}
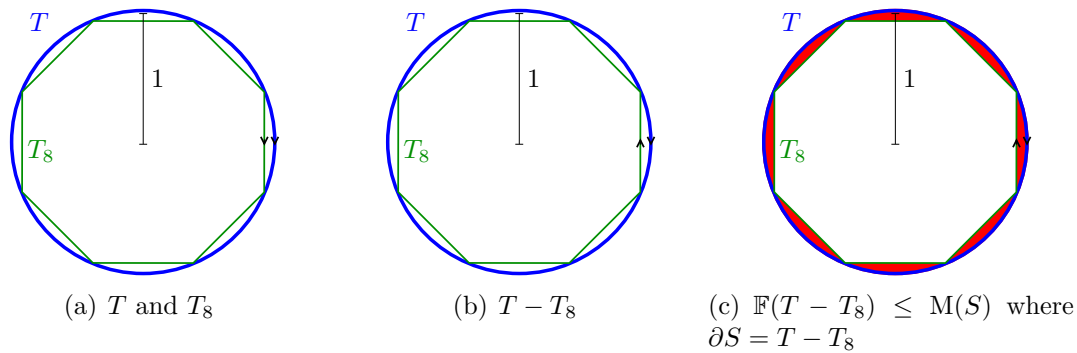

\begin{center}
\subfigure[$T$ and $T_8$]{\label{subfig:icflatdistance1} \makebox[\linewidth/10*3][c]{\input{polygonexampleb1.pspdftex}}}
\subfigure[$T - T_8$]{\label{subfig:icflatdistance2} \makebox[\linewidth/10*3][c]{\input{polygonexampleb2.pspdftex}}}
\subfigure[$\F(T-T_8) \leq \mass(S)$ where $\boundary S = T - T_8$]{\label{subfig:icflatdistance3} \makebox[\linewidth/10*3][c]{\input{polygonexampleb3.pspdftex}}}
%\subfigure[]{\subimport{figs}{polygonexampleb1.pspdftex}}
%\subfigure[]{\input{figs/polygonexampleb2.pspdftex}}
%\subfigure[]{\input{figs/polygonexampleb3.pspdftex}}
\end{center}
\caption[Flat distance example]{The flat norm indicates the unit circle $T$ and inscribed $n$-gon $T_n$ are close because the region they bound has small area.}
\end{figure}

The flat norm of a current $T$ is given by 
\[
\F(T) = \min \{\mass(X) + \mass(S) \mid T = X + \boundary S, X \in \mathcal{E}_d, S \in \mathcal{E}_{d+1}\}
\]
where $\mathcal{E}_d$ is the set of $d$-dimensional currents with compact support (see \cref{fig:icflatnorm}).
The Hahn-Banach theorem guarantees this minimum is attained\cite[p. 367]{Federer1969} so it makes sense to talk about particular $X$ and $S$ as a flat norm decomposition of $T$ (note, however, that the decomposition need not be unique).

For two currents, the flat distance between them is given by $\F(T, P) = \F(T-P)$.
This definition is useful because it is robust to small additions and perturbances (e.g., noise) and reflects when currents are intuitively close.
For example, given a current $T$ representing a unit circle in $\R^2$ and an inscribed $n$-gon $T_n$ (both oriented clockwise, see \cref{subfig:icflatdistance1}), one would like $T_n$ to converge to $T$ in some sense as $n \to \infty$ which the flat norm accomplishes (contrast with the mass norm $\mass(T_n - T)\rightarrow 4\pi$).

The flat norm can be usefully discretized as well.
Given a simplicial $(d+1)$-complex $K$ and a $d$-chain $T$ on $K$, the simplicial flat norm\cite{IbKrVi2013} of $T$ on $K$ is denoted by $\F_K(T)$ and defined analogously except that $X$ and $S$ are restricted to be chains on $K$.
%\[
%\F_K^\lambda(T) = \min_{\vs \in \Z^n} \cbr{\sum_{i=1}^m \operatorname{V}_d(\sigma_i)\abs{x_i} + \lambda \sum_{j=1}^n \operatorname{V}_{d+1}(\tau_j)\abs{s_j} \mid \vx = \vt - [\boundary_{d+1}]\vs, \vx \in \Z^m}
%\]

\subsection{Overview}
Our general technique is a standard notion: express the continuous problem as a limit of discrete problems for which the result holds.
\cref{thm:simpint} tells us that the {\em simplicial} flat norm of an integral chain in codimension 1 has an optimal integral current decomposition; by the compactness theorem from geometric measure theory, the limit of these decompositions is also integral.

\begin{figure}
\centering
\begin{tikzpicture}[scale=1]
    \tikzstyle{point}=[circle,thick,draw=black,fill=black,inner sep=0pt,minimum width=4pt,minimum height=4pt]
    \tikzstyle{current}=[thick]
%    \node (a)[point] at (0,0) {};
%    \node (b)[point] at (3,0) {};
%    \node (c)[point] at (2,2) {};

%    \begin{scope}[yshift=2cm]
%    \node (d)[point] at (1,1) {};
%    \node (e)[point] at (0,2) {};
%    \node (f)[point] at (4,2) {};
%    \end{scope}
    
    \node (p)[point,label={[label distance=0cm]90:$A$}] at (0,0) {};
    \node (K2)[below right = 1cm and 3.5cm,label={[label distance=0cm]90:$K_2$}] {};

    \node[above right = 1cm and 1.73cm] (t1)[point] {};
    \node[below = 2cm of t1] (t2)[point] {};
    \node[above right = 1cm and 1.73cm of t2] (t3)[point] {};
    \node[above right = 1cm and 1.73cm of t3] (t4)[point] {};
    \node[below = 2cm of t4] (t5)[point] {};
    \node[above right = 1cm and 1.73cm of t5] (q)[point,label={[label distance=0cm]5:$B$}] {};

     \draw (p.center) -- (t1.center) -- (t2.center) -- cycle;
     \draw (t1.center) -- (t2.center) -- (t3.center) -- cycle;
     \draw (t4.center) -- (t5.center) -- (t3.center) -- cycle;
     \draw (t4.center) -- (t5.center) -- (q.center) -- cycle;
     \draw[current,color=red,dashed] (p.center) -- node[label={[label distance=1cm]176:$T$}] {} (q.center);
     \draw[current,color=blue] (p.center) -- (t1.center) -- (t3.center) -- node[label={[label distance=0cm]88:$P_2$}] {} (t4.center) -- (q.center);

%    \draw[pattern=north east lines] (a.center) -- (p.center) -- (b.center) -- cycle;
%    \draw[pattern=north west lines] (a.center) -- (p.center) -- (c.center) -- cycle;
%    \draw[pattern=vertical lines]   (b.center) -- (p.center) -- (c.center) -- cycle;
%    \draw[pattern=dots] (d.center) -- (e.center) -- (f.center) -- cycle;
%    \draw (p.center) -- (d.center);
\end{tikzpicture}
\caption[Simplicial flat norm need not converge to continuous flat norm]{A sequence of simplicial chains that converges in the flat norm (i.e., $P_n \rightarrow T$) need not have convergent simplicial flat norm values (i.e., $\F_{K_n}(P_n) \rightarrow \F(T)$ need not hold).
The current $T$ is the segment from $A$ to $B$, the complex $K_n$ is the arrangement of $2n$ equilateral triangles of appropriate size stretching from $A$ to $B$ and $P_n$ is the top chain from $A$ to $B$ on $K_n$.
Clearly, $\F(T-P_n) \rightarrow 0$ but $\F_{K_n}(P_n) = \frac{2}{\sqrt{3}}\F(T) \not\rightarrow \F(T)$.}
\label{fig:flatnormconv}
\end{figure}
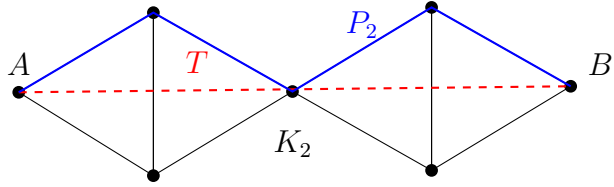

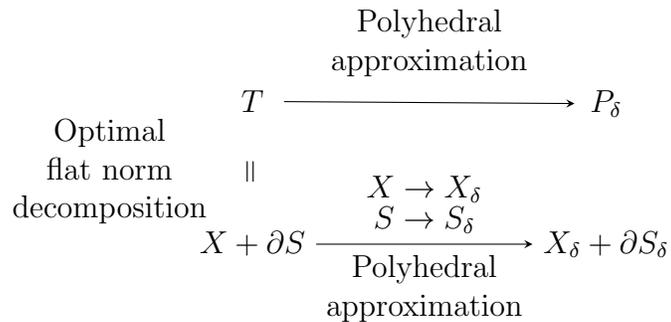
\begin{figure}
\centering
%\begin{tikzpicture}[scale=1,every node/.style={scale=1}]
%  \matrix (m) [matrix of math nodes,row sep=3em,column sep=3em,minimum width=2em]
%  {
%    T & P_\delta\\ % Add & Y_\delta + R_\delta to add simplicial decomposition of P
%    X+\boundary S & X_\delta + \boundary S_\delta\\
%  };
%  \path[-stealth]
%   (m-1-1) edge[-,draw opacity=0] node[rotate=90] {$=$} (m-2-1)
%            edge (m-1-2)
%   % Uncomment to fill in relationship between P and its simplicial decomposition (if added above)
%   %(m-1-2) edge[-,draw opacity=0] node {$=$} (m-1-3)
%   (m-2-1.east|-m-2-2) edge (m-2-2);
%\end{tikzpicture}
\begin{tikzpicture}[scale=1,every node/.style={scale=1}]
  \matrix (m) [matrix of math nodes,row sep=3em,column sep=7em,minimum width=2em]
  {
    T & P_\delta \\%&Y_\delta + \boundary R_\delta\\
    X+\boundary S & X_\delta + \boundary S_\delta\\
  };
  \path[-stealth]
   (m-1-1) edge[-,draw opacity=0] node[rotate=90] {$=$} node[left=1em,text width=6.5em,align=center] {Optimal\\flat norm\\decomposition}(m-2-1)
            edge node[above=0.5em,text width=7em,align=center] {Polyhedral\\approximation}(m-1-2)
   %(m-1-2) edge[-,draw opacity=0] node {$=$} node[above=0.5em,text width=7em,align=center] {Optimal simplicial\\flat norm\\decomposition} (m-1-3)
   (m-2-1.east|-m-2-2) edge node [above=1em] {$X \to X_\delta$}
            node [above] {$S \to S_\delta$} node [below,text width=7em,align=center] {Polyhedral\\approximation} (m-2-2);
\end{tikzpicture}
\label{fig:approxoverview}
\caption[Overview of useful approximations and decompositions]{Various approximations and decompositions used in our results.}
\end{figure}

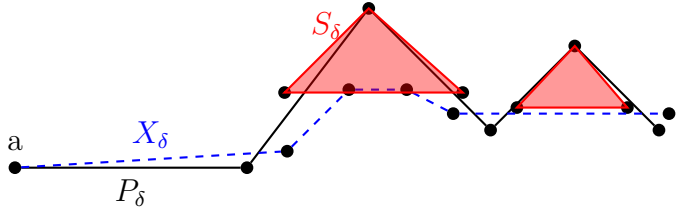
\begin{figure}
\centering
\begin{tikzpicture}[scale=0.5]
	\tikzstyle{point}=[circle,thick,draw=black,fill=black,inner sep=0pt,minimum width=4pt,minimum height=4pt]
	\tikzstyle{current}=[thick]

	\node (t1)[point,label={a}] at (0,0) {};
	\node[right = 3cm] (t2)[point] {};
	\node[above right = 2cm and 1.5cm of t2] (t3)[point] {};
	\node[below right = 1.5cm and 1.5cm of t3] (t4)[point] {};
	\node[above right = 1cm and 1cm of t4] (t5)[point] {};
	\node[below right = 1cm and 1cm of t5] (t6)[point] {};

	\node[above right=0.1cm and 3.5cm of t1] (u2)[point] {};
	\node[above right=0.7cm and 0.7cm of u2] (u3a)[point] {};
	\node[right=0.6cm of u3a] (u3)[point] {};
	\node[below right=0.2cm and 0.5cm of u3] (u4)[point] {};
	\node[right=2.7cm of u4] (u5)[point] {};

	\node[below left=1cm and 1cm of t3] (s2)[point] {};
	\node[right=2.2cm of s2] (s3)[point] {};
	\node[below left=0.7cm and 0.65cm of t5] (s4)[point] {};
	\node[right=1.3cm of s4][point] (s5) {};

    \fill[opacity=0.5,color=red]  (t3) -- (s2) -- (s3) -- cycle;
     %\draw[current,color=red,dashed] (p.center) -- node[label={[label distance=1cm]176:$T$}] {} (q.center);
	\draw[current] (t1) -- node[label={270:$P_\delta$}] {} (t2) -- (t3) -- (t4) -- (t5) -- (t6);
	\draw[current,color=blue,dashed] (t1) -- node[label={$X_\delta$}] {} (u2) -- (u3a) -- (u3) -- (u4) -- (u5);
    \draw[thick,color=red,fill=red,fill opacity=0.4,text opacity=1]  (t3.center) -- node[label={$S_\delta$}] {} (s2.center) -- (s3.center) -- cycle;
    \draw[thick,color=red,fill=red,fill opacity=0.4,text opacity=1]  (t5.center) -- (s4.center) -- (s5.center) -- cycle;
	%\node (b)[point,label={b}] at (5,0) {};
\end{tikzpicture}
\caption[Polyhedral approximation of \cref{fig:icflatnorm}]{A possible polyhedral approximation of the decomposition shown in \cref{fig:icflatnorm}.
Note that $P_\delta \neq X_\delta + \boundary S_\delta$.}
\label{fig:polyhedralization}
\end{figure}

In order to show that an integral current $T$ has integral flat norm decomposition, we therefore find suitable simplicial approximations to $T$ and take the limit of their simplicial flat norm decompositions to obtain an integral decomposition for $T$.

We must also show that this decomposition achieves the flat norm value for $T$ (that is, express $T$ using integral currents in such a way that it remains an optimal flat norm decomposition).
This is immediate if our simplicial approximations to $T$ have simplicial flat norm values that converge to the flat norm of $T$ but this is not necessary (see \cref{fig:flatnormconv}).
We wish to show 
\begin{align}
\label{eq:fteqgoal}
\lim_{\delta \downarrow 0}\F_{K_\delta}(P_\delta) &= \F(T)
\end{align}
 where $P_\delta$ is a simplicial approximation to $T$ on some complex $K_\delta$ with $\F(P_\delta - T) < \delta$.

This goal prevents us from simply using the simplicial deformation theorem to obtain $P_\delta$ since we may end up with the situation in \cref{fig:flatnormconv}.
Instead, we use a polyhedral approximation to $T$ which guarantees that the mass increases by at most $\delta$ (i.e., $\mass(P_\delta) < \mass(T) + \delta$ rather than the simplicial deformation theorem bound $\mass(P_\delta) < C_1\mass(T) + C_2\mass(\boundary T)$ with constants bounded away from 1).

The next step is to take an optimal (possibly nonintegral) decomposition of $T$ and approximate it with polyhedral chains (see \cref{fig:approxoverview}).
That is, approximate the decomposition $T = X + \boundary S$ with polyhedral $X_\delta$ and $S_\delta$.
If these approximations naturally form a decomposition (not necessarily optimal) of $P_\delta$ (i.e., $P_\delta = X_\delta + \boundary S_\delta$), then we would have $\F_{K_\delta}(P_\delta) \leq \mass(X_\delta) + \mass(S_\delta) < \F(T) + 2\delta$ for any complex $K_\delta$ containing $P_\delta$, $X_\delta$, and $S_\delta$.
This of course implies \cref{eq:fteqgoal}.
%the desired result $\lim_{n\rightarrow \infty} \F_{K_\delta}(P_\delta) = \F(T)$.
% from which the desired limit result follows.

However (as in \cref{fig:polyhedralization}), we need not have $P_\delta = X_\delta + \boundary S_\delta$.
Since we obtained these quantities by polyhedral approximation, it turns out that the extent to which this equation is violated is small (in the continuous flat norm).
That is, we have 
\begin{align}
\label{eq:icgoalrestate}
P_\delta = X_\delta + \boundary S_\delta + (P_\delta - T) + (\boundary S - \boundary S_\delta) + (X - X_\delta).
\end{align}

While \cref{eq:icgoalrestate} can be viewed as a decomposition of $P_\delta$, the added error terms mean it may not be a chain on a simplicial complex.
This means it cannot be used directly to bound the simplicial flat norm of $P_\delta$.

If we use the simplicial deformation theorem to push the error terms to some complex $K_\delta$ while preserving a pushed version of \cref{eq:icgoalrestate}, we can obtain a candidate simplicial decomposition of $P_\delta$.
In order to use this to bound $\F_{K_\delta}$, we must know that the deformation theorem didn't make the small error terms large enough to matter.
Unfortunately, the simplicial deformation theorem mass bounds rely on simplicial regularity so a sufficiently skinny simplex could mean the error terms become large.
If the simplicial irregularity in $K_\delta$ gets worse as $\delta \downarrow 0$, we will not be able to show \cref{eq:fteqgoal}.

Since we know exactly which currents we wish to push, the solution is to pick $K_\delta$ with these in mind: make sure the complex is as regular as possible overall (independently of $\delta$) with any irregularities (which may be required to embed $P_\delta$, $X_\delta$, and $S_\delta$) isolated in subcomplexes of small measure.
By making the irregular portions small enough (so they contain a negligible portion of the error terms, even considering the possible magnification from pushing), we establish a deformation theorem variant (\cref{thm:2dboundedsdt}) with constant mass expansion bounds, assuming a triangulation result that lets us isolate the irregularities as described (Shewchuk's Terminator algorithm\cite{Sh2002} provides this in $\R^2$).
The pushed version of \cref{eq:icgoalrestate} allows us to prove $\F_{K_\delta}(P_\delta) \leq \F(T) + O(\delta)$ from which \cref{eq:fteqgoal} and \cref{thm:icmainresult} follow.

\section{Preliminaries}
Our goal is to investigate conditions under which the flat norm decomposition of an integral current can be taken to be integral as well.
The corresponding statement for normal currents is true and useful in our development.

\begin{lemma}
\label{lem:normaldecomp}
If $T$ is a normal $m$-current and $X$ and $S$ are $m$- and $(m+1)$-currents such that $T = X + \boundary S$ and $\F(T) = \mass(X) + \mass(S)$ (i.e., $T = X + \boundary S$ is a flat norm decomposition of $T$), then $X$ and $S$ are normal currents.
\end{lemma}
\begin{proof}
By the definition of normal current, we have $\mass(T) + \mass(\boundary T) < \infty$.
Thus
\[
\mass(X) + \mass(S) = \F(T) \leq \mass(T) < \infty
\]
so $\mass(X) < \infty$ and $\mass(S) < \infty$.
Since $T = X + \boundary S$, we obtain
\[
\mass(\boundary X) = \mass(\boundary\del{X + \boundary S}) = \mass(\boundary T) < \infty.
\]
Lastly,
\[
\mass(\boundary S) \leq \mass(\boundary S - T) + \mass(T) = \mass(-X) + \mass(T)  < \infty.
\]
The currents $X$ and $S$ have compact support by the definition of the flat norm.
Thus $X$ and $S$ are normal by definition.
\end{proof}

Convergence in the flat norm is linear and commutes with the boundary operator as the following easy lemma shows.

\begin{lemma}
\label{lem:convprops}
Suppose that $T_n$ and $U_n$ are $m$-currents for $n = 1, 2, \dots$ and $T_n \rightarrow T$ and $U_n \rightarrow U$ in flat norm (i.e., $\F(T_n - T) \rightarrow 0$) for some $m$-currents $T$ and $U$.
%(that is, $\lim_{n\rightarrow \infty} \langle T_n - T, \phi\rangle = 0$ for every differential $m$-form $\phi$ with compact support) and $U_n \rightarrow U$ for some $m$-currents $T$ and $U$.
The following properties hold:
\begin{enumerate*}
\item[(a)] $\alpha T_n + \beta U_n \rightarrow \alpha T + \beta U$ for any constants $\alpha, \beta \in \R$,
\item[(b)] $\boundary T_n \rightarrow \boundary T$,
\end{enumerate*}
\end{lemma}
\begin{proof}
We apply properties of norms to obtain
\begin{align*}
\F((\alpha T_n + \beta U_n) - (\alpha T + \beta U)) &\leq \F(\alpha T_n - \alpha T) + \F(\beta U_n - \beta U)\\
&= \abs{\alpha}\F(T_n - T) + \abs{\beta}\F(U_n - U).
\end{align*}
Letting $n \rightarrow \infty$ yields the linearity result.
Now let $X_n$ and $S_n$ be $m$- and $(m+1)$-currents such that $X_n + \boundary S_n$ is a flat norm decomposition of $T_n - T$ for $n = 1, 2, \dots$, observing that
\begin{align*}
\F(\boundary T_n - \boundary T) &= \F(\boundary (X_n + \boundary S_n))
= \F(\boundary X_n)
\leq \mass(X_n)
\leq \F(T_n - T).
\end{align*}
The boundary result follows in the limit.
\end{proof}

In the case of the simplicial flat norm, an input integral chain is guaranteed an integral chain decomposition whenever the simplicial complex is totally unimodular\cite{IbKrVi2013}.
This occurs when the complex is free of relative torsion which is the case for any $(d+1)$-complex in $\R^{d+1}$ or when triangulating a compact, orientable $(d+1)$-dimensional manifold.

\begin{theorem}[Simplicial flat norm integral decomposition\cite{IbKrVi2013}]
\label{thm:simpint}
If $K$ is a simplicial $(d+1)$-complex embedded in $\R^{d+1}$, then for any integral $d$-chain $P$ on $K$, the optimal simplicial flat norm value for $P$ is attained by an integral decomposition.
\end{theorem}

We state the simplicial deformation theorem and sketch a portion of its proof.
We will later modify it to obtain a multiple current deformation theorem that preserves linearity (\cref{thm:multisdt}).

\begin{theorem}[Simplicial deformation theorem\cite{IbKrVi2013}]
\label{thm:sdt}
Suppose $K$ is a $p$-dimensional simplicial complex in $\R^q$ and $T$ is a normal $d$-current supported on the underlying space of $K$.
There exists a simplicial $d$-current $P$ supported on the $d$-skeleton of $K$ with boundary supported on the $(d-1)$-skeleton (i.e., a simplicial $d$-chain) such that $T-P=Q+\boundary R$ and there exists a constant $\upvartheta_K$ (depending only on simplicial regularity in $K$) such that the following controls on mass hold:
\begin{align}
\mass(P) &\leq (4\upvartheta_K)^{p-d} \mass(T) + \Delta (4\upvartheta_K)^{p-d+1} \mass(\boundary T)\\
\mass(\boundary P) &\leq (4\upvartheta_K)^{p-d+1}\mass(\boundary T)\\
\mass(Q) &\leq \Delta (4\upvartheta_K)^{p-d}(1+4\upvartheta_K)\mass(\boundary T)\\
\mass(R) &\leq \Delta(4\upvartheta_K)^{p-d}\mass(T)\\
\F(T,P) &\leq \Delta (4\upvartheta_K)^{p-d}(\mass(T)+(1+4\upvartheta_K)\mass(\boundary T))
\end{align}
where $\Delta$ is the diameter of the largest simplex in $K$.
The regularity constant $\upvartheta_K$ is given by
\begin{equation}
\label{eq:regularity}
\upvartheta_K = \sup_{\sigma \in K}\frac{\diam(\sigma)\perimeter(\sigma)}{B_\sigma} + 2 \sup_{\sigma \in K}\frac{\diam(\sigma)}{\inrad(\sigma)}
\end{equation}
where for each $l$-simplex $\sigma$, $\perimeter(\sigma)$ is the $(l-1)$-volume of $\boundary \sigma$ and $B_\sigma$ is the $l$-volume of a ball with radius $\inrad(\sigma)/2$ in $\R^l$.
\end{theorem}
\begin{proof}[Proof highlights]
The simplicial current $P$ is obtained by pushing $T$ and its boundary to the $d-$ and $(d-1)$-dimension skeletons of $K$ respectively.
This pushing is done one dimension at a time; that is, $T$ is pushed from the $p$-skeleton (i.e., the full complex $K$) to the $(p-1)$-skeleton, then to the $(p-2)$ and so on until the $d$-skeleton.
Pushing the current from the $\ell$-skeleton to the $(\ell-1)$-skeleton is done by picking a projection center in each $\ell$-simplex $\sigma$ and projecting the current in $\sigma$ outwards to $\boundary \sigma$ via straight-line projection.

A crucial step in the proof is to find a projection center that bounds the expansion of $T$ and $\boundary T$.
In particular, this is done by proving that over all possible centers, the average expansion is bounded and then showing that individual centers exist with bounded expansion.
We call out this particular step because we modify it to obtain the next theorem.

When projecting onto the skeleton of each simplex $\sigma$, we have\cite[Lemma 5.9]{IbKrVi2013}
\begin{equation}
\label{eq:avgbound}
 \int_{{\cal B}_\sigma} \int_{\sigma} J_d \phi(\vx,\va) \, {\rm
  d}\norm{T}(\vx) \, {\rm d} {\cal L}^{\ell}(\va) 
%= \int_{\sigma}
%  \int_{{\cal B}_\sigma} J_d \phi(\vx,\va) \, {\rm d}{\cal
%  L}^{\ell}(\va) \, {\rm d} \norm{T}(\vx) \,
 \leq \,
  \upvartheta_{\sigma} \mass(T|_{\sigma}).
\end{equation}
where ${\cal B}_\sigma$ is the set of possible centers in $\sigma$ and $\upvartheta_{\sigma}$ is a regularity constant for $\sigma$ related to $\upvartheta_K$ by $\upvartheta_K = \sup_{\sigma \in K} \upvartheta_\sigma$.
This shows that in each projection step the mass of $T$ expands by a factor of at most $\upvartheta_{K}$ averaged over all possible choices of centers.
As the average expansion over all centers is $\upvartheta_K$, we observe that at most $\frac{1}{4}$ of the possible centers can expand the mass of $T$ by a factor of $4\upvartheta_{K}$ or more.
Similarly, at most $\frac{1}{4}$ of the centers can expand $\boundary T$ by a factor of $4\upvartheta_{K}$ or more.
Therefore, at least $\frac{1}{2}$ of the possible centers bound the expansion of both $T$ and $\boundary T$ by at most a factor of $4\upvartheta_K$.
Choosing a center from this set for each simplex yields the bounds required in the theorem.
\end{proof}

The following theorem allows normal (or integral) currents to be approximated by polyhedral chains which are simplicial chains not necessarily contained in an a priori complex.
Note in particular that the mass bounds can be made arbitrarily tight by choice of $\epsilon$ in contrast with the larger bounds of the deformation theorems.

\begin{theorem}[Polyhedral approximation of currents\cite{Federer1969}, 4.2.21, 4.2.24]
\label{thm:polyapprox}
If $\rho > 0$ and $T$ is a normal $m$-current in $\R^n$ supported in the interior of a compact subset $K$ of $\R^n$, then there exists a polyhedral chain $P$ with
\begin{subequations}
\begin{align}
\F_K(P-T) \leq \rho,\\
\label{eq:Pbound}\mass(P) < \mass(T) + \rho,\\
\label{eq:BPbound}\mass(\boundary P) < \mass(\boundary T) + \rho.
\end{align}
\end{subequations}
If $T$ is integral, then $P$ can be taken to be integral as well.
\end{theorem}
\begin{proof}
This is a slight modification of Federer's theorems which do not state \cref{eq:Pbound,eq:BPbound} separately but rather a combined bound $\mass(P) + \mass(\boundary P) \leq \mass(T) + \mass(\boundary T) + \rho$.
We show only the derivation of the separated bounds.

In the normal current case~\cite[4.2.24]{Federer1969}, these bounds follow from Federer's proof.
In particular, we have currents $P_1$, $P_2$ and $Y$ such that $P = P_1 + Y$ and the following bounds hold:
\begin{subequations}
\label{eq:Fbounds}
\begin{align}
\label{eq:FB1}\mass(P_1) < \mass(T) + \rho/4,\\
\label{eq:FB2}\mass(P_2) < \mass(\boundary T) + \rho/4,\\
\label{eq:FB3}\mass(P_2 - \boundary P_1 - \boundary Y) + \mass(Y) < \rho/2.
\end{align}
\end{subequations}
The bounds in \cref{eq:Pbound,eq:BPbound} follow from the triangle inequality and \cref{eq:FB1,eq:FB2,eq:FB3}:
\begin{align*}
\mass(P) &\leq \mass(P_1) + \mass(Y)\\
&< \mass(T) + \rho/4 + \rho/2,\\
\mass(\boundary P) &= \mass(\boundary P_1 + \boundary Y)\\
&\leq \mass(P_2 - \boundary P_1 - \boundary Y) + \mass(P_2)\\
&< \rho/2 + \mass(\boundary T) + \rho/4.
\end{align*}

In the integral current case~\cite[4.2.21]{Federer1969}, Federer applies the approximation theorem 4.2.20 to obtain $P$ close to the pushforward of $T$ under a Lipschitz diffeomorphism $f$.
That is, for any fixed $\epsilon > 0$, there exist $P$ and $f$ such that
\begin{subequations}
\begin{align}
\label{eq:Fpushbound}\mass(P - f_\# T) + \mass(\boundary P - \boundary f_\# T) \leq \epsilon\\
\label{eq:Fpushlip1}\lip(f) \leq 1 + \epsilon\\
\label{eq:Fpushlip2}\lip(f^{-1}) \leq 1 + \epsilon
\end{align}
From \cref{eq:Fpushbound,eq:Fpushlip1,eq:Fpushlip2}, we obtain mass bounds on $P$ and $\boundary P$:
\end{subequations}
\begin{subequations}
\begin{align}
\mass(P) &\leq \mass(f_\# T) + \epsilon\\
&\leq (1+\epsilon)^m \mass(T) + \epsilon\\
\mass(\boundary P) &\leq \mass(\boundary f_\# T) + \epsilon\\
&\leq (1+\epsilon)^{m-1} \mass(\boundary T) + \epsilon
%\mass(P) &\leq \mass(f_\# T) + \epsilon
\end{align}
\end{subequations}
%\begin{align}
%\mass(P) &\leq \mass(f_\# T) + \epsilon\\
%& (1+\epsilon)^m \mass(T) + \epsilon
%\end{align}
The bounds in \cref{eq:Pbound,eq:BPbound} follow by choosing $\epsilon$ small enough.
\end{proof}
\section{Results} \label{sec:results}

The simplicial deformation theorem can be modified to allow multiple currents to be deformed simultaneously by projecting from the same centers.
As opposed to using \cref{thm:sdt} separately on each current (where the centers of projection need not be the same), this yields a linearity result: deformations of linear combinations are linear combinations of deformations.
Pushing multiple currents at the same time comes at the cost of looser bounds on the deformation (linear in the number of currents) although slightly tighter analysis allows the bounds to be reduced by approximately a factor of 2 (\cref{cor:multisdtsdt}).

\begin{theorem}
\label{thm:multisdt}
Suppose $\epsilon > 0$ and we have the hypotheses of \cref{thm:sdt} except that there are now $m$ $d$-currents $T_1, T_2, \dots, T_m$ and $n$ $(d+1)$-currents $S_1, S_2, \dots, S_n$ to push on to the complex to yield the corresponding simplicial chains $P_i$ and $O_j$.
There is a series of projection centers (as in the proof of \cref{thm:sdt} and depending on $\epsilon$, $K$, the $T_i$ and $S_j$) which can be used with every current $T_i$ and $S_j$ to obtain the bounds:
\begin{align*}
\mass(P_i) &\leq ((2m+2n +\epsilon)\upvartheta_K)^{p-d} \mass(T_i) + \Delta ((2m+2n+\epsilon)\upvartheta_K)^{p-d+1} \mass(\boundary T_i)\\
\mass(\boundary P_i) &\leq ((2m+2n+\epsilon)\upvartheta_K)^{p-d+1}\mass(\boundary T_i)\\
\F(T_i,P_i) &\leq \Delta ((2m+2n+\epsilon)\upvartheta_K)^{p-d}(\mass(T_i)+(1+(2m+2n+\epsilon)\upvartheta_K)\mass(\boundary T_i))\\
\mass(O_j) &\leq ((2m+2n+\epsilon)\upvartheta_K)^{p-d-1} \mass(S_j) + \Delta ((2m+2n+\epsilon)\upvartheta_K)^{p-d} \mass(\boundary S_j)\\
\mass(\boundary O_j) &\leq ((2m+2n+\epsilon)\upvartheta_K)^{p-d}\mass(\boundary S_j)\\
\F(S_j,O_j) &\leq \Delta ((2m+2n+\epsilon)\upvartheta_K)^{p-d-1}(\mass(S_j)+(1+(2m+2n+\epsilon)\upvartheta_K)\mass(\boundary S_j))
\end{align*}
Moreover, if we let $\pi_K$ denote the projection map that uses these centers to push $(d-1)-$, $d-$, and $(d+1)$-currents to chains on the complex, then we have that:
\begin{itemize}
\item $\pi_K$ commutes with the boundary operator (i.e., $\pi_K(\boundary A) = \boundary \pi_K(A)$ where $A$ is any $d$- or $(d+1)$-current)
\item $\pi_K$ is linear on the currents $T_i$, $\boundary T_i$, $S_j$ and $\boundary S_j$.
That is, for any scalars $a_i$ and $b_j$,
\begin{align*}
\pi_K\left(\sum_{i=1}^m a_i \boundary T_i\right) &= \sum_{i=1}^m a_i \pi_K(\boundary T_i)\\
%&= \sum_{i=1}^m a_i \boundary \pi_K(T_i)\\
\pi_K\left(\sum_{i=1}^m a_i T_i + \sum_{j=1}^n b_j \boundary S_j\right) &= \sum_{i=1}^m a_i \pi_K(T_i) + \sum_{j=1}^n b_j \boundary (\pi_K(S_j))\\
\pi_K\left(\sum_{j=1}^n b_j S_j \right) &= \sum_{j=1}^n b_j \pi_K(S_j)
\end{align*}
\end{itemize}
\end{theorem}
\begin{proof}
We must show that there are centers in the set of feasible centers ${\cal B}_\sigma$ (see the proof sketch of \cref{thm:sdt}) which simultaneously achieve the various bounds on the $2(m+n)$ relevant currents: $T_1, \dots, T_m$, $\boundary T_1, \dots \boundary T_m$, $S_1, \dots S_n$, $\boundary S_1, \dots, \boundary S_n$.

We consider the case of projecting currents from the $\ell$-skeleton to the $(\ell-1)$-skeleton in the $\ell$-simplex $\sigma$.
As in the proof of \cref{thm:sdt}, we again use the average bound in \cref{eq:avgbound}.
For each $k \in \Z^+$ and $i = 1, 2, \dots, m$, let 
\[
H_{T_i,k} = \left\{\va \in {\cal B}_\sigma \,\middle|\, \int_\sigma J_d \phi(\vx,\va) \, {\rm d}\norm{T_i}(\vx) > \left(2m+2n+\frac{1}{k}\right)\upvartheta_\sigma\mass(T_i)\right\}.
\] 
Then, using the same average-based argument as in \cref{thm:sdt}, we have that $\mathcal{H}^\ell(H_{T_i,k})/\mathcal{H}^\ell({\cal B}_\sigma) < \frac{1}{2m+2n}$ (i.e., the size of the set of poorly behaved centers with respect to each $T_i$ is a small fraction of the set ${\cal B}_\sigma$ of possible centers).
We similarly define $H_{\boundary T_i,k}$, $H_{S_j,k}$, and $H_{\boundary S_j, k}$ and obtain the same bound of $\frac{1}{2m+2m}$ on the bad centers.
For each $k \in \Z^+$, we are interested in the set of centers which are simultaneously good centers for all currents involved (i.e., points in ${\cal B}_\sigma$ but not any of the $H_{\cdot, k}$ sets).
Call this set $G_k$ and observe that it has positive measure:
\begin{align*}
{\cal H}^\ell(G_k) &={\cal H}^\ell\left({\cal B}_\sigma \backslash\left(\bigcup_{i=1}^m H_{T_i, k} \cup \bigcup_{i=1}^m H_{\boundary T_i, k}\cup \bigcup_{i=1}^n H_{S_i,k} \cup \bigcup_{i=1}^n H_{\boundary S_i,k}\right)\right)\\
&\geq {\cal H}^\ell({\cal B}_\sigma) - \sum_{i=1}^m{\cal H}^\ell(H_{T_i, k}) - \sum_{i=1}^m{\cal H}^\ell(H_{\boundary T_i, k}) - \sum_{j=1}^n{\cal H}^\ell(H_{S_j,k}) - \sum_{j=1}^n{\cal H}^\ell(H_{\boundary S_j,k})\\
&> {\cal H}^\ell({\cal B}_\sigma)\left(1 - \frac{m}{2m+2n} - \frac{m}{2m+2n} - \frac{n}{2m+2n} - \frac{n}{2m+2n}\right)\\
&=0.
\end{align*}
Thus for any $k > \frac{1}{\epsilon}$ we have that $G_k$ is a nonempty set of possible projection centers which simultaneously attain an expansion bound of at most $(2m+2n + \epsilon)\upvartheta_\sigma$ for all the pertinent currents.

The projection operator is clearly linear and commutes with the boundary operator as a consequence of properties\cite[4.1.6]{Federer1969} of the differential forms to which currents are dual.
\end{proof}

\begin{corollary}
\label{cor:multisdtsdt}
The bounds in \cref{thm:sdt} can all be tightened by replacing $4\upvartheta_K$ with $(2+\epsilon) \upvartheta_K$.
\end{corollary}
\begin{proof}
Simply taken $m = 1$ and $n = 0$ in \cref{thm:multisdt}.
\end{proof}

For a 2-complex $K$, the minimum angle over all triangles in the complex is easier to work with and can be used as a proxy for our simplicial regularity constant as \cref{lem:2danglereg} indicates.

\begin{lemma}
\label{lem:2danglereg}
A lower bound on the minimum angle of all triangles in a 2-complex implies an upper bound on the simplicial regularity constant.
That is, given a 2-complex $K$ with minimum angle at least $\theta$, we have $\upvartheta_K \leq C_\theta$ for some constant $C_\theta$.
\end{lemma}
\begin{proof}
The simplicial regularity constant $\upvartheta_K$ used for \cref{thm:sdt,thm:multisdt} in the case of triangles is given by
\[
\upvartheta_K = \frac{4}{\pi}\sup_{\sigma \in K}\frac{\diam(\sigma)\perimeter(\sigma)}{\inrad(\sigma)^2} + 2 \sup_{\sigma \in K}\frac{\diam(\sigma)}{\inrad(\sigma)}.
\]
We observe that bounding $\diam(\sigma)/\inrad(\sigma)$ and $\perimeter(\sigma)/\inrad(\sigma)$ for all triangles $\sigma \in K$ yields a bound for $\upvartheta_K$.
Suppose $\sigma$ has side lengths $a \geq b \geq c$ and angle $\gamma$ opposite $c$.
Using the law of cotangents, we obtain
\[
\frac{\diam(\sigma)}{\inrad(\sigma)} = \frac{a\cot(\gamma/2)}{(a+b)/2-c/2} \leq \frac{a\cot(\gamma/2)}{(a+b)/2-b/2} = 2\cot(\gamma/2) \leq 2\cot(\theta/2).
\]
The bound for $\perimeter(\sigma)/\inrad(\sigma)$ follows easily from this:
\[
\frac{\perimeter(\sigma)}{\inrad(\sigma)} \leq \frac{3\diam(\sigma)}{\inrad(\sigma)} < 6\cot(\theta/2).
\]
Thus we can take $C_\theta = \frac{48}{\pi}\cot(\theta/2)^2 + 4\cot(\theta/2)$.
\end{proof}

Our result relies on the ability to localize irregularities via subdivision, focusing on localization rather than removal because the latter is not possible.
For example, any subdivision of a 2-complex with a very small input angle will have an angle that is at least as small.
With that in mind, we require that subdivisions be possible which push the irregularities into the corners.
That is, the irregularity should be bounded by a constant (independent of the complex) away from the original complex skeleton and a complex-dependent constant (reflecting the necessity of some bad simplices) near the skeleton.
\cref{con:boundreg} formalizes this requirement and \cref{thm:2dboundreg} notes some cases where it holds.
We present our main theorem in such a way that proving \cref{con:boundreg} more generally will automatically extend our results.

\begin{conjecture}
\label{con:boundreg}
For any $p$-dimensional simplicial complex $K$ in $\R^q$ and $\epsilon > 0$, it is possible to subdivide $K$ so that all simplices are of bounded ``badness'' (with bound independent of $K$ or $\epsilon$) except possibly for simplices in a region of $p$-dimensional volume less than $\epsilon$ near the $(p-1)$-skeleton; even these simplices have bounded badness (dependent on $K$ but not $\epsilon$).
More precisely, there exists a subdivision $M_\epsilon$ of $K$ and a subcomplex $M_\epsilon'$ of $M_\epsilon$ (with simplicial regularity constants $\upvartheta_{M_\epsilon}$ and $\upvartheta_{M_\epsilon'}$) such that:
\begin{enumerate}
\item $M_\epsilon \backslash M_\epsilon' \subseteq \{x \in \R^q \mid \|x-y\| < \epsilon \text{ for some $y$ in the $(p-1)$-skeleton of $K$}\}$,
\item $\upvartheta_{M_\epsilon} \leq \alpha_K$ for some constant $\alpha_K$,
\item $\upvartheta_{M_\epsilon'} \leq \beta$ for some fixed constant $\beta$.
\end{enumerate}
In particular, $\alpha_K$ does not depend on $\epsilon$ and $\beta$ does not depend on $K$ or $\epsilon$.
The simplicial regularity constants are defined as in \cref{eq:regularity}.
\end{conjecture}

\begin{theorem}
\label{thm:2dboundreg}
\cref{con:boundreg} holds for:
\begin{itemize}
\item $q \geq p = 1$
\item $p = q = 2$
\end{itemize}
\end{theorem}
\begin{proof}
The $p = 1$ case is trivial as all 1-simplices have the same regularity so we have $\upvartheta_K = 8$ and can take $M_\epsilon = M_\epsilon' = K$.

For the $p = q = 2$ case, we proceed in two steps.
First we will superimpose a square grid on $K$ (orienting it to bound the minimum angle created between its edges and those of $K$), creating a cell complex which is a refinement of $K$.
Next we use Shewchuk's Terminator algorithm\cite{Sh2002} to further refine the cell complex back into a simplicial complex with bounds on the minimum angle and, crucially, restrictions on where these small angles can be so that we can obtain regularity bounds.

By superimposing a fine enough square grid, we can force the small angles (whether already present in the complex or newly created) to occur only in a small measure subset of the complex.
Pick $\delta > 0$ small enough that the set
\[
\{x \in \R^2 \mid y\text{ lies on the 1-skeleton of $K$, } \|x-y\| < 3\delta\}
\]
 has measure less than $\epsilon$.
Let $G$ be a finite square grid in $\R^2$ whose cells each have diameter $\delta$ such that $G$ covers the underlying space of $K$ in any rotation.
Note that there are only two directions present in $G$ so if we bound all possible angles created between these directions and the edges of $K$, we can bound the minimum new angle created by superimposing $G$.

Let $w \in \R^2$ be a fixed unit vector and define
\[
E = \left\{\phi, \phi + \frac{\pi}{2} \mid \phi \text{ is the angle between $u-v$ and $w$ for some edge $(u,v) \in K$}\right\}.
\]
Further let $E^\theta = \{\psi \in [0, 2\pi) \mid |\phi-\psi| < \theta\text{ for some }\phi \in E\}$.
This is the set of angles to avoid when rotating $G$ in order to guarantee all created angles will be $\theta$ or larger.

Denote by $N < \infty$ the cardinality of $E$ and note that $[0, 2\pi) \backslash E^\frac{\pi}{2N}$ has positive measure so there exist rotations of the square grid that create no new angles smaller than $\frac{\pi}{2N}$.

After superimposing a suitably rotated version of $G$, we obtain a new cellular complex which is a refinement of $K$.
This is a planar straight line graph which can be used as input to Shewchuk's Terminator algorithm\cite{Sh2002} which refines it into a simplicial complex $M_\epsilon$ with guarantees about the minimum angle bound of the resulting complex and where the small angles can occur.

In particular, if $\theta$ be the minimum angle in the cellular complex (either present originally or added in the square grid superposition), then the minimum angle of $M_\epsilon$ is at least $\arcsin((\sqrt{3}/2)\sin(\theta/2))$.
Furthermore, no angles less than $30^\circ$ are created by the algorithm except in the vicinity of angles less than $60^\circ$.
Specifically, newly created small angles must be part of a skinny triangle whose circumcenter encroaches upon a subsegment cluster bearing a small input angle.
As all such subsegment clusters must be contained within a distance of $2\delta$ of the 1-skeleton of $K$, we have that all small angles in $M_\epsilon$ are within $3\delta$ of the 1-skeleton of $K$.

Let $M_\epsilon'$ be the subcomplex of $M_\epsilon$ containing all triangles not fully contained in the $3\delta$ tube, noting that all angles in $M_\epsilon'$ are at least $30^\circ$ so by \cref{lem:2danglereg} we have
\[
\upvartheta_{M_\epsilon'} \leq \frac{48}{\pi}\cot(15^\circ)^2 + 4\cot(15^\circ) = \frac{4(2+\sqrt{3})(24+12\sqrt{3}+\pi)}{\pi}. %\approx 227.74
\]
We may take $\beta$ to be this quantity, noting that it is independent of $\epsilon$ and $K$.
The minimum angle bound $\theta$ for $M_\epsilon$ and \cref{lem:2danglereg} give us a bound $\alpha_K$ for $\upvartheta_{M_\epsilon}$ (independent of $\epsilon$).
\end{proof}

The following theorem shows that the bounds in \cref{thm:multisdt} may be replaced with constants independent of the complex and currents involved if we subdivide the complex by means of \cref{con:boundreg} (the subdivision does depend on the currents and complex, of course).

\begin{theorem}
\label{thm:2dboundedsdt}
Suppose we have integers $d < s \leq q$ and that \cref{con:boundreg} holds for the given $q$ and any $p$ such that $d-1 \leq p \leq s$ (that is, suppose we can isolate the irregularities of any $p$-complex in $\R^q$ by suitable subdivision).
Given a $s$-dimensional simplicial complex $K$ in $\R^q$ and a set of $d$-currents $T_1, \dots, T_m$ and $(d+1)$-currents $S_1, \dots, S_n$ in the underlying space of $K$ with $d < s$, there exists a complex $K'$ which is a subdivision of $K$ such that we have all of the conclusions of \cref{thm:multisdt} (i.e., mass and flat norm bounds and linear projection of the $T_i$ and $S_j$ to $K'$) except the simplicial irregularity constant $\upvartheta_{K'}$ in the various bounds can be replaced with a constant $L$ which does not depend on $K$.
\end{theorem}
\begin{proof}
In the simplicial deformation theorems, the current is projected step-by-step to lower dimensional skeletons (e.g., a $d$-current is projected from the initial $p$-complex to the $(p-1)$-skeleton, then the $(p-2)$-skeleton eventually down to the $d$-skeleton with one more step to push the current's boundary to the $(d-1)$-skeleton) with each projection being done by picking a center in each simplex and using it to project outward to the boundary of the simplex.
The simplicial regularity constant is used to bound the expansion of mass at each projection step and is defined by \cref{eq:regularity}, a bound on the regularity of all simplices in the complex.

However, this is a bit stronger than required as the projection is a local operation and the bound at each step depends only on the simplicial regularity of the simplex in question.
In addition, there is no reason in principle that we cannot subdivide the complex in between steps.
That is, after pushing to the $\ell$-skeleton, we can further subdivide the complex and then push to the newly refined $(\ell-1)$-skeleton.
In this case, the subdivision need not preserve the simplicial regularity of the $(\ell+1)$- or higher simplices as all subsequent pushing steps will take place in lower dimensional simplices.
Moreover, for a given portion of current we can use the maximum of the simplicial regularity constants of the simplices it encounters while being pushed (rather than the maximum over all simplices in the complex).

For all $\epsilon > 0$ and nonnegative integers $k < p$, let $N_{k}^\epsilon$ denote the set of all points in the $(k+1)$-skeleton of $K$ with positive distance less than $\epsilon$ from the $k$-skeleton of $K$ (i.e., all points in the interior of the $(k+1)$-simplices of $K$ which are close to the $k$-skeleton).
Let $T \currentrestr N_{p-1}^\epsilon$ denote the restriction of the current $T$ to the set $N_{p-1}^\epsilon$ and note that
\begin{align}
\label{eq:limrestr}
\begin{split}
\lim_{\epsilon \downarrow 0} \mass(T_i \currentrestr N_{p-1}^\epsilon) = 0, & \quad
\lim_{\epsilon \downarrow 0} \mass(S_j \currentrestr N_{p-1}^\epsilon) = 0, \\
\lim_{\epsilon \downarrow 0} \mass(\boundary T_i \currentrestr N_{p-1}^\epsilon) = 0, &\quad
\lim_{\epsilon \downarrow 0} \mass(\boundary S_j \currentrestr N_{p-1}^\epsilon) = 0.
\end{split}
\end{align}
Let 
\begin{equation}
\label{eq:currentdelta}
\delta = \frac{\beta}{\alpha_K}\min_{1\leq i \leq m,1 \leq j \leq n}\{\mass(T_i), \mass(\boundary T_i), \mass(S_j), \mass(\boundary S_j)\}
%\delta = \frac{2^{1/(p-d+1)}-1}{2(m+n)\alpha_K}\min_{1\leq i \leq m,1 \leq j \leq n}\{\mass(T_i), \mass(\boundary T_i), \mass(S_j), \mass(\boundary S_j)\}
\end{equation}
 where $\alpha_K$ and $\beta$ are as in (the assumed true) \cref{con:boundreg} and choose $\epsilon > 0$ to make each of the masses in \cref{eq:limrestr} less than $\delta$.
We can apply \cref{con:boundreg} with this $\epsilon$ to obtain a subdivision $M_\epsilon$ of $K$ and a subcomplex $M_\epsilon'$ such that the portion of each of our currents which lies in $M_\epsilon \backslash M_\epsilon'$ and is not already on the $(p-1)$-skeleton (so is not fixed by the first projection) has mass less than $\delta$.
This portion of each current increases in mass by a factor of at most $(2m+2n+\epsilon)\alpha_K$ when projecting to the $(p-1)$-skeleton (see proof of \cref{thm:multisdt}).
Letting $T_i'$ denote the result of projecting $T_i$ to the $(p-1)$-skeleton, we can bound its mass using \cref{eq:currentdelta}:
\begin{align*}
\mass(T_i') &\leq (2m+2n+\epsilon)\Big[\beta\mass(T_i\currentrestr M_\epsilon'\backslash \skel{p-1}(K)) + \alpha_K\mass(T_i \currentrestr M_\epsilon\backslash (M_\epsilon' \cup \skel{p-1}(K)))\Big] \\
&\quad+ \mass(T_i\currentrestr \skel{p-1}(K))\\
&\leq (2m+2n+\epsilon)(\beta\mass(T_i) + \alpha_K\delta)\\
% + \mass(T_i\currentrestr \skel{p-1}(K))\\
&\leq (2m+2n+\epsilon)(\beta\mass(T_i) + \beta\mass(T_i )) \\
&\leq (2m+2n+\epsilon)(2\beta)\mass(T_i).
\end{align*}
Similar inequalities hold for $S_j$, $\boundary T_i$, and $\boundary S_j$.
In the preceding, we have accomplished the goal of projecting all currents involved from the $p$-skeleton to the $(p-1)$-skeleton and can now consider them as currents in the underlying space of the $(p-1)$-complex $\skel{p-1}(K)$.
We can apply this procedure iteratively (use \cref{con:boundreg} to localize the irregularities and then project) to push to the $(p-2)$, etc. skeletons.

When we subdivide each $k$-skeleton using \cref{con:boundreg}, the higher dimension simplices are not subdivided by default but this is easy to fix.
After a $k$-simplex is subdivided, add a point to the interior of every $(k+1)$-simplex of which it was a face and connect the new point to every $k$-simplex on its boundary.
This will likely generate highly irregular simplices but since we've already pushed the currents down beyond their dimension, it isn't an issue.

This argument continues in the same way as \cref{thm:sdt,thm:multisdt} and establishes our result with $L = 2\beta$.
\end{proof}

\begin{theorem}
\label{thm:icmainresult}
If $T$ is an integral $d$-current in $\R^{d+1}$ and \cref{con:boundreg} holds for $d-$ and $(d+1)$-currents, then some flat norm minimizer for $T$ is an integral current.
That is, there is an integral $d$-current $X_I$ and integral $(d+1)$-current $S_I$ such that $\F(T) = \mass(X_I) + \mass(S_I)$ and $T = X_I + \boundary S_I$.
\end{theorem}
\begin{proof}
%TODO: Finite length
We let $X + \boundary S$ be an optimal flat norm decomposition of $T$.
That is, $X$ is a $d$-current and $S$ is a $(d+1)$-current such that
\begin{align}
\label{eq:Tdecomp}
T &= X + \boundary S, & \F(T) &= \mass(X) + \mass(S).
\end{align}
We note by \cref{lem:normaldecomp} that $X$ and $S$ are normal currents.

As a general outline of the proof, for each $\delta > 0$, we will choose a particular simplicial complex $K_\delta$ on which we have $d$-chains $P_\delta$ and $X_\delta$ and $(d+1)$-chain $S_\delta$ respectively approximating $T$, $X$, and $S$ with error at most $\delta$.
We convert the (possibly nonintegral) optimal flat norm decomposition of $T$ into a candidate simplicial decomposition of $P_\delta$ in order to show (\cref{claim:FTeq}) the simplicial flat norm of $P_\delta$ converges to the flat norm of $T$ (this step does not yet show that the flat norm decompositions converge).
We can take the optimal simplicial decomposition to be integral for each $P_\delta$ by \cref{thm:simpint}.
The compactness theorem from geometric measure theory along with the above convergence result allows us to take the limit of (a subsequence of) these integral simplicial decompositions and obtain an integral flat norm decomposition of $T$ (\cref{claim:Tintdecomp}).

%We require these be chosen such that the support of $T$ is in the underlying space of $K_\delta$, $K_\delta$ is totally unimodular and every simplex of $K_\delta$ has diameter less than $\delta$.

Suppose $\delta > 0$ and apply \cref{thm:polyapprox} to obtain polyhedral currents $P_\delta$, $X_\delta$ and $S_\delta$ with 
\begin{subequations}
\label{eq:polyhedralapproximations}
\begin{align}
\label{eq:Tapprox}
\F(T-P_\delta) &< \delta, & \mass(P_\delta) &< \mass(T) + \delta, & \mass(\boundary P_\delta) & < \mass(\boundary T) + \delta,\\
\label{eq:Xapprox}
% & \mass(\boundary P_\delta) &<\mass(\boundary T) + \delta/3\\
\F(X-X_\delta) &< \delta, & \mass(X_\delta) &< \mass(X) + \delta, & \mass(\boundary X_\delta) & < \mass(\boundary X) + \delta,\\
\label{eq:Sapprox}
\F(S-S_\delta) &< \delta, & \mass(S_\delta) &< \mass(S) + \delta, & \mass(\boundary S_\delta) & < \mass(\boundary S) + \delta.
%\mass(X_\delta) &< \mass(X) + \delta/3,\\
%\mass(S_\delta) &< \mass(S) + \delta/3.
\end{align}
\end{subequations}
We also require optimal flat norm decompositions of $P_\delta - T$, $X-X_\delta$ and $S-S_\delta$ so let $U_i^\delta$, $W_j^\delta$ and $V_2^\delta$ be $d$-, $(d+1)$- and $(d+2)$-dimensional currents such that:
\begin{subequations}
\label{eq:decompositions}
\begin{align}
\label{eq:PTdecomp}
P_\delta - T &= U_0^\delta + \boundary W_0^\delta, & \F(P_\delta - T) &= \mass(U_0^\delta) + \mass(W_0^\delta),\\
\label{eq:Xdecomp}
X - X_\delta &= U_1^\delta + \boundary W_1^\delta, & \F(X - X_\delta) &= \mass(U_1^\delta) + \mass(W_1^\delta),\\
\label{eq:Sdecomp}
S - S_\delta &= W_2^\delta + \boundary V_2^\delta, & \F(S - S_\delta) &= \mass(W_2^\delta) + \mass(V_2^\delta).
\end{align}
\end{subequations}
To clarify the notation, we adopt the convention that variables with a $\delta$ subscript are chains on the simplicial complex $K_\delta$ whereas a $\delta$ superscript merely indicates dependence on $\delta$.

Let $K_\delta$ be any simplicial complex that triangulates $P_\delta$, $X_\delta$ and $S_\delta$ separately as well as the convex hull of their union.
% such every simplex in $K_\delta$ has diameter bounded by $\delta$.
%TODO: cite a result saying this can be done.
We may assume (applying the subdivision algorithm of Edelsbrunner and Grayson\cite{EdGr2000} and \cref{thm:2dboundedsdt} if necessary) that the currents $U_0$, $U_1$, $W_0$, $W_1$, and $W_2$ can be pushed to $K_\delta$ with expansion bound at most $L$ and the maximum diameter $\Delta$ of a simplex of $K_\delta$ satisfies
\begin{align}
\label{eq:diameter}
\Delta &\leq \frac{\delta}{\max\{1, \mass(\boundary U_0^\delta), \mass(\boundary U_1^\delta), \mass(\boundary W_0^\delta), \mass(\boundary W_1^\delta), \mass(\boundary W_2^\delta)\}}.
\end{align}

\begin{subclaim}
\label{claim:FTleq}
$\F(T) \leq \lim_{\delta \downarrow 0} \F_{K_\delta}(P_\delta)$
\end{subclaim}
\begin{subproof}
%Let $X+\boundary S = T$ be an optimal flat norm decomposition for $T$ and $Y_\delta+\boundary R_\delta = P_\delta$ be an optimal simplicial flat norm decomposition for $P_\delta$ (that is, $\F(T) = \mass(X)+\mass(S)$ and $\F_{K_\delta}(P_\delta) = \mass(Y_\delta)+\mass(R_\delta)$)
%This is possible since TODO.
%
%Let $T-P_\delta = A_\delta + \boundary B_\delta$ be an optimal classical flat norm decomposition of the current $T-P_\delta$.
%
%The main idea here is that an optimal simplicial flat norm decomposition for $P_\delta$ can be turned into a candidate flat norm decomposition for $T$.
By the triangle inequality and since any simplicial flat norm decomposition is a candidate decomposition for the flat norm, we have
\begin{align*}
\F(T) &\leq \F(T-P_\delta) + \F(P_\delta)\\
 &\leq \F(T-P_\delta)+\F_{K_\delta}(P_\delta).
% &= \F(T-P_\delta)+\mass(Y_\delta) + \mass(R_\delta).
\end{align*}
The claim follows from letting $\delta \downarrow 0$ and noting that $\F(T-P_\delta) \to 0$.
\end{subproof}

\begin{subclaim} $\F(T) = \lim_{\delta \downarrow 0} \F_{K_\delta}(P_\delta)$
\label{claim:FTeq}
\end{subclaim}
\begin{subproof}
In light of \cref{claim:FTleq}, we must show that $\F(T) \geq \lim_{\delta \downarrow 0} \F_{K_\delta}(P_\delta)$.
%We will show that for $\delta$ sufficiently small, the optimal flat norm decomposition for $T$ induces a simplicial flat norm decomposition for $P_\delta$ which is strictly better than our assumption allows.

%There exists $\delta^* > 0$ such that $\F_{K_\delta}(P_\delta) - \F(T) > \Delta$ for all $0 < \delta \leq \delta^*$.

Recall that $X+\boundary S=T$ is an optimal flat norm decomposition of $T$ and $X_\delta$ and $S_\delta$ are polyhedral approximations to $X$ and $S$ on our simplicial complex $K_\delta$.
Using the decompositions in \cref{eq:Tdecomp,eq:decompositions}, we can write:
\begin{align}
\label{eq:Pdecomp}
\begin{split}
P_\delta &= T+U_0^\delta + \boundary W_0^\delta\\
&= X+\boundary S+U_0^\delta + \boundary W_0^\delta\\
&= X_\delta + U_0^\delta + U_1^\delta + \boundary(S_\delta + W_0^\delta + W_1^\delta + W_2^\delta).
\end{split}
\end{align}
Now apply \cref{thm:2dboundedsdt} with $\epsilon = 1$ to the currents $U_i^\delta$ and $W_j^\delta$ for all $i \in \{0, 1\}$ and $j \in \{0, 1, 2\}$ to obtain $U_{i,\delta}$ and $W_{j,\delta}$ on the simplicial complex $K_\delta$ with
\begin{subequations}
\label{eq:pushdecomp}
\begin{align}
\mass(U_{i,\delta}) &\leq (11L)^{p-d+1}\mass(U_i^\delta) + (11L)^{p-d}\Delta \mass(\boundary U_i^\delta), \\
\mass(W_{j,\delta}) &\leq (11L)^{p-d}\mass(W_j^\delta) + (11L)^{p-d-1}\Delta \mass(\boundary W_j^\delta)).
\end{align}
\end{subequations}
Applying \cref{eq:decompositions,eq:polyhedralapproximations,eq:diameter}, we obtain the following from \cref{eq:pushdecomp}:
\begin{subequations}
\label{eq:pushbounds}
\begin{align}
\begin{split}
\mass(U_{i,\delta}) & \leq (11L)^{p-d+1}\delta + (11L)^{p-d}\frac{\delta}{\mass(\boundary U_i^\delta)}\mass(\boundary U_i^\delta)\\
& = (11L)^{p-d}(1 + 11L)\delta
\end{split}\\
\begin{split}
\mass(W_{j,\delta}) & \leq (11L)^{p-d}\delta + (11L)^{p-d-1}\frac{\delta}{\mass(\boundary W_j^\delta)}\mass(\boundary W_j^\delta)\\
&= (11L)^{p-d-1}(1+11L)\delta
\end{split}
\end{align}
\end{subequations}
We apply the linearity result of \cref{thm:2dboundedsdt} to \cref{eq:Pdecomp} along with the fact that $P_\delta$, $X_\delta$, and $\boundary S_\delta$ are fixed by projection to the $d$-skeleton of $K_\delta$ to yield
\begin{align*}
P_\delta &=  (X_\delta + U_{0,\delta} + U_{1,\delta}) + \boundary(S_\delta + W_{0,\delta} + W_{1,\delta} + W_{2,\delta})
\end{align*}
which, as all quantities are chains on $K_\delta$, is a candidate simplicial flat norm decomposition of $P_\delta$.
Using this observation, the triangle inequality, and \cref{eq:pushbounds,eq:polyhedralapproximations}, we have
\begin{align*}
\F_{K_\delta}(P_\delta) &\leq \mass(X_\delta + U_{0,\delta} + U_{1,\delta}) + \mass(S_\delta + W_{0,\delta} + W_{1,\delta} + W_{2,\delta})\\
&\leq \mass(X_\delta) + \mass(U_{0,\delta}) + \mass(U_{1,\delta}) + \mass(S_\delta) + \mass(W_{0,\delta}) + \mass(W_{1,\delta}) + \mass(W_{2,\delta})\\
&\leq \mass(X) + \mass(S) + 2\delta + 2(11L)^{p-d}(1 + 11L)\delta + 3(11L)^{p-d-1}(1+11L)\delta\\
%&< \mass(X) + \mass(U_0^\delta) + \mass(U_1^\delta) + \mass(U_2,\delta) + \mass(S) + \mass(W_0^\delta) + \mass(W_1^\delta) + \frac{7\Delta}{12}\\
%&< \mass(X) + \mass(S) + \Delta\\
&= \F(T) + 2\delta + 2(11L)^{p-d}(1 + 11L)\delta + 3(11L)^{p-d-1}(1+11L)\delta.
\end{align*}
The claim follows from taking the limit as $\delta \downarrow 0$.

\end{subproof}

\begin{subclaim}
\label{claim:simpdecomp}
For each $\delta > 0$, there exist integral simplicial chains $Y_\delta$ and $R_\delta$ on $K_\delta$ such that $P_\delta = Y_\delta + \boundary R_\delta$ is an optimal simplicial flat norm decomposition (i.e., $\F_{K_\delta}(P_\delta) = \mass(Y_\delta) + \mass(R_\delta)$).
\end{subclaim}
\begin{subproof}
This follows from \cref{thm:simpint}.
\end{subproof}

\begin{subclaim}
\label{claim:simpdecompbounded}
There exists $c > 0$ such that for all $\delta \leq 1$, the currents $Y_\delta$, $\boundary Y_\delta$, $R_\delta$, and $\boundary R_\delta$ all have mass at most $c$.
%$There exists a $c > 0$ such that for all $\delta \leq 1$, 
%\begin{align*}
%\mass(Y_\delta) &\leq c, & \mass(\boundary Y_\delta) &\leq c,\\
%\mass(R_\delta) &\leq c, & \mass(\boundary R_\delta) &\leq c.
%\end{align*} 
\end{subclaim}
\begin{subproof}
Using the fact that $P_\delta = Y_\delta + \boundary R_\delta$ is an optimal simplicial flat norm decomposition and facts from \cref{eq:polyhedralapproximations}, we observe that
\begin{align*}
\begin{aligned}[c]
\mass(Y_\delta) &\leq \mass(P_\delta) \\
&< \mass(T) + \delta \\
&\leq \mass(T) + 1,\\
\\
\mass(R_\delta) &\leq \mass(P_\delta)\\
&< \mass(T) + \delta \\
&\leq \mass(T) + 1,\\
\end{aligned}
\qquad\qquad
\begin{aligned}
\mass(\boundary Y_\delta) &= \mass(\boundary (P_\delta - \boundary R_\delta))\\
 &= \mass(\boundary P_\delta)\\
 &< \mass(\boundary T) + \delta\\
 &\leq \mass(\boundary T) + 1,\\
\mass(\boundary R_\delta) &= \mass(P_\delta - Y_\delta)\\
&\leq \mass(P_\delta) + \mass(Y_\delta)\\
&< 2\mass(T) + 2.
\end{aligned}
\end{align*}
%\begin{align*}
%\mass(Y_\delta) &\leq \mass(P_\delta)\\
% &< \mass(T) + \delta \\
%&\leq \mass(T) + 1,\\
%\mass(R_\delta) &\leq \mass(P_\delta)\\
%&< \mass(T) + 1,\\
%\mass(\boundary Y_\delta) &= \mass(\boundary (P_\delta - \boundary R_\delta))\\
%& = \mass(\boundary P_\delta)\\
%& < \mass(\boundary T) + \delta\\
%& \leq \mass(\boundary T) + 1.\\
%\mass(\boundary R_\delta) &= \mass(P_\delta - Y_\delta)\\
%&\leq \mass(P_\delta) + \mass(Y_\delta)\\
%&< 2\mass(T) + 2 
%\end{align*}
So $c = \max\{2\mass(T)+ 2, \mass(\boundary T) + 1\}$ works.

\end{subproof}

\begin{subclaim}
\label{claim:Tintdecomp}
There is an optimal flat norm decomposition of $T$ with integral currents.
\end{subclaim}
\begin{subproof}
The compactness theorem\cite{Federer1969,Morgan2008} states that given any closed ball $K$ in $\R^n$ and nonnegative constant $c$, the set 
\[
\{I \text{ is an integral $p$-current in $\R^n$} \mid \mass(I) \leq c, \mass(\boundary I) \leq c, \spt I \subseteq K\}
\]
 is compact with respect to the flat norm.
In light of \cref{claim:simpdecompbounded}, this means there is a compact set of integral currents containing $Y_\delta$ for all $\delta \leq 1$ (and similarly for $R_\delta$).

Let $\delta_n = \frac{1}{n}$ and consider the sequences $\{Y_{\delta_n}\}$ and $\{R_{\delta_n}\}$.
By compactness, there exists a subsequence $\{\delta^*_n\}$ of $\{\delta_n\}$ and integral currents $Y^*$ and $R^*$ such that $Y_{\delta^*_n} \to Y^*$ and $R_{\delta^*_n} \to R^*$ in the flat norm.
By \cref{lem:convprops}, we have $Y_{\delta^*_n} + \boundary R_{\delta^*_n} \to Y^* + \boundary R^*$.
Applying \cref{claim:simpdecomp} and \cref{claim:FTeq}, we get $\mass(Y_{\delta^*_n}) + \mass(R_{\delta^*_n}) = \F_{K_\delta} (P_{\delta^*_n}) \to \F(T)$.

Since $Y_\delta + \boundary R_\delta = P_\delta \to T$, we also have $Y_{\delta^*_n} + \boundary R_{\delta^*_n} \to T$.
That is, $T = Y^* + \boundary R^*$.
As mass is lower semicontinuous with respect to convergence in the flat norm and by \cref{claim:FTeq}, we have that 
\begin{align*}
\mass(Y^*) + \mass(R^*) &\leq \lim_{n\to\infty} \mass(Y_{\delta^*_n}) + \mass(R_{\delta^*_n}) \\
&= \lim_{n\to\infty} \F_{K_{\delta^*_n}}(P_\delta)\\
&= \F(T).
\end{align*}
Thus $\mass(Y^*) + \mass(R^*) = \F(T)$ and $Y^* + \boundary R^*$ is an optimal flat norm decomposition of $T$.
\end{subproof}
\end{proof}

\chapter{Nonasymptotic densities\footnote{Previously published as \cite{ibrahim-2014-1}}}
\label{ch:nad}

\newcommand{\tc}[2]{T_{#1}(#2)}

\newcommand{\disk}[2]{D(#1, #2)}
\renewcommand{\cir}[2]{C(#1, #2)}
\newcommand{\shape}{\Omega}
\newcommand{\bd}{\partial \shape}
\newcommand{\tgl}{tangentially graph-like}
\newcommand{\tcgl}{tangent-cone graph-like}
\newcommand{\mono}{horizontally strictly monotonic}
\newcommand{\gm}[1]{\gamma (#1)}
\newcommand{\dgm}[1]{{\gamma'}(#1)}
\newcommand{\dgmm}{{\gamma'}}
\newcommand{\delep}{\delta_{\epsilon}}
\newcommand{\enu}{\epsilon_\nu}
\newcommand{\rhh}{\hat{\rho}}
\newcommand{\prff}{\noindent {\bf Proof: }\setcounter{clm}{0}}
\newcommand{\prfclm}[1]{\noindent {\bf Proof of Claim \ref{#1}: }}
\newcommand{\epf}{$\blacksquare \;\;$}
\newcommand{\epfclm}{$\Box \;\;$}

\title{Nonasymptotic densities for shape reconstruction}
\author{Sharif Ibrahim\footnote{Department of Mathematics, Washington State University, Pullman, WA 99164-3113, USA} \and Kevin Sonnanburg\footnote{Department of Mathematics, University of Tennessee Knoxville, Knoxville, TN 37996-1320}
\and Thomas J. Asaki\footnotemark[1]
\and Kevin R. Vixie\footnotemark[1]
}

\newenvironment{mylist}{\begin{list}{*}{\itemsep=.01in \parsep=.01in \topsep=.1in \parskip=.01in}}{\end{list}}
\newenvironment{mybiglist}{\begin{list}{*}{\itemsep=.1in \parsep=.1in \topsep=.2in \parskip=.1in}}{\end{list}}

\section{Introduction}
\label{sec:intro}

This work discusses the integral area invariant introduced by Manay et al.\cite{manay-2004-1}, particularly with regard to reconstructability of shapes.
This topic has been considered previously by Fidler et al.\cite{fidler-2008-1}\cite{fidler-2012-1} for the case of star-shaped regions.
Recent results have shown local injectivity in the neighborhood of a circle \cite{bauer-2012-1} and for graphs in a neighborhood of constant functions \cite{calder-2012-1}.

The present work does not assume a star-shaped condition but does make use of a \tcgl{} condition which is local to the integral area circle.
We also present an interpretation of the integral area invariant as a nonasymptotic density.
This is based on a poster presented by the authors\cite{ibrahim-poster-2010}.

Our \tgl{} and \tcgl{} conditions (definitions \ref{def:tgl} and \ref{def:tcgl} in section \ref{sec:notation}) restrict our attention to shapes with boundaries that can locally (i.e., within radius $r$) be viewed as graphs of functions in a Cartesian plane in one particular orientation (in the case of \tgl{}) or a particular set of orientations (for \tcgl{}).
Intuitively, these conditions guarantee that the boundary does not turn too sharply within the given radius and that working locally in Euclidean space is the same as working locally on the boundary of our shapes (i.e., the shape boundary does not pass through any given invariant circle multiple times, section \ref{sec:tgl_2arc}).
These simplifying assumptions allow us to explicitly analyze what happens when we move along the boundary and to work locally without worrying about global effects.

We show that the \tcgl{} property can be preserved when approximating a shape with a polygon (section \ref{sec:tgl-ncgl}) and discuss what the derivatives of these nonasymptotic densities represent (section \ref{sec:derivatives}) and show that all \tgl{} boundaries can be reconstructed (modulo translations and rotations) given sufficient information about the nonasymptotic density and its derivatives (section \ref{sec:reconstruct-T} and appendix \ref{sec:appen}).

The main contribution of this paper is to show (under our \tcgl{} condition) that all polygons (theorem \ref{thm:tcgl-poly-reconstruct} in section \ref{sec:polygon-no-tail}) and a $C^1$-dense set of $C^2$ boundaries (theorem \ref{thm:tcgl-reconstruct} in section \ref{sec:generic}) are reconstructible (modulo translations and rotations).
We briefly discuss and sketch the proofs of these two theorems.

{
\renewcommand{\thethm}{\ref{thm:tcgl-poly-reconstruct}}
\begin{thm}
For a polygon $\shape$ which is \tcgl{} with radius $r$, suppose that we have the integral area invariant $g(s,r)$ where $s$ is parameterized by arc length.
Suppose that for all $s$ we know $g(s,r)$ and its first derivatives with respect to $r$ (disk radius) and $s$ (position along the boundary).
This information is sufficient to completely determine $\shape$ up to translation and rotation; that is, we can recover the side lengths and angles of $\shape$.
\end{thm}
\addtocounter{thm}{-1}
}

The proof of this theorem uses the discontinuities in the $s$ derivative to determine the locations of vertices (and thus the side lengths between them).
We combine the $r$ derivative and the one-sided $s$ derivative information when centered on a vertex to recover the angles at which the polygon enters and exits the circle (which might not be the polygon vertex angle if the circle contains another vertex).
Doing this with the other one-sided $s$ derivative gives the same thing but using the orientation determined by the other polygon side incident to the vertex.
The combination of these yields the polygon's angle at each vertex.

{
\renewcommand{\thethm}{\ref{thm:tcgl-reconstruct}}
\begin{thm}
  Define $\Bbb{G} \equiv \{ \gamma| \gamma$ is a $C^2$ simple closed
  curve and \tgl{} for $r=\hat{r}\}$.  Suppose that, for $r=\hat{r}$, for
  all $s\in[0,L]$, and for each $\gamma\in\Bbb{G}$, we know the first-,
  second-, and third-order partial derivatives of $g_{\gamma}(s,r)$.
  Then the set of reconstructible $\gamma\in\Bbb{G}$ is
  $C^1$ dense in $\Bbb{G}$ where reconstructability is modulo reparametrization, translation, and rotation.
\end{thm}
\addtocounter{thm}{-1}
}

The first part of the proof shows that the derivative information can be used to obtain the curvature.
However, it is not the curvature at the boundary point where the circle is centered but rather the curvature at each of the points where the boundary enters and exits the circle.
Although the Euclidean distance to these points is known, the arc length distances are not and can vary from point to point.
Thus the sequences of curvatures we obtain also lose the arc length parameterization of our area invariant.
The rest of the proof is concerned with finding the arc length distance from the center to the entry and exit points which effectively recovers the curvature for all points.
This relies on matching up the unique features of exit angle sequences with each other which in turn relies on the existence of unique maxima and minima in these sequences.
While this is not true in general, it can be arranged to be so by a suitable small perturbation of the boundary (which is why our result is one of density rather than for all shapes).

This is a theoretical paper about a measure that is useful in applications: we do
not pretend that the reconstruction techniques in our proofs are practically
useful. In fact, the reconstructions we use to show uniqueness would be
seriously disturbed by the noise that any practical application would
encounter.
We do, however, comment on some possible approaches to
reconstruction (section \ref{sec:numerics}) using the \textsc{OrthoMads} direct search algorithm\cite{AbAuDeLe2009} to successfully reconstruct shapes which are not predicted by our theory.

\section{Notation and Preliminaries}
\label{sec:notation}

\begin{figure}[!]
\centering
\input{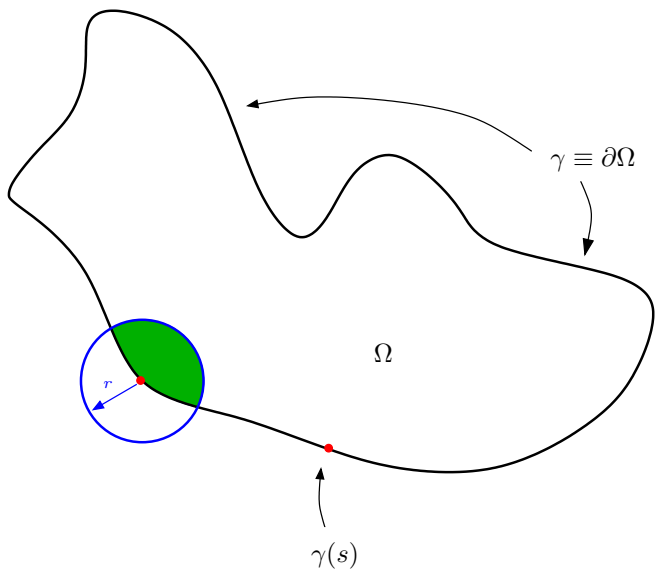}
\caption{Notation and basic setup}
\label{fig:shape-notation}
\end{figure}

Unless otherwise specified, we will be assuming throughout this paper that $\shape \subset \mathbb{R}^2$ is a compact set with simple closed, piecewise continuously differentiable boundary $\bd$ of length $L$.
Let $\gamma : [0, L] \rightarrow \bd$ be a continuous arclength parameterization of $\bd$ (see Figure~\ref{fig:shape-notation}).
We will adopt the convention that $\gamma$ traverses $\bd$ in a counterclockwise direction so it always keeps the interior of $\shape$ on the left (there is no compelling reason for this particular choice, but adopting a consistent convention allows us to avoid some ambiguities later).
Note that $\gamma(0) = \gamma(L)$ and that $\gamma$ restricted to $[0,L)$ is a bijection.
Denote by $\disk{p}{r}$ the closed disk and $\cir{p}{r}$ the circle of radius $r$ centered at the point $p \in \mathbb{R}^2$.

In geometric measure theory, the $m$-dimensional density of a set $A \subseteq \mathbb{R}^n$ at a point $p \in \mathbb{R}^n$ is given by
\[
\Theta^m(A, p) = \lim_{r\downarrow 0}\frac{\mathcal{H}^m(A \cap \disk{p}{r})}{\alpha_m r^m}
\]
where $\mathcal{H}^m$ is the $m$-dimensional Hausdorff measure
and $\alpha_m$ is the volume of the unit ball in $\mathbb{R}^m$\cite{morgan-gmt}.
In the current context, the $2$-dimensional density of $\shape$ at $\gamma(s)$ is simply
\[
\Theta^2(\shape, \gamma(s)) = \lim_{r\downarrow 0}\frac{\area(\shape \cap \disk{\gamma(s)}{r})}{\pi r^2}.
\]

While we can evaluate this for all $s \in [0,L)$, just knowing the density at every point along the boundary is generally insufficient to reconstruct the original shape.
If $\gamma'(s)$ exists, then $\area(\shape \cap \disk{\gamma(s)}{r})$ is approximated arbitrarily well for sufficiently small $r$ by replacing $\bd$ with its tangent line (which gives us an area of exactly $\frac{\pi r^2}{2}$).
Hence, we have $\Theta^2(\shape, \gamma(s)) = \frac{1}{2}$ at any point where $\gamma$ is differentiable.
That is, just knowing $\Theta^2$ (i.e., the limit) is insufficient to distinguish any two shapes with $C^1$ boundary.

Contrast this with the situation where we know $\area(\shape \cap \disk{\gamma(s)}{r})$ for every $s \in [0,L)$ and $r > 0$ (i.e., we have all of the values needed to compute the limit as well).
This added information is sufficient to uniquely identify $C^2$ curves by recovering their curvature at every point (see Appendix \ref{sec:appen}).

One natural question to ask (and the focus of the present work) is whether failing to pass to the limit (i.e., using some fixed radius $r$ instead of the limit or all $r > 0$) and collecting the values for all points along the boundary preserves enough information to reconstruct the original shape.
That is, can a nonasymptotic density (perhaps along with information about its derivatives) be used as a signature for shapes?

\subsection{Definitions}

\begin{dfn}
\label{def:geomeasure}
In the current context, the integral area invariant\cite{manay-2004-1} is denoted by $g : [0, L) \times \mathbb{R}^+ \rightarrow \mathbb{R}^+$ and given by
\[
g(s, r) = \int_{\disk{\gamma(s)}{r}\cap\shape} \,dx = \area(\shape \cap \disk{\gamma(s)}{r}).
\]
\end{dfn}

\begin{rem}
Note the lack of the normalizing factor $\pi r^2$ in the definition of $g(s,r)$.
Since we presume that $r$ is fixed and known for the situations we study, it's trivial to convert data between the forms $g(s,r)$ and $\frac{g(s,r)}{\pi r^2}$;
we choose to leave out the normalizing factor in the definition of $g(s,r)$ as it is the integral area invariant of Manay et al.\cite{manay-2004-1} and this form proves useful when computing derivatives in section \ref{sec:derivatives}.
\end{rem}

We introduce the \tgl{} condition as a simplifying assumption for the shapes we consider.

\begin{dfn}
\label{def:gl}
\label{def:tgl}
For a fixed radius $r$, we say that $\bd$ is graph-like (GL) at a point $p \in \bd$ (or graph-like on $\disk{p}{r}$) if it is possible to impose a Cartesian coordinate system such that the set of points $\bd \cap \disk{p}{r}$ is the graph of some function $f$ in this coordinate system.
Without loss of generality, we adopt the convention that $p$ is the origin so that $f(0) = 0$.
We define \tgl{} (TGL) in the same way but further require that $\bd$ be continuously differentiable and $f'(0) = 0$ (noting that $f$ is $C^1$ because $\bd$ is).
This is illustrated in figure \ref{fig:tgl}.
Without loss of generality (and in keeping with our convention that $\gamma$ traverses $\bd$ counterclockwise), we assume that the interior of $\shape$ is ``up'' in the circle (i.e., that $(0,\epsilon) \in \shape$ for sufficiently small $\epsilon > 0$).
If $\bd$ is (tangentially) graph-like on $\disk{p}{r}$ for all $p \in \bd$, we say that $\bd$ is (tangentially) graph-like for radius $r$.
\end{dfn}
\begin{figure}
\makebox[\textwidth][c]{
\subfigure[]{
\input{TGL.pdf_t}
\label{fig:tgl}
}
\subfigure[]{
\input{TCGL.pdf_t}
\label{fig:tcgl}
}}
\caption{\subref{fig:tgl} Tangentially and \subref{fig:tcgl} tangent cone graph-like}
\end{figure}
\begin{figure}
\subfigure[]{
\label{fig:glfail}
\begin{tikzpicture}[scale=.9]
\draw[thick=2pt,fill=teal!40] (-2, 4) -- (-2,0) -- (2,0) -- (2,4) -- cycle;
\draw[color=blue] (0,0) circle (3);
\node[label=below:{$\gamma(s)$},style=circle,fill=red,minimum size=4pt,inner sep=0pt] () at (0,0) {};
\end{tikzpicture}}
\subfigure[]{
\label{fig:tglfail}
\begin{tikzpicture}[scale=.9]
\draw [thick=2pt,rounded corners=5mm,fill=teal!40] (-4,0)--(4,0)--(4,2)--(-4,2)--cycle;
\draw[color=blue] (0,0) circle (3);
\node[label=below:{$\gamma(s)$},style=circle,fill=red,minimum size=4pt,inner sep=0pt] () at (0,0) {};
\end{tikzpicture}}
\caption[Graph-like examples]{\subref{fig:glfail} The square is not graph-like with the indicated radius (no orientation makes it a graph).
\subref{fig:tglfail} The rounded rectangle is graph-like but not \tgl{} with the indicated center and radius.}
\end{figure}

It is instructive to consider what is {\em not} graph-like or \tgl{}.
Violations of the graph-like condition are generally due to a radius that is too large (certainly, choosing a radius so large that all of $\shape$ is in the disk will do it).
For example, a unit side length square is not graph-like with radius $\frac{1}{2} + \epsilon$ for any $\epsilon > 0$ (position the circle at the center of a side; see figure \ref{fig:glfail}).
Notice that the same square is graph-like with any radius $\frac{1}{2}$ or below.
A shape can fail to be \tgl{} while still being graph-like if it fails to be a graph in the required orientation but works in some other (see figure \ref{fig:tglfail}).

We would like to consider shapes with corners but our \tgl{} condition requires that the boundary be differentiable everywhere.
The following definitions allow us to generalize the \tgl{} condition to this situation by using one-sided derivatives.

\begin{dfn}
\label{def:tcone}
Given a piecewise $C^1$ function $\gamma : [0,L] \rightarrow \mathbb{R}^2$, we define the tangent cone of $\gamma$ at a point $s$ (which we denote by $\tc{\gamma}{s}$) in terms of the one-sided derivatives.
In particular, we let $\tc{\gamma}{s} = \{\alpha \Gamma^- + \beta \Gamma^+ \mid \alpha,\beta \geq 0, \alpha + \beta > 0\}$ where $\Gamma^- = \lim_{t\uparrow s} \gamma'(t)$ and $\Gamma^+ = \lim_{t\downarrow s} \gamma'(t)$.
\end{dfn}

\begin{dfn}
\label{def:tcgl}
We extend the \tgl{} notion to boundaries that are piecewise $C^1$ by defining $\bd$ to be \tcgl{} (TCGL) at a point $\gamma(s) \in \bd$ if it is graph-like at $\gamma(s)$ for every orientation in the tangent cone of $\bd$ at $s$.
More precisely, for every 
$w \in \tc{\gamma}{s}$ and every pair of distinct points $u,v\in\bd\cap\disk{p}{r}$, we have $\langle w, u-v \rangle \neq 0$ (see figure \ref{fig:tcgl}).
\end{dfn}

\begin{rem}
It is clear that $\tc{\gamma}{s}$ in definition \ref{def:tcone} is a convex cone.
The tangent cone is dependent on the direction in which $\gamma$ traverses $\bd$ (which by convention was counterclockwise) since an arc-length traversal $\hat{\gamma}(s,r) = \gamma(L-s,r)$ would have different tangent cones (namely, $w \in \tc{\gamma}{s}$ iff $-w \in \tc{\hat{\gamma}}{s}$).
However, these differences are irrelevant to the application of definition \ref{def:tcgl}.
\end{rem}

\begin{rem}
Note that when $\bd$ is $C^1$, there is only one direction in $\tc{\gamma}{s}$ for each $s$ (i.e., the tangent to $\bd$ at $\gamma(s)$).
Thus, the definitions of \tgl{} and \tcgl{} coincide when $\bd$ is $C^1$ and every \tgl{} boundary is \tcgl{}.
\end{rem}

\subsection{Two-Arc Property}
\label{sec:tgl_2arc}

The graph-like condition implies (in proof of the following lemma) that $\shape$ will never be entirely contained in the disk, no matter where on the boundary we center it.
That is, some part of $\shape$ lies outside of $\disk{p}{r}$ for every $p \in \bd$.

\begin{lem}
\label{2intlemma}
Let $r \in \mathbb{R}^+$ and $p \in \bd$. If $\bd$ is graph-like on $\disk{p}{r}$, then 
$|\bd \cap \cir{p}{r}| \geq 2$.
\begin{proof}
Suppose by way of contradiction that $|\bd \cap \cir{p}{r}| < 2$.
Since $\bd$ is a simple closed curve, we have $\bd \subseteq \disk{p}{r}$.
As $\bd$ is graph-like at $p$ with radius $r$, there exists some orientation for which $\bd \cap \disk{p}{r} = \bd$ is the graph of a well-defined function.
However, $\bd$ is a simple closed curve so it is not the graph of a function in any orientation, yielding a contradiction.
\end{proof}
\end{lem}

The next result is the reason we find the \tcgl{} condition useful.
It says that if $\bd$ is \tcgl{} with radius $r$, then, for every $p \in \bd$, the disk $\disk{p}{r}$ has only two points of intersection with $\bd$ and these are transverse.
In other words, this means that when working locally in the disk $\disk{p}{r}$ we need only consider a single piece of $\bd$.

\begin{thm}
\label{thm:two-intersections}
If $\bd$ is \tcgl{} with radius $r \in \mathbb{R}^+$ at $p \in \bd$,
then $|\bd \cap \cir{p}{r}| = 2$ and $\bd$ crosses $\cir{p}{r}$ transversely at these points.
As a result, for every $q_1, q_2 \in \bd \cap \disk{p}{r}$, there is a unique arc along $\bd$ between them in $\disk{p}{r}$.
\begin{proof}
By Lemma \ref{2intlemma}, we have that $|\bd \cap \cir{p}{r}| \geq 2$.
Note that $\bd$ contains an interior point ($p$) and at least two boundary points of the disk $\disk{p}{r}$ (since $|\bd \cap \cir{p}{r}| \geq 2$).
As $\bd$ is connected and simply closed, there must exist an arc of $\bd$ within the disk going from some point on $\cir{p}{r}$ through $p$ to another point on $\cir{p}{r}$.

Suppose $|\bd \cap \cir{p}{r}| > 2$; that is, there are other points of intersection.
Letting $q$ denote one of these, there are two cases to consider (illustrated in Figure \ref{fig:nointersect}).
\begin{figure}
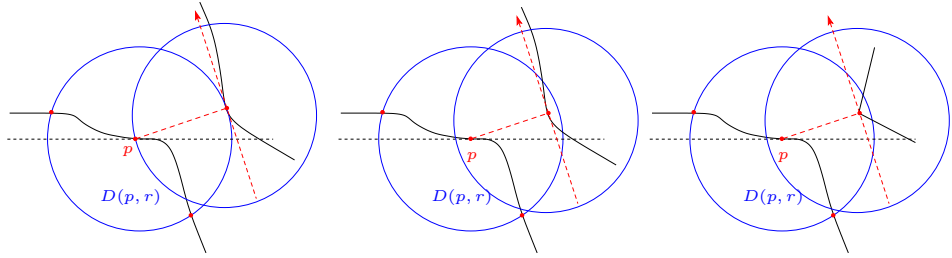

\centering
\subfigure{
\input{TCGLto2Arc1b.pdf_t}
}
\input{TCGLto2Arc2b.pdf_t}
\input{TCGLto2Arc3b.pdf_t}
\caption{Additional points of intersection violate the TCGL condition.}
\label{fig:nointersect}
\end{figure}
\begin{enumerate}
\item [(a)]
\label{case:nocross}
$\bd$ does not cross $\cir{p}{r}$ at $q$.

As $\bd$ is \tcgl{} at $q$, then $\bd\cap\cir{q}{r}$ is a graph in every orientation in the tangent cone of $\bd$ at $q$.
In particular, note that the tangent line to $\cir{p}{r}$ at $q$ is in this cone.
However, the line from $p$ to $q$ is normal to this line and thus $\bd\cap\cir{q}{r}$ is not graph-like in this orientation, a contradiction.
Therefore, this case cannot occur.
This argument applies to all points in $\bd \cap \cir{p}{r}$ so we immediately have the result that $\bd$ always crosses $\cir{p}{r}$ transversely.

\item [(b)]
$\bd$ crosses $\cir{p}{r}$ at $q$.

There exists $q' \in \bd \cap \cir{p}{r}$ such that there is a path along $\bd$ in $\disk{p}{r}$ from $q$ to $q'$.
That is, there exist $s_1, s_2 \in [0,L)$ (without loss of generality, $s_1<s_2$) such that $\gamma(s_1) = q$, $\gamma(s_2) = q'$ and the image of $[s_1,s_2]$ under $\gamma$ is contained in $\disk{p}{r}$ (but does not include $p$, since it is on another arc and $\bd$ is simple).
Thus $\gamma$ enters $\cir{p}{r}$ at $s_1$ and exits at $s_2$.

If we can find $s \in [s_1, s_2]$ and $w$ in the tangent cone of $\bd$ at $\gamma(s)$ satisfying $\langle w, p - \gamma(s) \rangle  = 0$, we will contradict that $\bd$ is \tcgl{}.

Define $v : [s_1,s_2] \rightarrow \mathbb{R}^2$ by
\[
v(s) = \begin{cases}
\lim_{t\downarrow s_1} \gamma'(s),& s=s_1,\\
\lim_{t\uparrow s} \gamma'(s),& s \in (s_1, s_2].
\end{cases}
\]
Note that $v(s)$ is in the tangent cone of $\bd$ at $\gamma(s)$ so that $\bd\cap\disk{\gamma(s)}{r}$ is graph-like using the orientation given by $v(s)$.

Define $\phi(s) : [s_1, s_2] \rightarrow \mathbb{R}$ by $\phi(s) = \langle v(s), p - \gamma(s)\rangle$.
Note that from $\gamma(s_1)$ both $v(s_1)$ and $p - \gamma(s_1)$ are directions pointing into the circle so $\phi(s_1) > 0$.
Similarly, $v(s_2)$ points out and $p - \gamma(s_2)$ points in so that $\phi(s_2) < 0$.

Observe that $v$ (and therefore $\phi$) is piecewise continuous since $\gamma$ is piecewise $C^1$.
By a piecewise continuous analogue of the intermediate value theorem, there exists $\bar{s} \in [s_1, s_2]$ such that
\[
\lim_{t\rightarrow \bar{s}^-} \phi(t) \leq 0 \leq \lim_{t\rightarrow \bar{s}^+} \phi(t).
\]

By continuity of the inner product and $\gamma$, we have
\[
\lim_{t\rightarrow \bar{s}^-} \phi(t) %
= \langle \lim_{t\rightarrow \bar{s}^-} \gamma'(t), p-\gamma(\bar{s})\rangle.
\]
Similarly, $\lim_{t\rightarrow \bar{s}^+} \phi(t) = \langle \lim_{t\rightarrow \bar{s}^+} \gamma'(t), p-\gamma(\bar{s})\rangle$

If $\gamma$ is differentiable at $\bar{s}$, then $\phi(\bar{s}) = \lim_{t \rightarrow\bar{s}}\phi(t) = 0$ and we have our contradiction.
Otherwise, let $w_1 = \lim_{t\rightarrow \bar{s}^-} \gamma'(t)$ and $w_2 = \lim_{t\rightarrow \bar{s}^+} \gamma'(t)$.
As both $w_1$ and $w_2$ are in the convex tangent cone of $\bd$ at $\gamma(\bar{s})$, any positive linear combination of them is as well.
Letting $\psi(\lambda) = \lambda w_1 + (1-\lambda)w_2$, we have 
\[
\langle \psi(0), p - \gamma(\bar{s})\rangle \leq 0 \leq \langle \psi(1), p-\gamma(\bar{s})\rangle.
\]
Noting that $\psi$ is continuous in $\lambda$, we apply the intermediate value theorem to obtain $\bar{\lambda} \in (0,1)$ such that $\langle \psi(\bar{\lambda}), p-\gamma(\bar{s})\rangle = 0$.
Letting $w = \psi(\bar{\lambda})$, we obtain our contradiction.
\end{enumerate}
Therefore, there are no other points of intersection and $|\bd \cap \cir{p}{r}| = 2$.
\end{proof}
\end{thm}

\begin{dfn}
We say that $\shape$ has the two-arc property for a given radius $r$ if for every point $p \in \bd$, we have that $\disk{p}{r}$ divides $\bd$ into two connected arcs: $\bd \cap \disk{p}{r}$ and $\bd \backslash \disk{p}{r}$.
Instead of considering how $\disk{p}{r}$ divides $\bd$, we can equivalently frame the definition in terms of how $\bd$ divides $\cir{p}{r}$.
That is, $\shape$ has the two-arc property if the circle $\cir{p}{r}$ is divided into two connected arcs by $\bd$ for every $p \in \bd$.  \end{dfn}

\begin{cor}
If $\shape$ is \tcgl{} for some radius $r$, then it has the two-arc property.
\label{lem:tcgl-2arc}
\begin{proof}
This is a trivial consequence of Theorem \ref{thm:two-intersections}.
\end{proof}
\end{cor}

\begin{cor}
If $\shape$ is \tgl{} for some radius $r$, then it has the two-arc property for radius $r$.
\label{cor:tgl-2arc}
\end{cor}

\begin{figure}[htp!]
  \centering
  \input{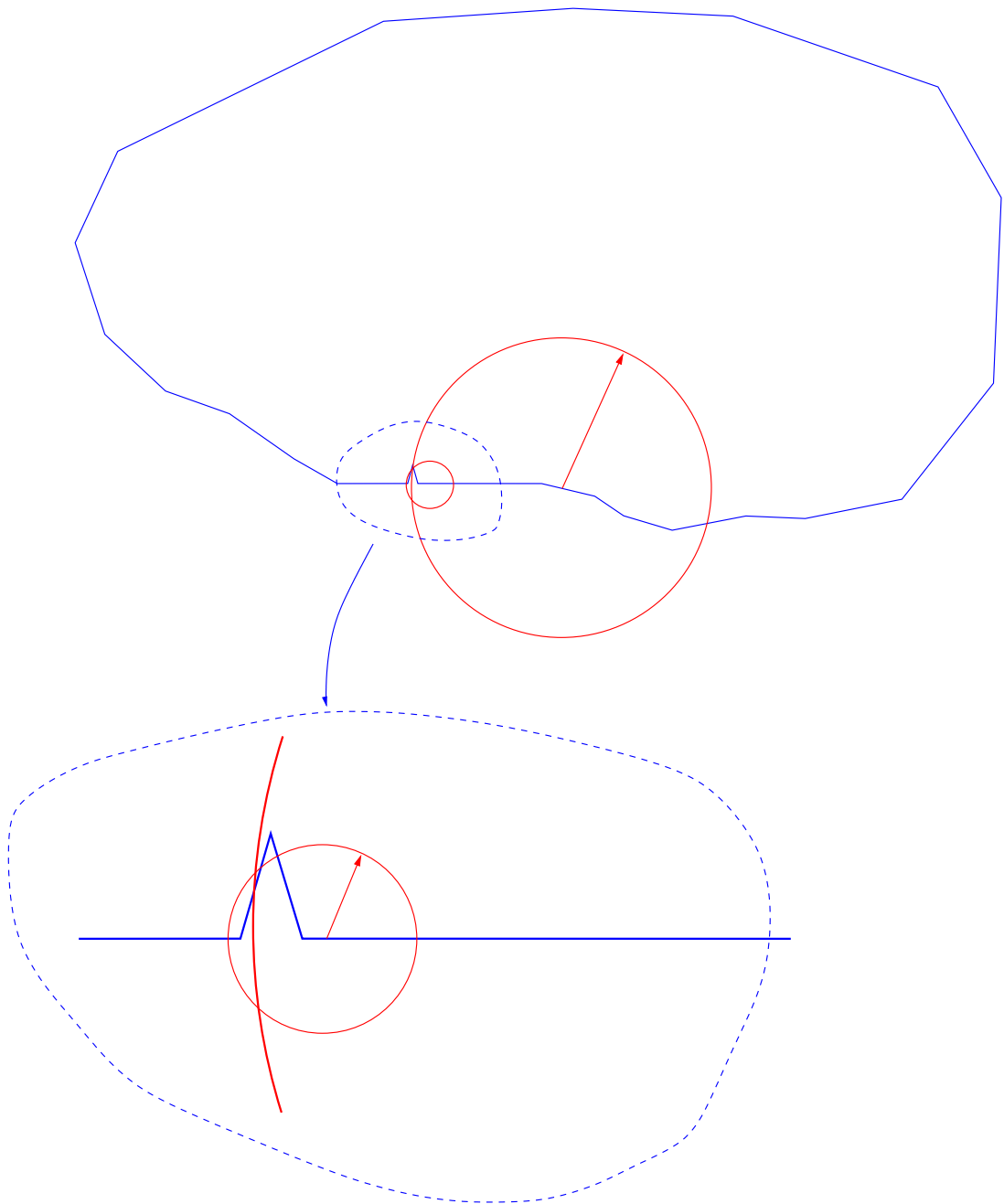}
  \caption{The two-arc property for $r = \hat{r}$
    does not imply that it holds for all $r < \hat{r}$}
 \label{fig:tgl_counter_ex1}
\end{figure}
\begin{rem}
\label{rem:2arc}
  While the assumption of the two-arc property for disks of radius $r =
  \hat{r}$ does not imply the two-arc property for all $r < \hat{r}$
  (see Figure~\ref{fig:tgl_counter_ex1}), it is the case that TGL for
  $r = \hat{r}$ does imply that $\gamma$ is TGL for all $0 < r <
  \hat{r}$. The fact that $\gamma$ is TGL for all $0 < r < \hat{r}$ follows
  easily from the definition of TGL and the fact that $\disk{p}{r} \subsetneq \disk{p}{\hat{r}}$.
\end{rem}

\subsection{Notation}
\begin{figure}[htp!]
  \centering
  \input{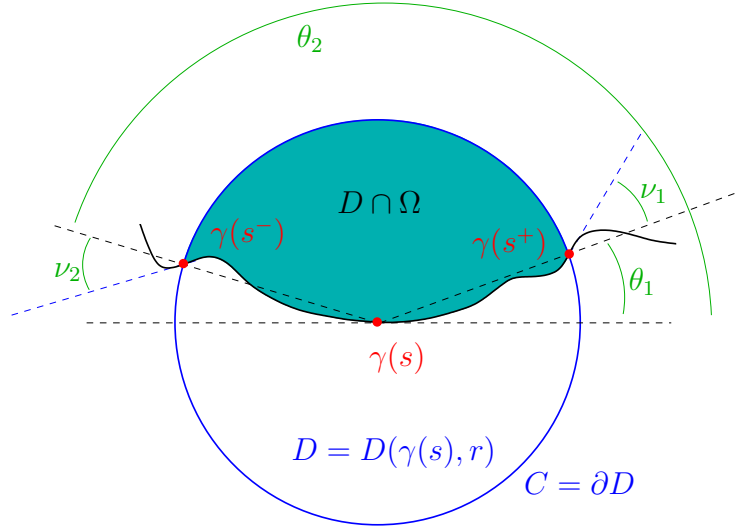}
  \caption{Notation and conventions}
\label{fig:step0}
\label{fig:notation}
\end{figure}

Suppose that $\bd$ is \tcgl{} with radius $r$ and we have some $s \in [0, L)$ such that $\bd$ is \tgl{} at $\gamma(s)$ with radius $r$.
Since $\bd$ is TGL at $\gamma(s)$, it has two points of intersection with $\cir{\gamma(s)}{r}$ by theorem \ref{thm:two-intersections}.
In the orientation forced by the TGL condition, one of these points of intersection must be on the right side of the circle and one must be on the left side.

With reference to figure \ref{fig:step0} we define $s^{+}(s)$ and $s^{-}(s) \in [0,L)$ so that $\gamma(s^{+}(s))$ is the point of intersection on the right and $\gamma(s^{-}(s))$ is the point of intersection on the left.
The notation is motivated by the fact that $0 < s^{-}(s) < s < s^{+}(s) < L$ in general due to our convention that $\gamma$ traverses $\bd$ counterclockwise.
The only case where this is not true is when $\gamma(L) = \gamma(0)$ is in the disk but even then it will hold for a suitably shifted $\hat{\gamma}$ that starts at some point outside the current disk.

The quantities $\theta_1(s)$ and $\theta_2(s)$ are the angles that the rays from the origin to the right and left points of intersection, respectively, make with the positive $x$ axis.
We can assume $\theta_1(s) \in (-\frac{\pi}{2},\frac{\pi}{2})$ and $\theta_2(s) \in (\frac{\pi}{2},\frac{3\pi}{2})$.

We define $\nu_1(s)$ as the angle between the vector $\gamma(s^{+}(s))-\gamma(s)$ and the vector $\lim_{t\downarrow s^{+}(s)}\gamma'(t)$, the one-sided tangent to $\bd$ at the point of intersection on the right.
That is, we are measuring the angle between the outward normal to the disk at the point of intersection and the actual direction $\gamma$ is going as it exits the disk.
We define $\nu_2(s)$ similarly.
We have $\nu_1, \nu_2 \in (-\frac{\pi}{2},\frac{\pi}{2})$ due to the fact that all circle crossings are transverse by theorem \ref{thm:two-intersections}.

When the proper $s$ to use is implied by context, we will often simply write $s^{+}$, $s^{-}$, $\theta_1$, $\theta_2$, $\nu_1$ and $\nu_2$ in place of $s^{+}(s)$, $s^{-}(s)$, and so forth.
\subsection{Calculus on Tangent Cones}

The following result is a version of the intermediate value theorem for elements of the tangent cones.
\begin{lem}
\label{lem:intvalue}
Suppose $\bd$ is \tcgl{} on $\disk{\gamma(s)}{r}$ and $s_1 < s_2$ such that $\gamma(s_1), \gamma(s_2) \in \disk{\gamma(s)}{r}$.
Further suppose that $w_1 \in \tc{\gamma}{s_1}$, $w_2 \in \tc{\gamma}{s_2}$, $\alpha \in (0, 1)$, and let $w' = \alpha w_1 + (1-\alpha)w_2$.
Then, there exists $s' \in [s_1, s_2]$ such that either $w'$ or $-w'$ is in $\tc{\gamma}{s'}$.
\begin{proof}
Let $n$ be a unit vector in $\mathbb{R}^2$ with $n \perp (\alpha w_1 + (1-\alpha)w_2)$.
We have $\alpha\langle n, w_1 \rangle = -(1-\alpha)\langle n, w_2 \rangle$.
It suffices to consider only $\langle n, w_1 \rangle \leq 0 \leq \langle n, w_2 \rangle$ as the argument is identical in the other case.
Note that since $0 \leq \langle n, w_2 \rangle = c_1\langle n, \lim_{t\uparrow s_2}\gamma'(t)\rangle + c_2\langle n, \lim_{t\downarrow s_2}\gamma'(t)\rangle$ for some nonnegative constants $c_1, c_2$ not both zero, at least one of the inner products on the right is nonnegative.
Using the notation of definition \ref{def:tcone}, we define $M_2 = \argmax_{\Gamma \in \{\Gamma^+, \Gamma^-\}} \langle n, \Gamma\rangle$ and have $\langle n, M_2 \rangle \geq 0$.
We similarly define $M_1$ with respect to $w_1$ such that $\langle n, M_1 \rangle \leq 0$.%

Define 
\[
v(t) = \begin{cases}
M_i, & t = s_i, i = 1, 2\\
\lim_{t\uparrow t} \gamma'(s)
\end{cases}
\]
and $\phi(t) = \langle n, v(t) \rangle$.
Since $\phi(s_1) \leq 0 \leq \phi(s_2)$, the argument proceeds as in theorem \ref{thm:two-intersections} to yield $\bar{s} \in [s_1, s_2]$ and $\bar{w} \in \tc{\gamma}{\bar{s}}$ such that $\langle n, \bar{w} \rangle = 0$.
Thus $\bar{w} = kw'$ for some $k \neq 0$.
In particular, $w' = \frac{1}{k}\bar{w}$ so either $w' \in \tc{\gamma}{\bar{s}}$ or $-w' \in \tc{\gamma}{\bar{s}}$ (depending on the sign of $k$).
\end{proof}
\end{lem}

In addition to the intermediate value theorem, we have an analogous mean value theorem for tangent cone elements.

\begin{lem}
\label{lem:meanvalue}
Suppose $\gamma:[a,b] \rightarrow \mathbb{R}^2$ is a simple, arc-length parameterized curve with piecewise continuous derivative defined on $(a,b)$ except possibly on finitely many points.
Further suppose that the image of $\gamma$ has no cusps.
Then there exists $c$ in $(a,b)$ such that either $\gamma(b)-\gamma(a)$ or $-(\gamma(b)-\gamma(a))$ is in $T_\gamma(c)$.
\begin{proof}
Let $n$ be a unit vector with $\langle \gamma(b)-\gamma(a), n \rangle = 0$.
Consider $\psi(t) = \langle \gamma(t) , n \rangle$ and note that $\psi'(t) = \langle \gamma'(t), n \rangle$ is defined wherever $\gamma(t)$ is differentiable.
We have $\int_a^b \psi'(t) = \psi(b) - \psi(a) = \langle \gamma(b) - \gamma(a), n \rangle = 0$.
Thus, either $\psi'(t) = 0$ everywhere it is defined or it takes on both positive and negative values.
In particular, there exists a point $c \in (a,b)$ such that either $\psi'(c) = 0$ or $\lim_{t\uparrow c} \psi'(t) \leq 0 \leq \lim_{t\downarrow c} \psi'(t)$.

If $\psi'(c) = 0$, then we have $\langle \gamma'(c), n \rangle = 0$ so that $\gamma'(c) = k(\phi(b) - \phi(a))$ for some $k \neq 0$.
As $\gamma'(c) \in \tc{\gamma}{c}$, we have $\frac{k}{|k|}(\phi(b) - \phi(a)) \in \tc{\gamma}{c}$ which gives us our conclusion.

If $\lim_{t\uparrow c} \psi'(t) \leq 0 \leq \lim_{t\downarrow c} \psi'(t)$, there exists $\alpha \in (0, 1)$ such that $0 = \alpha \lim_{t\uparrow c} \psi'(t) + (1-\alpha)\lim_{t\downarrow c} \psi'(t)$.
Note that $\lim_{t\uparrow c} \psi'(t) = \langle w_1, n \rangle$ and $\lim_{t\downarrow c} \psi'(t) = \langle w_2, n \rangle$ for some $w_1, w_2 \in \tc{\gamma}{c}$ and let $w' = \alpha{w_1} + (1-\alpha){w_2}$.

By the convexity of $\tc{\gamma}{c}$, we have $w' \in \tc{\gamma}{c}$ with $\langle w', n\rangle = 0$ which follows as in the previous case.
\end{proof}
\end{lem}

The following lemma tells us that the \tcgl{} condition is sufficient to apply lemma \ref{lem:meanvalue}.

\begin{lem}
If $\bd$ is \tcgl{} for some radius $r$, then $\bd$ has no cusps.
\begin{proof}
Suppose $\bd$ has a cusp at $\gamma(s)$.
Then, using the terminology of definition \ref{def:tcone} and the fact that $\gamma$ is arc length parameterized, we have $\Gamma^+ = -\Gamma^-$.
We let $w = 0$ and note that $w = \Gamma^+ + \Gamma^- \in \tc{\gamma}{s}$.
Letting $u,v\in\bd\cap\disk{\gamma(s)}{r}$ with $u \neq v$, we have $\langle w, u-v\rangle = 0$, contradicting the fact that $\bd$ is \tcgl{}.
Therefore, $\bd$ has no cusps.
\end{proof}
\end{lem}

\subsection{TCGL Boundary Properties}

The following technical lemmas allow us to bound various distances and areas encountered in \tcgl{} boundaries.

\begin{lem}
Suppose that $\bd$ is \tcgl{} with radius $r$ and points $p_1, p_2 \in \bd$ with $d(p_1, p_2) < r$.
Then one of the arcs (call it $P$) along $\bd$ between $p_1$ and $p_2$ is such that, for any two points $q_1, q_2 \in P$, we have $d(q_1, q_2) < r$.
\begin{proof}
Note that $p_2 \in \disk{p_1}{r}$ so that there is an arc along $\bd$ from $p_1$ to $p_2$ which is fully contained in the interior of $\disk{p_1}{r}$ by theorem \ref{thm:two-intersections}.
We will call this arc $P$.

For all $x$ on $P$, let $P_x$ denote the subpath of $P$ from $p_1$ to $x$ (so $P = P_{p_2}$).
We claim that $P_x$ is contained in $\disk{x}{r}$ for all $x$ on $P$ (thus, $P$ is contained in $\disk{p_2}{r}$).
Indeed, if this were not the case, then there must be some $\hat{x}$ on $P$ such that $P_{\hat{x}}$ is contained in $\disk{\hat{x}}{r}$ but $\cir{\hat{x}}{r} \cap P_{\hat{x}}$ is nonempty (i.e., we can move the disk along $P$ until some part of the subpath hits the boundary).
That is, the subpath $P_{\hat{x}}$ has a tangency with the disk $\disk{\hat{x}}{r}$ which is impossible because of theorem \ref{thm:two-intersections}.

Let $q_1 \in P$ and note that since $P_x$ is contained in $\disk{x}{r}$ for all $x$ on $P$, we have that $P$ is contained in $\disk{q_1}{r}$.
Therefore, $d(q_1, q_2) < r$ for all $q_1, q_2 \in P$ as desired.
\end{proof}
\end{lem}

\begin{lem}
\label{lem:arcbound}
If $q_1 = \gamma(s_1), q_2 = \gamma(s_2) \in P$ where $P$ is as in the previous lemma, then the arc length between $q_1$ and $q_2$ along $P$ is at most $\sqrt{2}d(q_1, q_2)$.
\begin{proof}
Since $\shape$ is \tgl{}, for any $w_1 \in \tc{\gamma}{s_1}, w_2 \in \tc{\gamma}{s_2}$, the angle between $w_1$ and $w_2$ is at most $\frac{\pi}{2}$.
Since this is true for all $q \in P$, there is a point $q' = \gamma(s') \in P$ and $w' \in \tc{\gamma}{s'}$ such that the angle between $w'$ and tangent vectors for any other point $q \in P$ is at most $\frac{\pi}{4}$.

This means that $P$ is the graph of a Lipschitz function $g$ of rank 1 in the orientation defined by $w'$.
This does not necessarily imply that $\disk{q'}{r} \cap \bd$, $\disk{p_1}{r}\cap \bd$ or $\disk{p_2}{r}\cap \bd$ is the graph of a Lipschitz function; we explore a Lipschitz condition for the disks in section \ref{sec:tgl-ncgl}.
Let $x_1, x_2 \in [-r, r]$ with $p_1 = (x_1, g(x_1))$, $p_2 = (x_2, g(x_2))$.
Then the arclength from $p_1$ to $p_2$ is given by 
\[
\int_{x_1}^{x_2} \sqrt{1 + g'(x)^2}\,dx \leq \int_{x_1}^{x_2} \sqrt{2}\,dx = \sqrt{2}(x_2-x_1) \leq \sqrt{2}d(p_1, p_2).\qedhere
\]
\end{proof}
\end{lem}

\begin{lem}
\label{lem:areabound}
If $\gamma$ is \tcgl{} with radius $r$ and $0 \leq s_1 \leq s_2 < L$ with $d(\gamma(s_1), \gamma(s_2)) = \delta < r$, then the image of $[s_1, s_2]$ together with the straight line from $\gamma(s_1)$ to $\gamma(s_2)$ enclose a region with $O(\delta^2)$ area.
\begin{proof}
By Lemma \ref{lem:arcbound}, we have that the image of $[s_1, s_2]$ under $\gamma$ has arc length $s_2 - s_1 \leq \sqrt{2}\delta$.
Therefore, the region of interest has perimeter at most $(\sqrt{2}+1)\delta$ so by the isoperimetric inequality has area at most $\frac{(\sqrt{2}+1)^2}{4\pi}\delta^2$ from which the conclusion follows.
\end{proof}
\end{lem}

\section{TCGL polygonal approximations}
\label{sec:tgl-ncgl}

If $\shape$ is \tcgl{} with radius $r$, it can sometimes be nice to know that there is an approximating polygon to $\shape$ which is also \tcgl{}.
The following lemmas explore this idea.

\begin{lem}
If $\bd$ is TCGL with radius $r$ then for each $\epsilon \in (0, r)$, then there exists a polygonal approximation to $\bd$ that is TCGL with radius $r - \epsilon$ and such that every point on $\bd$ is within distance $\frac{\epsilon}{6}$ of the polygon.
\label{lem:tcgl-poly}
\begin{proof}
First, choose a finite number of points along the boundary such that the arc length along $\gamma$ between any two neighboring points is no more than $\frac{\epsilon}{3}$.
These will be the vertices of our polygon.
Similarly to $\gamma$, we let $\phi$ be an arclength parameterization of this polygon so that they both encounter their common points in the same order.

The fine spacing between vertices guarantees that we obtain the $\frac{\epsilon}{6}$ bound.
Indeed, given any point $p \in \bd$ and its neighboring vertices $v_1$ and $v_2$, the arc length along $\bd$ from $v_1$ to $p$ plus that from $p$ to $v_2$ is at most $\frac{\epsilon}{3}$ by assumption.
Since Euclidean distance is bounded above by arc length, we have $d(p, v_1) + d(p, v_2) \leq \frac{\epsilon}{3}$.
This bound in turn implies that at least one of $d(p,v_1)$ and $d(p,v_2)$ is bounded above by $\frac{\epsilon}{6}$.

Consider a point $p = \phi(t)$ on a side of the polygon (i.e., not a vertex) and its neighboring vertices $v_1 = \phi(t_1) = \gamma(s_1)$ and $v_2 = \phi(t_2) = \gamma(s_2)$ (chosen with $t_1 < t < t_2$ and $s_1 < s_2$).
By lemma \ref{lem:meanvalue}, there exists $s \in (s_1, s_2)$ such that $v_2 - v_1 \in T_\gamma(s)$.
Note that this is the only member of $T_\phi(t)$ up to positive scalar multiplication.

Combining the arcs along $\gamma$ and $\phi$ between $v_1$ and $v_2$, we obtain a closed curve with total length at most $\frac{2\epsilon}{3}$, so that the distance between any two points on the curve is at most $\frac{\epsilon}{3}$.
That is, for any $s' \in [s_1, s_2]$ and $t' \in [t_1, t_2]$, we have $d(\gamma(s'), \phi(t')) \leq \frac{\epsilon}{3}$.

Let $x \in \disk{\phi(t)}{r-\epsilon}$.
Then $d(x, \gamma(s)) \leq d(x, \phi(t)) + d(\phi(t), \gamma(s)) \leq r - \frac{2\epsilon}{3}$ so that $\disk{\phi(t)}{r-\epsilon}$ is contained in $\disk{\gamma(s)}{r - \frac{2\epsilon}{3}}$.

Let $a, b$ be distinct points on the polygon in $\disk{\phi(t)}{r-\epsilon}$ and consider the line connecting them.
This line also intersects $a', b'$ on $\gamma$ such that we have $a'\neq b'$, $d(a, a') \leq \frac{\epsilon}{3}$ and $d(b,b') \leq \frac{\epsilon}{3}$ so that $a', b' \in \bd \cap \disk{\gamma(s)}{r}$.
As $a - b = c(a'-b')$ for some scalar $c > 0$, we have
\[
\langle v_2 - v_1, a - b \rangle = c\langle v_2 - v_1, a' - b' \rangle \neq 0
\]
since $\gamma$ is TCGL at $\gamma(s)$ with radius $r$ and $v_2 - v_1 \in T_\gamma(s)$.
Thus $\phi$ is TCGL at $p$ with radius $r - \epsilon$.

The case where $p = \phi(t)$ is a vertex is similar but we must consider an arbitrary vector $w \in T_\phi(t)$ in the inner product.
We wish to show that, for every $w \in T_\phi(t)$, there is a $s'$ such that either $w$ or $-w \in T_\gamma(s')$ and $d(p, \gamma(s')) \leq \frac{\epsilon}{3}$, after which the proof follows as in the first case with $w$ (or $-w$) in place of $v_2-v_1$.
We let $\gamma(s) = \phi(t) = p$ and let $v_1 = \phi(t_1) = \gamma(s_1)$ and $v_2 = \phi(t_2) = \gamma(s_2)$ be the neighboring vertices (so $t_1 < t < t_2$ and $s_1 < s < s_2$).

As above, there exist $s_1', s_2'$ such that $s_1 \leq s_1' \leq s \leq s_2' \leq s_2$, $\gamma(s) - \gamma(s_1) \in \tc{\gamma}{s_1'}$ and $\gamma(s_2) - \gamma(s) \in \tc{\gamma}{s_2'}$.
Note that $\tc{\phi}{t}$ is exactly the set of positive linear combinations of these vectors.
By lemma \ref{lem:intvalue}, for every $w \in \tc{\phi}{t}$, there is a $s' \in [s_1', s_2']$ such that $w \in \tc{\gamma}{s'}$.
As $d(p, \gamma(s')) < \frac{\epsilon}{3}$, the proof is complete.
\end{proof}
\end{lem}
\begin{dfn}
We say that $\shape$ is tangentially graph-like and Lipschitz (TGLL) with radius $r$ if $\shape$ is \tgl{} with radius $r$ and there is some constant $0 < K < \infty$ such that for every $p \in \bd$, the arc $\disk{p}{r} \cap \bd$ is the graph of a Lipschitz function (in the same orientation used by the \tgl{} definition) and that the Lipschitz constant is at most $K$.
\end{dfn}

\begin{rem}
Note that \tgl{} does not imply tangentially graph-like and Lipschitz: taking
$\gamma$ to be a square with side length $5$ whose corners are
replaced by quarter circles of radius $1$ and then considering disks
of radius $\sqrt{2}$ centered on $\gamma$ yields one example.
\end{rem}

Because $\gamma$ is arclength parameterized by $s$, $||\dgm{s}|| = 1$
for all $s$.  Since $\gamma$ is assumed $C^1$ on its compact domain
$[0,L]$, $\dgmm$ is uniformly continuous: for any $\epsilon > 0$,
there is a $\delep$ such that if $|s_2-s_1| < \delep$ then
$||\dgm{s_2} - \dgm{s_1}||< \epsilon$.

We will use the fact that $\gamma$ always crosses $\partial D$ transversely to prove that $\gamma$ is in fact TGLL on slightly
bigger disks of radius $r+\delta$ as long as one takes a somewhat
bigger Lipschitz constant $\hat{K}$.
It is then an immediate result of lemma \ref{lem:tcgl-poly} that we can find an approximating polygon that is TCGL with radius $r$.

\begin{figure}[htp!]
  \centering
  \input{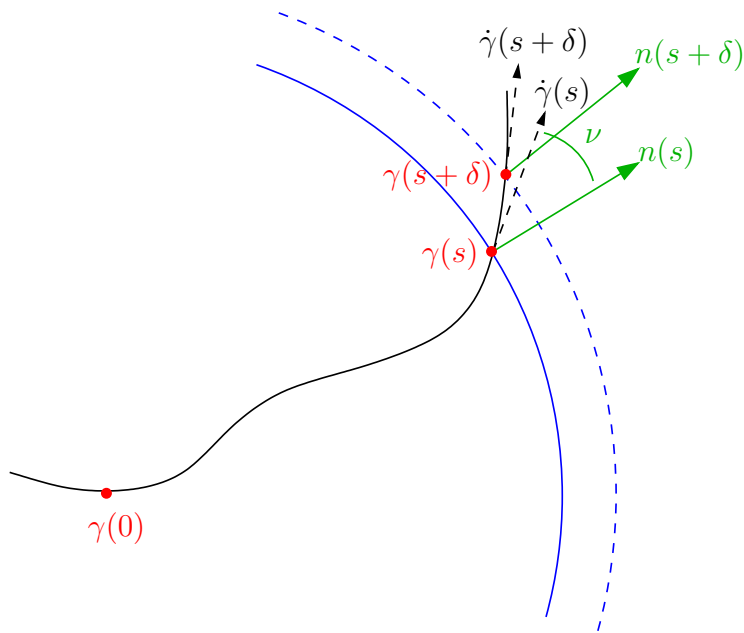}
  \caption{TGLL implies TCGL: Step one}
 \label{fig:step1}
\end{figure}

\begin{lem}
  If $\gamma$ is TGLL with radius $r$, then it is TGLL with radius $r+\delta$ for some $\delta > 0$ and there is an 
  approximating polygon $P_{\gamma}$ which is TCGL with radius $r$.
\end{lem}

\begin{proof}
\noindent{\bf Step 1: Show that the quantities $\nu_1$ and $\nu_2$
are continuous as a function of $s \in [0,L]$.(see Fig.~\ref{fig:step0})}

 Define $R^2(s,t)
\equiv ||\gm{s} - \gm{t}||^2$.  Taking the derivative, we get
\[DR = \left[ \left\langle \frac{\gm{s} - \gm{t}}{R(s,t)} ,
\dgm{s}\right\rangle , \left\langle \frac{\gm{t}-\gm{s}}{R(s,t)},
\dgm{t}\right\rangle \right]. \] Because $\nu_1$ and $\nu_2$ are both less than
$\pi/2$ and $\gamma$ is graph-like in the disk, we have that both
elements of this derivative are nowhere zero. By the implicit function
theorem, we get that $s^{-}(s)$ and $s^{+}(s)$ are continuous
functions of $s$. From this it follows that $\nu_1$ and $\nu_2$ are
continuous on $[0,L]$.

\noindent{\bf Step 2:} From the previous step and the compactness of $[0,L]$ we
get that $\nu_1(s)$ and $\nu_2(s)$ are both bounded by $M_\nu <
\pi/2$.
We define $\enu \equiv \pi/2 - M_{\nu} > 0$.
Fix a $t\in[0,L]$. Define $\rhh(s)$ by $\rhh^2(s) = R^2(s,t) =
||\gm{s} - \gm{t}||^2$. Then $\dot{\rhh}(s) = \langle \frac{\gm{s}
  -\gm{t}}{\rhh}, \dgm{s}\rangle = \langle n_t(s), \dgm{s}\rangle $
where $n_t(s)= \frac{\gm{s}-\gm{t}}{||\gm{s}-\gm{t} ||} =
\frac{\gm{s}-\gm{t}}{\rhh}$, the external normal to $\partial
D(\gamma(t),\rhh)$ at $\gm{s}$ (see Figure \ref{fig:step1}).  On any interval in $s$ where
$\dot{\rhh}(s) > 0$ we have that $\rhh(s)$ is one to one and strictly increasing.
Define $s^{*} \equiv s^{+}(t)$ and $s_{*} \equiv s^{-}(t)$. We showed
above that $\dot{\rhh}(s^*) = \langle n_t(s^*), \dgm{s^*}\rangle \geq
\cos(M_{\nu}) > 0$.

For $\langle n_t(s),\dgm{s}\rangle = 0$, $n_t(s)$ and $\dgmm$ will
have to have together turned by at least $\pi/2 - M_{\nu}$ radians.
And until they have turned this far, $\langle n_t(s),\dgm{s}\rangle >
0$. But $\dot{n}_t(s) \leq \frac{1}{\rho} \leq \frac{1}{r_{min}}$ for
some $r_{min} > 0$. (Choosing $r_{min} = \frac{r}{2}$ works.) And
$\dgmm$ is uniformly continuous on $[0,L]$. Therefore, there is a
$\delta_s$ such that on $[s^* , s^* + \delta_s]$, $n_t(s)$ and $\dgmm$
both turn by less than $\enu/3$.  Therefore, for $s \in [s^* , s^* +
\delta_s]$, we have that $\langle n_t(s), \dgm{s}\rangle > \cos(\pi/2
- \enu/3)$ and $\gamma([s^* , s^* + \delta_s))$ intersects $C =
\partial D(\gm{t},\rho)$ once for each $\rho \in [r,r+\delta_r]$, where
$\delta_r \equiv \delta_s\cos(\pi/2 - \enu/3)$.

A completely analogous argument works to show that $\gamma([s_*
-\delta_s, s_*])$ intersects $C = \partial D(\gm{t},\rho)$ once for
each $\rho \in [r,r+\delta_r]$. 

Define $d(t)$ to be the distance from $D(\gm{t},r)$ to
$\gamma\setminus\gamma([s_*-\delta_s,s^*+\delta_s])$. Since $\gamma$
is TGL, d(t) is greater than zero for all $t$ and is continuous in
$t$. Therefore, there is a smallest distance $\delta_d$ such that
$d(t) \geq \delta_d$ for all $t$. Define $\delta_{\gamma_o} =
\min(\delta_d/2,\delta_r/2)$.

Therefore, $\partial D(\gm{t},\rho)$ intersects $\gamma$ exactly
twice for $\rho \in [r,r+\delta_{\gamma_o}]$ for any $t\in[0,L]$. 

A similar argument shows that $\partial D(\gm{t},\rho)$ intersects
$\gamma$ exactly twice for $\rho \in [r-\delta_{\gamma_i},r]$ for any
$t\in[0,L]$. Defining $\delta_\gamma \equiv \min(\delta_{\gamma_i},
\delta_{\gamma_o})$ we get that $\partial D(\gm{t},\rho)$ intersects
$\gamma$ exactly twice for $\rho \in [
r-\delta_{\gamma},r+\delta_{\gamma} ]$, with the additional fact that
$\langle n_t(s), \dgm{s}\rangle > \cos(\pi/2
- \enu/3)$ at all those intersections.

\noindent{\bf Step 3:} TGLL implies that there is a constant $K <\infty$
such that $\gamma\cap\D(\gm{t},r)$ is the graph of a function whose
x-axis direction is parallel to $\dgm{t}$ and this function is
Lipschitz with Lipschitz constant $K$.

Since $\dgmm$ is uniformly continuous, there will be a $\delta_1$ such
that if $|u-v| < \delta_1$, then $\angle(\dgm{u}, \dgm{v}) <
\arctan{2K} - \arctan{K}$. Define $\delta_{K,s} = \min(\delta_s,
\delta_1)$. Define $\delta_{K,r} = \min(\delta_\gamma,
\delta_{K,s}\cos(\pi/2-\enu/3))$. Then
$\gamma\cap\D(\gm{t},r+\delta_{K,r})$ is the graph of a Lipschitz
function with Lipschitz constant at most $2K$ when $\dgm{t}$ is used
as the x-axis direction. That is, for all $t$, $\gamma$ is TGLL with
Lipschitz constant $2K$ for disks of radius $r+\delta_{K,r}$.
The result follows by lemma \ref{lem:tcgl-poly}.
\end{proof}

\section{Derivatives of $g(s,r)$}
\label{sec:derivatives}

\begin{figure}[htp!]
  \begin{center}
    \input{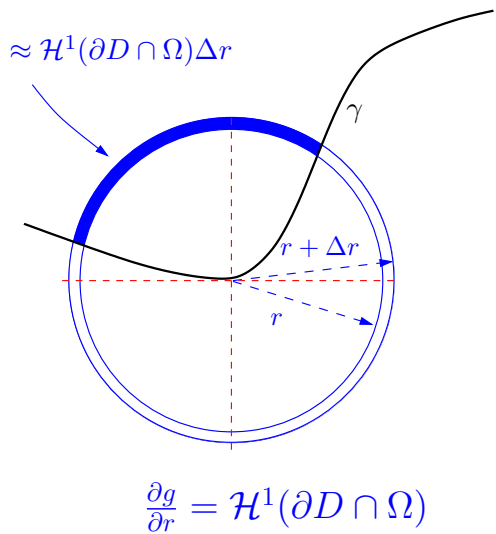}
  \end{center}
\caption{Deriving $\pd{g}{r}$ as the arclength of the circular segment.}
\label{fig:dgdr}
\end{figure}

\begin{lem}
\label{lem:dr}
Using the notation of figure \ref{fig:step0}, we have $\pd{}{r}g(s, r) = (\theta_2 - \theta_1)r$.
That is, the derivative exists and equals the length of the curve $\cir{\gamma(s)}{r} \cap \shape$.
\begin{proof}
We have (see figure~\ref{fig:dgdr})
\[
\pd{}{r}g(s,r) 
= \lim_{\Delta r \rightarrow 0}\frac{\area (\shape \cap \disk{\gamma(s)}{r+\Delta r})  - \area (\shape \cap \disk{\gamma(s)}{r})}{\Delta r}.
\]

This difference of areas can be modeled by the difference in the circular sectors of $\disk{\gamma(s)}{r+\Delta r}$ and $\disk{\gamma(s)}{r}$ with angle $\theta_1 - \theta_2$.
The actual area depends on the image of $\gamma$ outside of $\disk{\gamma(s)}{r}$, but this correction will be a subset of the circular segment of $\disk{\gamma(s)}{r+\Delta r}$ which is tangent to $\disk{\gamma(s)}{r}$ at the point $\gamma$ exits.
This has area $O(\Delta r^2)$ by lemma \ref{lem:areabound}.

Thus we have
\[
\pd{}{r}g(s,r) = \lim_{\Delta r \rightarrow 0} \frac{(\theta_1 - \theta_2)r\Delta r + \frac{1}{2}(\theta_1 - \theta_2)\Delta r^2 + O(\Delta r^2)}{\Delta r} = (\theta_1 - \theta_2)r.\qedhere
\]
\end{proof}
\end{lem}

\begin{figure}[htp!]
  \begin{center}
    \input{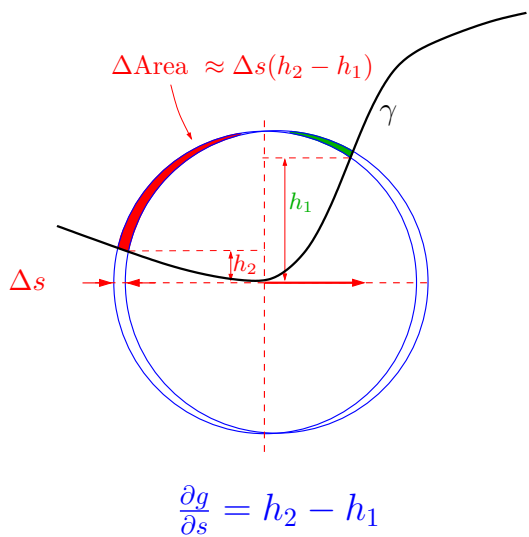}
  \end{center}
\caption{Deriving $\pd{g}{s}$ as the difference in heights of the entry and exit points}
\label{fig:dgds}
\end{figure}
\begin{lem}
Using the notation of figures \ref{fig:step0} and \ref{fig:dgds}, we have $\pd{}{s}g(s,r) = h_2 - h_1 = r\sin(\theta_2) - r\sin(\theta_1)$.
\label{lem:ds}
\begin{proof}
We have
\[
\pd{}{s} g(s,r) = \lim_{\Delta s \rightarrow 0} 
\frac{\area (\shape \cap \disk{\gamma(s+\Delta s)}{r})  - \area (\shape \cap \disk{\gamma(s)}{r})}{\Delta s}.
\]
The situation is illustrated in figure~\ref{fig:dgds} where we can see that the area being added as we go from $s$ to $s + \Delta s$ is the shaded region on the right with height $r - h_1$ and, considering first-order terms only, uniform width $\Delta s$ so has area $(r-h_1)\Delta s$.
Similarly, we are subtracting the area $(r-h_2)\Delta s$ on the left.
Therefore, we have
\[
\pd{}{s}g(s,r) = \lim_{\Delta s\rightarrow 0} \frac{(r-h_1)\Delta s - (r-h_2)\Delta s}{\Delta s} = h_2 - h_1.
\]
\end{proof}
\end{lem}

\section{Reconstructing shapes from T-like data}
\label{sec:reconstruct-T}

\begin{figure}[htp!]
  \centering
  \input{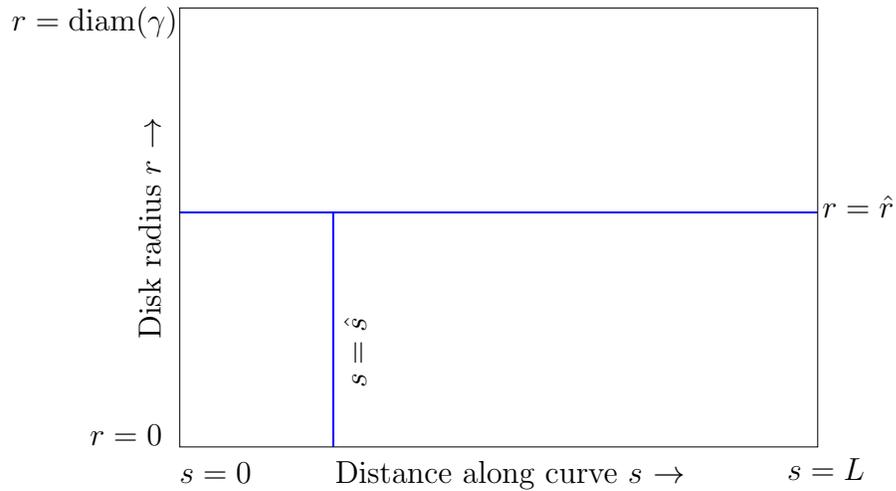}
  \caption[Restricting $g(s,r)$ to a T-like domain]{T-like data: we restrict the domain of $g(s,r)$ to a fixed radius
    $\hat{r}$ plus any vertical segment from $r=0$ to $r = \hat{r}$}
 \label{fig:t_like_data1}
\end{figure}

In this section, we consider the case where nonasymptotic densities and first derivatives are
known along a T-shaped set (i.e., for all $s$ with a fixed radius $\hat{r}$ and for all $r \leq \hat{r}$ with a fixed $\hat{s}$).
We show that this information is sufficient to guarantee reconstructability modulo reparametrizations, translations, and rotations.

\begin{lem}
  Assume that $\gamma$ is TGL for $\hat{r}$ (and thus all $r \leq \hat{r}$). Then if we know
  $g(s,r)$, $g_s(s,r) = \pd{g(s,r)}{s}$, and $g_r(s,r) =
  \pd{g(s,r)}{r}$ for $(s,r)\in([0,L]\times\{\hat{r}\}) \cup
  (\{\hat{s}\}\times (0,\hat{r}])$, we can reconstruct
  $\gamma(s)\in\Bbb{R}^2$ for all $s\in[0,L]$ modulo reparametrizations, translation, and rotations. (See
  figure~\ref{fig:t_like_data1}.)
\end{lem}

\noindent \emph{Proof:} As was shown in section~\ref{sec:derivatives},
$g_r$ gives us the length of the arc $\partial D(s,\hat{r})\cap\Omega$
and $g_s$ tells us precisely what position this arc is along $\partial
D(s,\hat{r})$ with respect to the direction $\dgm{s}$. The assumption
of TGL for $r=\hat{r}$ implies TGL for $0<r<\hat{r}$ (see
remark~\ref{rem:2arc}) and this implies that $\gamma$ has the 2 arc
property and transverse intersections with $\partial \D(s,r)$ for all
disks corresponding to $(s,r)\in([0,L]\times\{\hat{r}\}) \cup
(\{\hat{s}\}\times [0,\hat{r}])$. Since we care only about
reconstructing a curve $\gamma$ isometric to the original curve, we
choose $\gm{\hat{s}} = (0,0)\in\Bbb{R}^2$ and $\dgm{\hat{s}} = (1,0)$.
Taken together, $g_s(\hat{s},r)$ and $g_r(\hat{s},r)$ locate both
points in $\partial D(\hat{s},r) \cap \gamma$ for all $r\in
[0,\hat{r}]$. This yields $\gamma\cap\D(\gm{\hat{s}},\hat{r})$. Now,
simply increase $s$, sliding the center of a disk of radius $\hat{r}$
along $\gamma\cap\D(\gm{\hat{s}},\hat{r})$, using $g_r(s,\hat{r})$ to
find the element of $\gamma\cap\D(\gm{s},\hat{r})$ outside
$D(\gm{\hat{s}},\hat{r})$, using the fact that the other element of
$\gamma\cap\D(\gm{s},\hat{r})$ is inside $\D(\gm{\hat{s}},\hat{r})$
and known. This process can be continued until the entire curve is
traced out in $\Bbb{R}^2$.

\section{TCGL Polygon Is Reconstructible from $g_r$ and $g_s$ without
  tail}
\label{sec:polygon-no-tail}

\begin{thm}
\label{thm:tcgl-poly-reconstruct}
For a \tcgl{} polygon $\shape$, knowing $g(s,r)$, $g_r(s,r)$ and $g_s(s,r)$ for all $s \in [0,L)$ and a particular $r$ for which $\bd$ is \tcgl{} is sufficient to completely determine $\shape$ up to translation and rotation; that is, we can recover the side lengths and angles of $\shape$.
\begin{proof}
For a given $s$ and $r$ where $g_r$ and $g_s$ exist, we can use them to obtain $r(\theta_2-\theta_1)$ as the length of the circular arc between the entry and exit points by Lemma \ref{lem:dr} and $r(\sin\theta_2-\sin\theta_1)$ as the difference in heights of the entry and exit points by Lemma \ref{lem:ds}.

We wish to recover $\theta_1$ and $\theta_2$ from these quantities.
Note that if $(\theta_1, \theta_2) = (\phi_1, \phi_2)$ is one possible solution, then so is $(\theta_1, \theta_2) = (2\pi - \phi_2, 2\pi - \phi_1)$ so solutions always come in pairs.

We can imagine placing a circular arc with angle $\frac{g_r}{r}$ on our circle and sliding it around until the endpoints have the appropriate height difference, yielding our $\theta_1$ and $\theta_2$.
Note that since $\shape$ is \tcgl{}, one endpoint must be on the left side of the circle and the other must be on the right and we cannot slide either endpoint to or beyond the vertical line through the center of the circle.

Therefore, as we slide the right endpoint down, the left endpoint slides up so that the height difference as a function of the slide is strictly monotonic.
Therefore, the slide that gives us $\theta_1$ and $\theta_2$ is unique for a given starting arc placement.
However, there are two starting arc placements: the first calls the angle for the right endpoint $\theta_1$ and the left endpoint $\theta_2$ (so the interior of $\shape$ is ``up'' in the circle) and the second swaps these (so the interior of $\shape$ is ``down'').
Since we have adopted the convention that $\bd$ is traversed in a counterclockwise direction (so the interior of $\shape$ is up in the circles) we therefore pick the first option; this gives us a unique solution for $\theta_1$ and $\theta_2$.

This procedure works whenever $g_r$ and $g_s$ exist which is certainly true whenever the density disk does not touch a vertex of $\shape$ either at its center or on its boundary because if we avoid these cases, then there is only one graph-like orientation to deal with and $\bd$ is $C^\infty$ for all the points that enter into the computation.
In fact, with a moment's thought, we can make a stronger statement than this: $g_r$ always exists and $g_s$ exists as long as the center of the density disk is not a vertex of the polygon.

We can identify the $s$ values at which $g_s(s,r)$ does not exist to obtain the arc length positions of the vertices (and therefore obtain side lengths).
For a given $s$ corresponding to a vertex, we can find $g_r$ and the one-sided derivatives $g_{s-}$ and $g_{s+}$.
These correspond to the graph-like orientations required by the polygon sides adjacent to the current vertex.

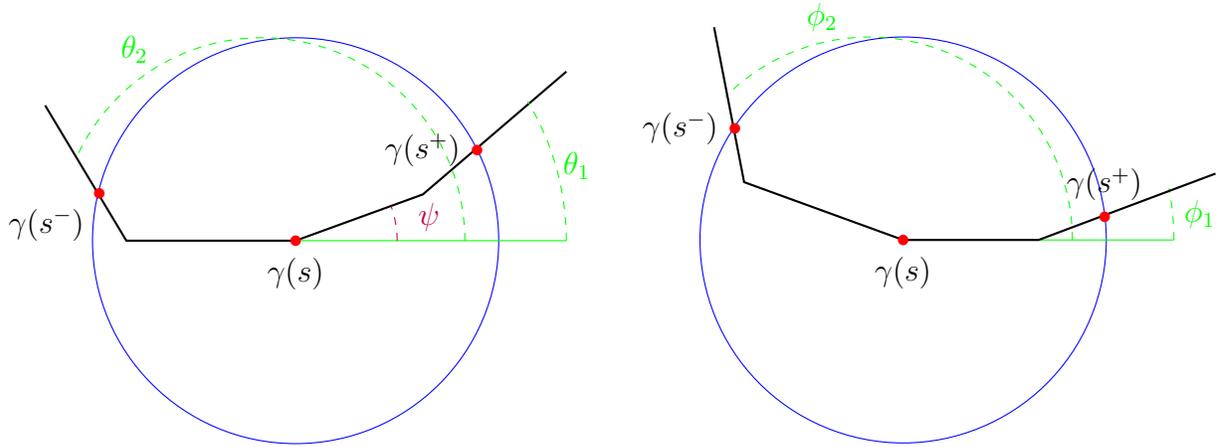
\begin{figure}
\makebox[\textwidth][c]{
\begin{tabular}{cc}
\begin{tikzpicture}[scale=0.9]
\draw[color=blue] (0,0) circle (3);

\draw[thick=2pt] (-3.7, 2) -- (-0:-2.5) -- (0,0) -- (20:2) -- (4,2.5);
\node[label=below left:{$\gamma(s^{-})$},style=circle,fill=red,minimum size=4pt,inner sep=0pt] (nenter) at (-2.9,.7) {};
\node[label=left:{$\gamma(s^{+})$},style=circle,fill=red,minimum size=4pt,inner sep=0pt] (nexit) at (2.68,1.34) {};
\node (nright) at (3,0) {};
\node (nleft) at (-3,0) {};

\draw[color=green] (0,0) -- (4,0);
\draw[dashed,color=green] (0:4) arc (0:31:4);
\node[color=green] () at (15:4.3) {$\theta_1$};
\draw[dashed,color=green] (0:2.5) arc (0:155:3);
\node[color=green] () at (130:3.7) {$\theta_2$};
\draw[dashed,color=purple] (0:1.5) arc (0:20:1.5);
\node[color=purple] () at (10:2) {$\psi$};
\node () at (270:3.23) {};
\node[label=below:{$\gamma(s)$},style=circle,fill=red,minimum size=4pt,inner sep=0pt] () at (0,0) {};

\end{tikzpicture}
&
\begin{tikzpicture}[scale=0.9,rotate=-20]
\draw[color=blue] (0,0) circle (3);
\draw[color=green] (200:0) -- (20:4);
\draw[dashed,color=green] (20:4) arc (20:31:4);
\node[color=green] () at (25:4.4) {$\phi_1$};
\draw[dashed,color=green] (20:2.5) arc (20:155:3);
\node[color=green] () at (130:3.5) {$\phi_2$};

\draw[thick=2pt] (-3.7, 2) -- (-0:-2.5) -- (0,0) -- (20:2) -- (4,2.5);
\node[label=left:{$\gamma(s^{-})$},style=circle,fill=red,minimum size=4pt,inner sep=0pt] (nenter) at (-2.9,.7) {};
\node[label=above:{$\gamma(s^{+})$},style=circle,fill=red,minimum size=4pt,inner sep=0pt] (nexit) at (2.68,1.34) {};
\node (nright) at (3,0) {};
\node (nleft) at (-3,0) {};
\node[label=below:{$\gamma(s)$},style=circle,fill=red,minimum size=4pt,inner sep=0pt] () at (0,0) {};

\end{tikzpicture}
\end{tabular}
}
\caption{Using $g_{s-}$ and $g_{s+}$ to obtain the polygon angle at $s$.}
\label{fig:tcgl-poly}
\end{figure}

Referring to Figure \ref{fig:tcgl-poly}, the one-sided derivatives along with the argument at the beginning of the proof yield the angles $\theta_1$, $\theta_2$, $\phi_1$, and $\phi_2$.
Thus we can calculate $\psi = \theta_1 - \phi_1$ which means that the polygon vertex at $s$ has angle $\pi - \psi$.

Doing this for all $s$ corresponding to vertices, we can determine all of the angles of the polygon.
With the side lengths identified earlier, this completely determines the polygon $\shape$ up to translation and rotation.
\end{proof}
\end{thm}

\section[Simple closed curves are generically reconstructible using fixed radius data]{Simple closed curves are generically\\ reconstructible using fixed radius data}
\label{sec:generic}

We will assume that $\gamma$ is TGL for the radius $\hat{r}$.  We will also
assume that we know the first, second, and third derivatives of
$g(s,r)$ for $r=\hat{r}$.  Under these assumptions, $\gamma$ is
generically reconstructible. By generic we mean the \emph{admittedly weak}
condition of density -- reconstructible curves are $C^1$ dense in
the space of $C^2$ simple closed curves.

\begin{thm}
\label{thm:tcgl-reconstruct}
  Define $\Bbb{G} \equiv \{ \gamma| \gamma$ is a $C^2$ simple closed
  curve and TGL for $r=\hat{r}\}$.  Suppose that, for $r=\hat{r}$, for
  all $s\in[0,L]$, and for each $\gamma\in\Bbb{G}$ we know the first-,
  second-, and third-order partial derivatives of $g_{\gamma}(s,r)$.
  Then the set of reconstructible $\gamma\in\Bbb{G}$ is
  $C^1$ dense in $\Bbb{G}$ where reconstructability is modulo reparametrization, translation, and rotation.
\end{thm}

\begin{figure}[htp!]
  \centering
  \input{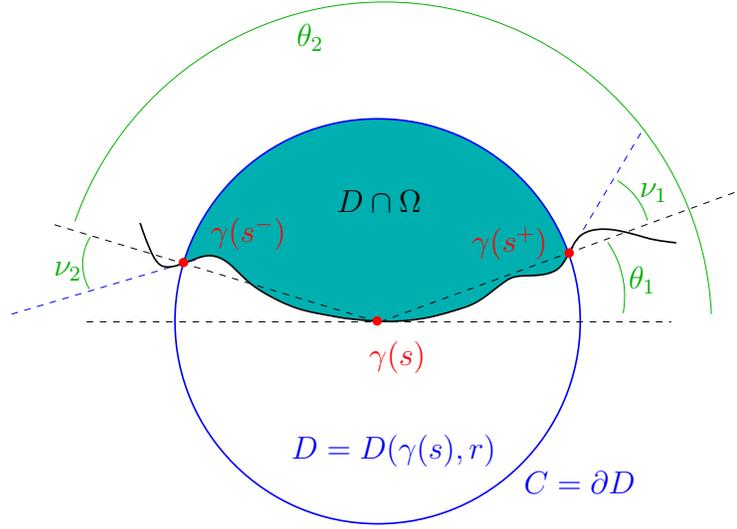}
  \caption{Figure~\ref{fig:step0} again as a reminder}
\label{fig:ref-fig}
\end{figure}

\prff In section \ref{sec:derivatives} we showed that $\frac{\partial
  g(s,r)}{\partial r} = r(\theta_2 - \theta_1)$ and $\frac{\partial
  g(s,r)}{\partial s} = r(\sin(\theta_2) - \sin(\theta_1))$, where the notation is as in Figure \ref{fig:ref-fig}.
Because $\gamma$ is TGL, we can solve for $\theta_1$ and $\theta_2$ from these two derivatives as in the proof of Theorem \ref{thm:tcgl-poly-reconstruct}.

\begin{clm}
\label{clm:2ndderiv}
The following equations hold:
  $\frac{\partial^2 g(s,r)}{\partial r^2} = \theta_2 - \theta_1 +
  r(\frac{\partial\theta_2}{\partial r} -
  \frac{\partial\theta_1}{\partial r})$ and $\frac{\partial^2
    g(s,r)}{\partial r\partial s} = \sin(\theta_2) - \sin(\theta_1) +
  r(\cos(\theta_2)\frac{\partial\theta_2}{\partial r} -
  \cos(\theta_1)\frac{\partial\theta_1}{\partial r} )$.
\end{clm}

\prfclm{clm:2ndderiv} Simply differentiate the expressions we already have for
$\frac{\partial g(s,r)}{\partial r}$ and $\frac{\partial
  g(s,r)}{\partial s}$. \epfclm 

We wish to express this in terms of $\nu_1$ and $\nu_2$.
Note that if we expand the circle radius by $\Delta r$, the right exit point $s_+(s)$ moves approximately (i.e., considering first-order terms only) a distance of $k \equiv \Delta r \sec(\nu_1)$ (so $\pd{k}{r} = \sec{\nu_1}$, a fact we will use later to compute curvature).
Therefore,
\[
\pd{\theta_1}{r} = \lim_{\Delta r \rightarrow 0} \frac{\arctan\left( \frac{r\sin\theta_1+k\sin(\theta_1+\nu_1)}{r\cos\theta_1 + k\cos(\theta_1 + \nu_1)}\right)-\theta_1}{\Delta r}.
\]
Straightforward techniques yield 
$
\pd{\theta_1}{r} %
= \frac{\tan{\nu_1}}{r}
$ and a similar calculation shows that $\pd{\theta_2}{r} = \frac{\tan\nu_2}{r}$.

Therefore, rewriting the second derivatives of $g(s,r)$ in terms of $\nu_1$ and $\nu_2$, we get:
\begin{eqnarray*}
  \frac{\partial^2 g(s,r)}{\partial r^2} & = & \theta_2 - \theta_1 +
  \tan(\nu_2) - \tan(\nu_1) \\
  \frac{\partial^2
    g(s,r)}{\partial r\partial s} & = & \sin(\theta_2) - \sin(\theta_1) +
  \cos(\theta_2)\tan(\nu_2) -
  \cos(\theta_1)\tan(\nu_1)
\end{eqnarray*}
Using these 2 derivatives, together with the previous two, we can
solve for $\nu_1 = \arctan(r \frac{\partial\theta_1}{\partial r})$ and
$\nu_2 = \arctan(r \frac{\partial\theta_1}{\partial r})$ whenever
$\cos(\theta_1) \neq \cos(\theta_2)$. Since we are assuming that the
curve is a simple closed curve, $\cos(\theta_1) \neq \cos(\theta_2)$ is
always true.

\begin{clm}
\label{clm:3rdderiv}
  Knowing $\frac{\partial^3 g(s,r)}{\partial r^3}$ and
  $\frac{\partial^3 g(s,r)}{\partial r^2\partial s}$ gives us
  $\kappa(s^+(s))$ and $\kappa(s^-(s))$, the curvatures of $\gamma$ at
  $s^+(s)$ and $s^-(s)$.
\end{clm}
\prfclm{clm:3rdderiv} Computing, we get
\begin{eqnarray*}
  \frac{\partial^3 g(s,r)}{\partial r^3} & = & \frac{\partial \theta_2}{\partial r} - \frac{\partial \theta_1}{\partial r} +
  \sec^2(\nu_2)\frac{\partial \nu_2}{\partial r} - \sec^2(\nu_1)\frac{\partial \nu_1}{\partial r} \\
  \frac{\partial^3
    g(s,r)}{\partial r^2\partial s} & = & \cos(\theta_2)\frac{\partial \theta_2}{\partial r} - \cos(\theta_1)\frac{\partial \theta_1}{\partial r} 
  - \sin(\theta_2)\frac{\partial \theta_2}{\partial r}\tan(\nu_2) \\ & + &
  \sin(\theta_1)\frac{\partial \theta_1}{\partial r}\tan(\nu_1) + \cos(\theta_2)\sec^2(\nu_2)\frac{\partial \nu_2}{\partial r} -
  \cos(\theta_1)\sec^2(\nu_1)\frac{\partial \nu_1}{\partial r}.
\end{eqnarray*} 
Since $\nu_2'\equiv\frac{\partial \nu_2}{\partial r}$ and
$\nu_1'\equiv\frac{\partial \nu_1}{\partial r}$ are the only unknowns,
we end up having to invert
\begin{equation*}
  \left[
    \begin{array}{cc}
      1 & -1 \\
     \cos(\theta_2) & \cos(\theta_1)
    \end{array}
\right]
\end{equation*} 
again and this is always nonsingular, giving us $\nu_1'$ and $\nu_2'$
as a function of s, the coordinate of the center of the disk.

Relative to the horizontal, the angle of the curve at $s^+(s)$ is $\theta_1 + \nu_1$ so the rate of change in angle as we expand the circle is $\pd{\theta_1}{r}+\nu_1'$.
Recalling that rate of movement of this exit point as we expand the circle is given by $\pd{k}{r} = \sec \nu_1$, we have that the curvature is given by $\kappa(s^+(s))= \pd{k}{r}(\pd{\theta_1}{r}+\nu_1') = \sec \nu_1(\pd{\theta_1}{r}+\nu_1')$.
Similarly, $\kappa(s^-(s)) = \sec(\nu_2)(\pd{\theta_2}{r}+\nu_2')$.
\epfclm
\begin{clm}
  Generically, we can deduce $s^+(s)$ from knowledge of $\nu_1(s)$,
  $\nu_2(s)$, $\theta_1(s)$ and $\theta_2(s)$.
\end{clm}

\prff We outline the proof without some of the explicit
constructions that follow without much trouble from the outline. We
have that $\theta_1(s^-(s)) + \nu_1(s^-(s))= \pi - \theta_2(s) -
\nu_2(s)$ and $\theta_1(s) + \nu_1(s) = \pi - \theta_2(s^+(s)) -
\nu_2(s^+(s))$. All four of these quantities (the left- and right-hand sides of each of
the 2 equations) are the turning angles between the tangent to the
curve at the center of the disk and the tangent to the curve at a
point $r$ away from the center of the disk.

Now we use this correspondence between the $\theta+\nu$ curves to
solve for $s^-(s)$ and $s^+(s)$. But these curves can differ by a
homeomorphism of the domain. Thus, we can only find the correspondence
if there is a distinguished point on those curves as well as no places
where the values attained are constant. The turning angle curves
having isolated critical points and a unique maximum or minimum is
sufficient for our purposes.

To get isolated extrema, start by approximating the curve $\gamma$
with another one, $\hat{\gamma}$, that agrees in $C^1$ at a large but
finite number of points $\{s_i\}_{i=1}^N$ (i.e. agrees in tangent
direction as well as position) and has isolated critical points in
the derivative of the tangent direction. Now perturb $\hat{\gamma}$ to
one that is $C^1$ close (but not $C^2$ close) by using oscillations
about the curve so that the 2nd and 3rd derivatives are never
simultaneously below the bounds on the 2nd and 3rd derivatives of the
curve we started with. We do this in a way that alternates around the
curve. See Figure~\ref{fig:corre_perturb}.  In a bit more detail,
suppose that $\max\{d^2\hat{\gamma}/ds^2, d^3\hat{\gamma}/ds^3\} <
L_1$. Choose a starting point on the curve; $s=0$ works. Now begin
perturbing $\hat{\gamma}$ at the point $s_{\hat{r}}$ in the positive
$s$ direction such that $|\hat{\gamma}(s_{\hat{r}}) - \hat{\gamma}(0)|
= \hat{r}$. We name the newly perturbed curve $\hat{\hat{\gamma}}$ and
we keep $L_1 < \max\{ d^2\hat{\hat{\gamma}}/ds^2,
d^3\hat{\hat{\gamma}}/ds^3 \} < L_2$.  We continue perturbing until we
have reached $s_{2r}$ defined by $|\hat{\gamma}(s_{2\hat{r}}) -
\hat{\gamma}(s_{\hat{r}})| = \hat{r}$. We begin perturbing again when
we reach $s_{3\hat{r}}$. Continue in this fashion around
$\hat{\gamma}$. The last piece, shown in green in the figure, will
require a perturbation that is distinct in size due to the fact that
it will interact with the perturbation that starts at $s_{\hat{r}}$.
On this last piece, we enforce $L_2 < \max\{
d^2\hat{\hat{\gamma}}/ds^2, d^3\hat{\hat{\gamma}}/ds^3 \} < L_3$. All
these perturbations can be chosen with isolated singularities in
derivatives, thus giving us $\theta +\nu$ curves that are monotonic
between isolated singularities.  (In fact, we might as well choose all
perturbations to be piecewise polynomial perturbations. This
immediately gives us the isolated singularities and monotonicity that
we want.)

\begin{figure}[htp!]
  \centering
  \input{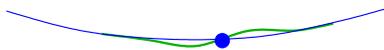}
  \caption[A twist perturbation]{A twist perturbation. Notice that if the twist is applied
    precisely at a global max of the turning angle (as measured by the
    tangent here and the one lagging it in $s$), we will increase the
    turning angle there and will end up with a unique global maximum.
  }
 \label{fig:twist_perturb}
\end{figure}

Finally, if there is not a distinct maximum, we can choose one of the
maxima and add a small twist to the curve at that point. See
Figure~\ref{fig:twist_perturb}. The idea is that a small twist,
applied to the leading edge of the tangents we are comparing to get the
turning angle, will increase the angle most at the center of the
twist. If this corresponds to a nonunique global maximum, we end up
with a unique global maximum.

\begin{figure}[htp!]
  \centering
  \input{corre_perturb.pdf_t}
  \caption[Alternating perturbation around the curve]{In this schematic figure, we illustrate the alternating
    perturbation around the curve, keeping the curve $C^1$ close to
    and messing with the second and third derivatives to eliminate any
    critical points other than isolated maxima and minima. Here the
    perturbation is of course greatly exaggerated. }
 \label{fig:corre_perturb}
\end{figure}

Now the correspondence scheme works. That is, we know that the global
maximums must match, and because the turning angle curves are
monotonic between isolated critical points, we can find the
homeomorphisms in $s$ that move the turning angle curves into
correspondence.  
\epfclm

Taken together, the last two claims give us the curvature as a
function of arclength. This determines $\gamma$ up to translations and
rotations.  \epf

\section{Numerical experiments}
\label{sec:numerics}

In this section, we consider a numerical curve reconstruction for the situation in which $g(s,r)$ is known for a given radius $r$ but no derivative information is available.  This reconstruction is more strict than the scenarios of sections \ref{sec:reconstruct-T}--\ref{sec:generic}.  Our motivation is to explore whether any $\gamma$ can be uniquely and practically reconstructed with this limited information.  

We consider $\gamma_a(\bar{s})\in \mathcal{P}^N$, the set of \textit{simple} polygons of $N$ ordered vertices $\{(x_1,y_1), \dots,\allowbreak (x_{N},y_{N})\}$ parameterized by the set $\{\bar{s}_k\}_{k=1}^{N}$ with $\bar{s}_k=k/N$ as 
\begin{equation}
\begin{array}{l}
 \displaystyle x_k = \sum_{j=0}^{m-1}{a_{1,j} \cos(2\pi j \bar{s}_k/N)+a_{2,j} \sin(2\pi j \bar{s}_k/N)},\\
 \displaystyle y_k = \sum_{j=0}^{m-1}{a_{3,j} \cos(2\pi j \bar{s}_k/N)+a_{4,j} \sin(2\pi j \bar{s}_k/N)},
\end{array}
\end{equation}
for some coefficients $a_{i,j}\in\mathbb{R}$.  In this way, the polygon $\gamma$ is a discrete approximation of a $C^{\infty}$ curve.  The sides of $\gamma_a(\bar{s})$ are not necessarily of equal length.

We take the vector signature $g_a(\bar{s},r)\in\mathbb{R}^N$ to be the discrete area densities of $\gamma_a(\bar{s})$ computed at each vertex.  Given such a signature for fixed radius $r$ and fixed partition $\bar{s}$, we seek $a^*$ satisfying
\begin{equation}
\label{equ:numopt}
\begin{array}{rl}
 \displaystyle  a^* \in &  \displaystyle \arg\min_{b\in\mathbb{R}^{4m}}\|g_b(\bar{s},r)-g_a(\bar{s},r)\|_2^2 \\
  & \\
  &  \displaystyle \text{s.t.  } \gamma_b\in\mathcal{P}^N
\end{array}
\end{equation}
Equation (\ref{equ:numopt}) represents a nonlinearly constrained optimization problem with continuous nonsmooth objective.  The constraint ensures that polygons are simple though any optimal reconstruction $\gamma_{a^*}$ is not expected to lie on the feasible region boundary except in cases of noisy signatures.  This approach to reconstructing curves seeks a polygon that matches a given discrete signature, rather than an analytic sequential point construction procedure.

We use the direct search \textsc{OrthoMads} algorithm\cite{AbAuDeLe2009} to solve this problem.  \textsc{Mads} class algorithms do not require objective derivative information\cite{AbAuDeLe2009,AuDe2006} and converge to second-order stationary points under reasonable conditions on nonsmooth functions\cite{AbAu2006}.  We implement our constraint using the extreme barrier method\cite{AuDeLe2009} in which the objective value is set to infinity whenever constraints are not satisfied. 
We utilize the standard implementation with partial polling and minimal spanning sets of $4m+1$ directions.

\begin{figure}[h!]
\centering
\includegraphics[width=.6\textwidth]{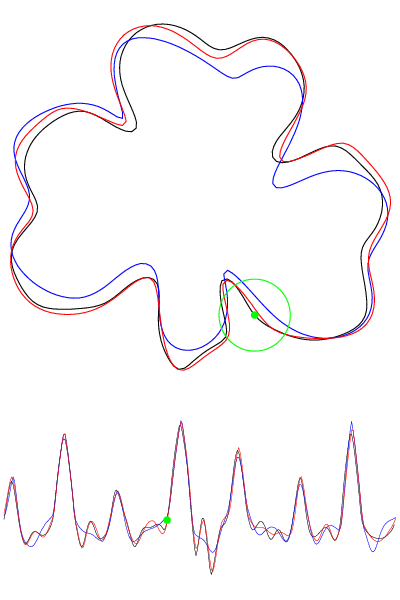}
\caption[Shamrock reconstruction]{Shamrock reconstruction: comparing the original curve with those found for {\color{blue}$m = 12$} and {\color{red}$m=18$}.  Curves for $m \geq 20$ are visually indistinguishable from the original curve.  The shape signatures are given at the bottom.}
\label{fig:shamrock}
\end{figure}

We performed a series of numerical tests using the synthetic shamrock curve shown in black in the upper portion of Figure~\ref{fig:shamrock}.  This curve is given as a polygon in $\mathcal{P}^{256}$ with discretization coefficients $a\in\mathbb{R}^{4\times 20}$ ($m=20$).  A sequence of reconstructions was performed with all integer values $8\leq m \leq 20$.  The $m=8$ reconstruction begins with initial coefficients, $a_{i,j}$, which determine a regular $256$-gon with approximately the same interior area as the shamrock (as determined by the signature $g_a(\bar{s},r)$.
In particular, the value(s) $a_{i,j}$ supplied initially are those which define the best fit circle ($m=1$),
which can be computed directly.
That is, only $a_{1,0}$ and $a_{4,0}$ are nonzero.
Subsequent reconstructions begin with initial coefficients optimal to the previous relatively coarse reconstruction.  Curve reconstructions for $m=12$ (blue) and $m=18$ (red) are compared to the shamrock in the upper portion of Figure~\ref{fig:shamrock}.   Reconstructions for $m\geq20$ are visually indistinguishable from the actual curve and are not shown.  Corresponding area density signatures are shown in the lower portion of Figure~\ref{fig:shamrock}.  A representative disk of radius $r$ is shown in green along with corresponding location in the signature; note that the shamrock is not \tcgl{} with this radius.

When comparing and interpreting the shamrock curves, it is important to note that the scale of the curves is determined entirely by the fit parameters $a_{i,j}$.  On the other hand, as the density signature is independent of curve rotation, the rotation is eyeball adjusted for easy visual comparison.  
Also note that the two-arc property does not hold for this example so our reconstructability results do not apply.
The accuracies of both the curve reconstruction and area density signature fit suggest that somewhat more general reconstructability results hold.
In particular, we speculate that general simple polygons may be reconstructible from $g(s,r)$ for fixed $r$ and no derivative information.

%\bibliographystyle{plain}
%\bibliography{TomSectionBib.bib}

\section{Conclusions}
We have studied the integral area invariant with particular emphasis on the \tcgl{} condition.
In particular, we have shown that all TCGL polygons and a $C^1$-dense set of $C^2$ TGL curves are reconstructible using only the integral area invariant for a fixed radius along the boundary and its derivatives.

We also showed that TCGL boundaries can be approximated by TCGL polygons, determined what the derivatives represented, and commented on other sets of data sufficient for reconstruction (namely, both T-like and all radii in a neighborhood of 0).

These reconstructions are all modulo translations, rotations, and reparametrizations.
The arc length parameterization plays a special role here since any two such parameterizations of a boundary will differ only by a shift and can easily be placed into correspondence.
The situation becomes more complicated in higher dimensions as boundaries are no longer canonically parameterized by a single variable which is a fundamental assumption of our results and methods.
It is not immediately obvious how to resolve the issues created by higher dimensions except that it may be possible to modify some of the machinery to work with star convex regions which restore some semblance of canonical representation.

Another space which is open for further development is that of reconstruction algorithms.
This is doubly true since our theoretical reconstructions are unstable and the numerical examples in the present work do not have guaranteed reconstruction.
However, even without these guarantees, the numerical examples hint at more expansive reconstructability results.

\section{Acknowledgments}
The authors would like to thank David Caraballo for introducing us to this topic as well as Simon Morgan and William Meyerson for initial discussions and work on related topics that are not in this paper.
This research was supported in part by National Science Foundation grant DMS-0914809.

\section{Appendix: Easy Reconstructability}
\label{sec:appen}

For completeness, we include a short proof of the fact that knowing
$g(s,r)$ for all $s$ and $r$ very easily gives us reconstructability.
This follows from the fact that knowing the asymptotic behavior of
$g(s,r)$ as $r\rightarrow 0$ for any $s$ gives us $\kappa(s)$. That in
turn implies that knowing $g(s,r)$ in any neighborhood of the set
$(s,r)\in [0,L] \times \{r=0\}$ also gives us $\kappa(s)$ and
therefore the curve.

\begin{figure}[htp!]
\centering
\subfigure[]{
\label{fig:pos-osculating-circle}
\begin{tikzpicture}[scale=0.7]
\fill[color=teal!40] (0,3) circle (3);
\draw[color=black] (0,3) circle (3);
\draw[color=blue] (0,0) circle (2);
\draw[ultra thick] plot[smooth] coordinates {(195:3) (190:2.5) (185:2) (180:1.5) (175:1) (175:0.5) (0,0) (5:0.5) (10:1) (15:1.8) (20:2.5) (25:3) (30:3.5) (35:4)};
\node () at (25:4) {$\bd$};
\draw[<->] (0:0) -- (90:3) node[pos=0.5,left] {$R$};
\draw[color=blue,<->] (0:0) -- (45:2) node[pos=0.5,left] {$r$};
\fill[color=red] (0,0) circle (2pt);
\node[color=red] () at (0,-.5) {$\gamma(s)$};
\end{tikzpicture}}
\subfigure[]{
\label{fig:neg-osculating-circle}
\begin{tikzpicture}[scale=0.7,rotate=180]
\fill[color=teal!40] (0,0) circle (2);
\fill[color=white] (0,3) circle (3);
\draw[color=black] (0,3) circle (3);
\draw[color=blue] (0,0) circle (2);
\draw[ultra thick] plot[smooth] coordinates {(195:3) (190:2.5) (185:2) (180:1.5) (175:1) (175:0.5) (0,0) (5:0.5) (10:1) (15:1.8) (20:2.5) (25:3) (30:3.5) (35:4)};
\node () at (25:4) {$\bd$};
\draw[<->] (0:0) -- (90:3) node[pos=0.5,left] {$R$};
\draw[color=blue,<->] (0:0) -- (45:2) node[pos=0.5,left] {$r$};
\fill[color=red] (0,0) circle (2pt);
\node[color=red] () at (0,-.5) {$\gamma(s)$};
\end{tikzpicture}
}
\caption[Using the osculating circle as a surrogate for $\bd$]{Using the osculating circle as a surrogate for $\bd$ in the \subref{fig:pos-osculating-circle} positive and \subref{fig:neg-osculating-circle} negative curvature cases.}
\label{fig:osculating-circle}
\end{figure}
\begin{thm}
Suppose $\bd$ is $C^2$ and there exists $\epsilon > 0$ such that we know $g(s,r)$ for all $(s,r) \in [0,L) \times (0,\epsilon)$.
This information is enough to determine the curvature of every point on $\bd$.
In particular, if $\gamma : [0,L) \rightarrow \bd$ is a counterclockwise arclength parameterization of $\bd$, then $\kappa(\gamma(s)) = -3\pi \lim_{r\rightarrow 0} \pd{}{r} \frac{g(s,r)}{\pi r^2}$.
\begin{proof}
Fix $s \in [0,L)$.
If the curvature of $\gamma$ at $s$ is positive, we consider what happens if we replace $\shape$ with the disk whose boundary is the osculating circle of $\bd$ at $\gamma(s)$ (call its radius $R$).
We have the following expression for the new normalized nonasymptotic density (see Figure \ref{fig:pos-osculating-circle}):
\[
\frac{g(s,r)}{\pi r^2} = \frac{1}{\pi r^2}\int_{-p}^{p} \sqrt{r^2-x^2}-(R-\sqrt{R^2-x^2})\, dx.
\]
where $x = p$ is the positive solution to $\sqrt{r^2-x^2} = R-\sqrt{R^2-x^2}$.
Differentiating with respect to $r$ and then taking the limit as $r$ goes to 0 gives us $-\frac{1}{3\pi R}$.
That is, for the case where $\shape$ is locally a disk, the curvature at $\gamma(s)$ is given by $-3\pi\lim_{r\rightarrow 0}\pd{}{r}\frac{g(s,r)}{\pi r^2}$.

If the curvature of $\bd$ at $\gamma(s)$ is negative, we can set up a similar surrogate (see figure \ref{fig:neg-osculating-circle}) and again obtain that $\kappa(\gamma(s)) = -3\pi \lim_{r\rightarrow 0} \pd{}{r} \frac{g(s,r)}{\pi r^2}$.

Lastly, this calculation gives the right result in the curvature 0 case when $\bd$ is locally a straight line (so $\frac{g(s,r)}{\pi r^2} = \frac{1}{\pi r^2}\int_{-r}^{r}\sqrt{r^2-x^2}\,dx = \frac{1}{2}$ for sufficiently small $r$ and $-3\pi\lim_{r\rightarrow 0} \pd{}{r} \frac{g(s,r)}{\pi r^2} = 0$).

For the case where $\bd$ is not locally a circle or straight line, the corrections to the integrals are of order $O(x^3)$ as $r$ goes to 0 and have no impact on the final answer so the curvature at $\gamma(s)$ is always given by $-3\pi\lim_{r\rightarrow 0}\pd{}{r}\frac{g(s,r)}{\pi r^2}$.
The available data (the values $g(s,r)$ for all $s \in [0,L)$ and all $r \in (0, \epsilon)$) are sufficient to compute the relevant derivative and limit so we can use this process to determine the curvature of every point on the $C^2$ curve $\bd$.
\end{proof}
\end{thm}

%\appendix
\bibliography{dissertation}
\bibliographystyle{plain}
\end{document}